\documentclass[10pt, twocolumn]{article}
\usepackage{cite}
\usepackage{amsmath,amssymb,amsfonts}
\usepackage{algorithmic}
\usepackage{graphicx}
\usepackage{textcomp}
\usepackage{lipsum,multicol}
\usepackage{float}
\usepackage{color,soul}
\usepackage{lipsum,multicol}
\usepackage{float}
\usepackage{hyperref}

\usepackage[
left=2.00cm,
right=2.00cm,
top=2.00cm,
bottom=2.00cm
]{geometry}

\begin{document}

\title{Low--Order Nonlinear Animal Model of Glucose Dynamics for a Bihormonal Intraperitoneal Artificial Pancreas\footnote{This article has been accepted for publication in IEEE Transactions on Biomedical Engineering, \url{http://doi.org/10.1109/TBME.2021.3125839}. This work was carried out during the tenure of an ERCIM `Alain Bensoussan' Fellowship Programme at the Norwegian University of Science and Technology (NTNU). This research is also funded by The Norwegian Research Council (Project no.: 248872/O70), Central Norway Regional Health Authority and Johan Selmer Kvanes Endowment for Research and Combating of Diabetes.
}}

\author{Claudia Lopez-Zazueta\thanks{Claudia Lopez-Zazueta was with the Department of Engineering Cybernetics, Faculty of Information Technology and Electrical Engineering, Norwegian University of Science and Technology  (NTNU), O. S. Bragstads Plass 2D, 7034 Trondheim, Norway. She is now with MaIAGE INRAE, Université Paris-Saclay, Domaine de Vilvert, 78352 Jouy-en-Josas, France (e-mail: claudia.lopez-zazueta@inrae.fr)
}, Øyvind Stavdahl and Anders Lyngvi Fougner \thanks{Øyvind Stavdahl and Anders Lyngvi Fougner are with the Department of Engineering Cybernetics, Faculty of Information Technology and Electrical Engineering, Norwegian University of Science and Technology  (NTNU), O. S. Bragstads Plass 2D, 7034 Trondheim, Norway (e-mail: oyvind.stavdahl@ntnu.no, anders.fougner@ntnu.no)
}}

\maketitle
\date{}

\begin{abstract}
Objective: The design of an Artificial Pancreas (AP) to regulate blood glucose levels requires reliable control methods. Model Predictive Control has emerged as a promising approach for glycemia control. However, model--based control methods require computationally simple and identifiable mathematical models that represent glucose dynamics accurately, which is challenging due to the complexity of glucose homeostasis. Methods: In this work, a simple model is deduced to estimate blood glucose concentration in subjects with Type 1 Diabetes Mellitus (T1DM). Novel features in the model are power--law kinetics for intraperitoneal insulin absorption and a separate glucagon sensitivity state. Profile likelihood and a method based on singular value decomposition of the sensitivity matrix are carried out to assess parameter identifiability and guide a model reduction for improving the identification of parameters. Results: A reduced model with 10 parameters is obtained and calibrated, showing good fit to experimental data from pigs where insulin and glucagon boluses were delivered in the intraperitoneal cavity. Conclusion: A simple model with power--law kinetics can accurately represent glucose dynamics submitted to intraperitoneal insulin and glucagon injections. The reduced model was found to exhibit local practical as well as structural identifiability. Importance: The proposed model facilitates intraperitoneal bi-hormonal model-based closed-loop control in animal trials.
\end{abstract}

Keywords: Artificial pancreas (AP), power--law kinetics, model validation, parameter identification.

\section{Introduction}\label{sec:introduction}

Type 1 Diabetes Mellitus (T1DM) is the condition resulting from a deficiency in the production of insulin from $\beta$--pancreatic cells. A common therapy to control this disease consists of exogenous insulin infusions given several times per day. Insulin injections have to be controlled to some extent manually, which is a laborious task and a major concern for people with T1DM and their families \cite{christiansen2017review, cinar2018artificial}. 

In order to help patients with T1DM, the idea of a fully--automated Artificial Pancreas (AP) has been studied for decades \cite{goodwin2019performance}. Basically, an AP is a device that uses blood glucose measurements (collected with a sensor from a subject) in a decision--making algorithm that estimates the necessary amount of insulin to be administered in the subject \cite{cinar2018artificial}. The quantity of insulin to be infused must be precisely calculated in order to keep glucose levels in a safe and optimal range \cite{christiansen2017review}. Additionally, glucagon infusions can be used to counteract severe hypoglycemia \cite{lv2013pharmacokinetics, dirnena2018intraperitoneal}.

AP can be classified in hybrid or fully--automated AP systems. Hybrid systems demand the user to announce known disturbances (e.g. meals, physical activity, etc.), while for fully--automated AP the subject does not need to take part in the control \cite{beneyto2018new}. A fully--automated AP is expected to provide a better control and reduce hyperglycemia and hypoglycemia occurrences \cite{magdelaine2015long, messori2018individualized, bianchi2019invalidation}.  

However, the slow insulin absorption is one obstacle to achieve stable and safe control of blood glucose levels. For instance, in the case of subcutaneous insulin infusions there is a significant delay of insulin transport through the subcutaneous tissue. Then, to keep glucose levels within a target range, the patient is required to announce in advance events such as the ingestion of meals with an estimation of carbohydrate contents or physical activities \cite{christiansen2017review, cinar2018artificial}. Nevertheless, this relies on the ability of the patient to announce events, which is cumbersome to repeat multiple times per day and patients struggle to properly estimate meal content.

Consequently, a fully--automated AP would be of great help to any patient and, if its operation is accurate, it could improve the control of blood glucose levels compared to manual glucose regulation. However, there are still problems to face and efforts must be made to achieve a fully--automated AP.

The control of blood glucose concentration in T1DM is a challenging problem because of the complex, non--linear, and time--varying dynamics of glucose homeostasis  \cite{cinar2018artificial}, including delays in insulin infusions \cite{christiansen2017review}, inter and intra--subject variability \cite{heise2004lower, toffanin2019multiple}, etc. Model Predictive Control (MPC) has emerged as a promising approach for the control of blood glucose levels \cite{cobelli2009diabetes}. Effective MPC requires an accurate and individualized model of patient glucose--insulin dynamics  \cite{garcia2019closed}. However, the identification of such models is a tough task and the methods to measure individual parameters are invasive and expensive \cite{messori2018individualized}.

Several proposed models attempt to describe the glucose metabolism in some detail \cite{cobelli1982integrated, hovorka2004nonlinear, man2014uva}, and are thus useful for simulating glucose dynamics. However, these nonlinear, high--order models have a large number of parameters that in general cannot be identified from easily obtainable data \cite{cobelli2009diabetes}. Moreover, their usage in MPC is computationally demanding, making their integration on an AP impracticable  \cite{messori2018individualized} and does not assure a better closed loop control  \cite{bianchi2019invalidation}.

In contrast, it has been shown that low--order linear models with parameters estimated from clinical data can be suitable for MPC with subcutaneous insulin delivery  \cite{percival2011development}. However, more computationally tractable, minimalistic, and rather linear models \cite{bolie1961coefficients, yipintsoi1973mathematical, bergman1979quantitative, candas1994adaptive} often do not represent well the non--linear dynamics that have been seen after intraperitoneal insulin infusions \cite{lopez2019simple}. 

The purpose of this article is to present a relatively simple nonlinear model to approximate T1DM glucose dynamics when insulin and glucagon boluses are introduced in the intraperitoneal (IP) cavity. The conception of a new model is also motivated by the fact that most of the work developed for AP systems is adapted to subcutaneous insulin delivery \cite{chakrabarty2019new}. 

A tight glycemic control is difficult to attain with an AP operating with subcutaneous insulin infusion, due to delays in insulin absorption and slow insulin--clearance rates \cite{christiansen2017review, chakrabarty2019new}. For these reasons, the intraperitoneal route to infuse insulin has been investigated, since insulin--glucose kinetics are significantly faster in the IP cavity than subcutaneously \cite{fougner2016intraperitoneal, chakrabarty2019new}. Furthermore, it has been observed that blood glucose increases faster when glucagon is delivered in the IP cavity, compared with subcutaneous glucagon infusions \cite{dirnena2018intraperitoneal}. 

 Therefore, the problem of not having a stable and safe control due to the slow subcutaneous hormone (insulin and glucagon) absorption can be overcome with an intraperitoneal AP, because the delay of hormone absorption would be considerably reduced.

In addition, a small identifiable model to predict blood glucose dynamics during intraperitoneal insulin and glucagon infusions is required for future experimental purposes. The objective is to calibrate the model with subcutaneous glucose measurements and / or intravenous (IV) blood samples analyzed on a blood gas machine \cite{dirnena2019intraperitoneal, dirnena2018intraperitoneal, am2020intraperitoneal} from the first 2--3 hours of an experiment, to  obtain a personalized model and then decide a model--based controller to normalize blood glucose during the rest of the animal trial using intraperitoneal boluses. This implies that the parameters of the model have to be estimated as fast as possible, in order to test the control method throughout the rest of experiment (which, with anesthetized pigs, may last for not more than 10--12 hours). 

This basic procedure will mimic the operation of an AP, which requires a model feasibly adaptable (i.e. a model that can be personalized) to the glucose dynamic scenario of a determined time frame. Additionally, having a simple model that can be adapted multiple times a day will allow for intra--subject variability to be addressed.

In order to elucidate whether and how the proposed model can be simplified, its  local practical and structural parameter identifiability has been analyzed. Parameter profile likelihoods as well as a method based on Singular Value Decomposition (SVD) of the sensitivity matrix have been used. The analysis has led to detect non--identifiable parameters and to define a reduced model with less parameters to estimate, while keeping a satisfactory accuracy when approximating blood glucose measurements.

The paper is organized as follows:  In Section \ref{Chapter:complete_model}, a model to describe T1DM glucose dynamics is introduced. The model accounts for intraperitoneal infusion of insulin and glucagon, as well as for an exogenous IV glucose input. 
In Section \ref{sec:reduced_model}, a reduced model to describe T1DM glucose dynamics is presented. The parameter identification analysis of the reduced model is also exposed. In Section \ref{section:data}, the experiments from where glucose data were obtained are described, as well as the method to calibrate the model with the data. In Section \ref{Section:reduced_simulations}, the results obtained after calibrating the reduced model are shown. Finally, the discussion and conclusion about the work are exposed in Sections  \ref{Section:discussion} and \ref{Section:conclusion}, respectively.

\section{Bihormonal--Glucose Model}\label{Chapter:complete_model}
\begin{table*}[!]\centering\scriptsize
\caption{Description of the states, inputs, and parameters of the bihormonal--glucose model \eqref{glucose}--\eqref{glycogen} and the reduced bihormonal--glucose model \eqref{red:glucose}. Blood glucose concentration is the only output of the system and the rest of states are considered dimensionless.
Abbreviations: CS: circulatory system; IC: intermediate compartment; IP: Intraperitoneal Cavity; IV: intravenously.}\label{table:bimodel}
\begin{tabular}{cccc}
\hline
	\bf{Description}				&	\bf{State}		&	\bf{Units}			&	\bf{Compartment}\\
\hline
Blood Glucose	&	$G$			&	mmol/L			&	CS\\
Blood Insulin 				&	$I$			&	dimensionless				&	CS \\
Insulin IC					&	$i_1$		&	dimensionless			&	IC \\
Insulin IP			 	 	&	$i_2$		&	dimensionless			&	IP\\
Blood Glucagon 				&	$H$			&	dimensionless				&	CS\\
Glucagon IP 					&	$h_1$		&	dimensionless			&	IP \\
Glucagon sensitivity			&	$\xi$		&	dimensionless			&	IC\\
\hline
							&	\bf{Input}	&					&	\\
\hline
Exogenous Glucose 	&	$Ra_G$	&	mmol/h 			&	CS\\	
infusion IV&	&	&\\
Insulin IP bolus				&	$u_I$		&	U 				&	IP\\
Glucagon IP bolus			&	$u_H$	&	$\mu$g				&	IP\\
\hline
										&	\bf{Parameter}						&				&	\\
\hline
Blood Insulin					&	$I_b$							&	dimensionless				&	CS \\
basal value	&	&	&\\
Blood Glucagon			 	&	$H_b$							&	dimensionless				&	CS\\
basal value	&	&	&\\
Insulin--independent	&	$k_1$							&	1/d				&	CS	\\
removal rate of glucose	&								&					&		\\
Insulin--dependent	&	$k_I, k_{i_1}$					&	1/d				&	CS	\\
removal rates of glucose	&&	&		\\
Exogenous glucose 	&	$r_G$							&	h/(L$\cdot$ d)			&	CS	\\
rate of appearance		&	&	&\\
Glucose response 			&	$k_H$							&	1/d				&	CS	\\
to glucagon rate	&	&	&\\
Consumption, degradation, 				&	$m_1, m_ 2$			&	1/d				&	CS, IC, IP \\
and transport rates						&	$m_3, m_4,$					&					&	 \\
	&	$ n, n_1, n_2$					&					&	 \\

Decrease rate of		&	$x_1$							&	1/d 				&	IC\\
glucagon sensitivity		&	&	&\\

Restoration rate		&	$x_2$							&	L/(mmol$\cdot$ d)		&	IC\\
glucagon sensitivity		&	&&\\

Powers 									&	$p, q$							&	dimensionless				&	CS, IC\\
\hline
\end{tabular}
\end{table*}

In this Section, a model to simulate T1DM glucose dynamics  is described. This model accounts for insulin and glucagon boluses introduced in the IP cavity. The model is based on two different models, which were deduced to approximate experimental data from experiments with either intraperitoneal insulin boluses, or intraperitoneal and subcutaneous glucagon boluses. These two initial models are presented in Supplementary material S.1 and Supplementary material S.2

The bihormonal--glucose model is the following:
\begin{align}
\label{glucose}
\frac{dG}{dt}  = & - [k_1 + k_I\cdot ( I + I_b ) + k_{i_1}\cdot i_1]\cdot G  \\\notag
& 	+ k_H\cdot ( H +H_b )  \cdot \xi+ r_G\cdot Ra_G	\\
\label{insulin}
\frac{dI}{dt} = & - m_1\cdot I + m_2\cdot i_1^p 	\\
\label{insulinIC}
\frac{di_1}{dt} = & { - m_3\cdot i_1^q}+ m_4\cdot i_2 	\\	
\label{insulinIP}
\frac{di_2}{dt} = & - m_4\cdot  i_2 + u_I	\\	
\label{glucagon}
\frac{dH}{dt} = & - n\cdot H + {n_2 \cdot h_1} 	\\
\label{glucagonIP}
\frac{dh_1}{dt} = & - n_1\cdot h_1 + u_H	\\
\label{glycogen}
\frac{d\xi}{dt} = & - x_1\cdot H\cdot \xi + {x_2\cdot G\cdot I}
\end{align}
The states, inputs, and parameters of the bihormonal--glucose model \eqref{glucose}--\eqref{glycogen} are described in Table \ref{table:bimodel}.

In the bihormonal--glucose model \eqref{glucose}--\eqref{glycogen}, glucose consumption can be insulin-independent with rate $k_1\cdot G$, or insulin-dependent relying on insulin states $I$ and $i_1$. On the other hand, the increase in blood glucose concentration is regulated by glucagon levels $H$ and the exogenous IV glucose input $Ra_G$.

The several insulin and glucagon states account for the transport of each hormone through different compartments, from the IP cavity (where the hormone boluses are released) to the circulatory system. The infusion of exogenous insulin and glucagon is represented by the external inputs $u_I$ and $u_H$, respectively.

In the model it is assumed that the transport of insulin between compartments is nonlinear. This assumption was made based on experimental data \cite{dirnena2019intraperitoneal}. For more details see \cite{lopez2019simple} and Section Supplementary material S.1. This was described in \eqref{insulin} and \eqref{insulinIC} with power--law kinetics, including the exponents $p$ and $q$. 

Power--law approximation or synergistic systems (S--systems) are used to describe in a non canonical form reactions with particular non--linearities \cite{savageau1969biochemical, voit2015150, crampin2004mathematical}. The power--law formulation accounts for the change rate of a state as the difference of two products of states raised to non--integer powers \cite{savageau1969biochemicalII, savageau1970biochemical, crampin2004mathematical}. 
In this work, given the uncertainty of insulin concentration in compartments where measurements cannot be non--invasively obtained, power--law was used to simulate insulin dynamics using few parameters and simple equations. 

According to experimental data where glucagon boluses were administered in the IP cavity and subcutaneously \cite{am2020intraperitoneal}, the changes in blood glucose concentration induced by glucagon boluses are not always linearly proportional to bolus sizes (see Section Supplementary material S.2). For this reason, the state $\xi$ to represent the sensitivity to glucagon boluses is included in \eqref{glycogen}. Since glucagon sensitivity might be linked to the amount of glycogen store in the liver or to the hepatic responsivity to glucagon \cite{visentin2016one}, it is assumed that glucagon sensitivity decrease is proportional to glucagon concentration and it can be restored when glucose is stored in the liver in presence of insulin. For more details about the bihormonal--glucose model \eqref{glucose}--\eqref{glycogen} see Supplementary material S.3.

\section{Reduced Bihormonal--Glucose Model }\label{sec:reduced_model}

The reduced bihormonal--glucose model \eqref{red:glucose} is
\begin{align}\label{red:glucose}
\frac{dG}{dt}  = & - [k_1 + \overline{k_{i_1}}\cdot \overline{i_1}]\cdot G \\\notag 
&	+ \overline{k_H}\cdot (\overline H + \overline{H_b}) \cdot \overline\xi + r_G\cdot Ra_G \\\notag
\frac{d\overline{i_1}}{dt} = & { - \overline{m_3}\cdot \overline{i_1}^q}+ i_2 	\\\notag	
\frac{di_2}{dt} = & - m_4\cdot  i_2 + u_I	\\\notag	
\frac{d\overline{H}}{dt} = & - {n}\cdot \overline{H} + { h_1} 	\\\notag
\frac{dh_1}{dt} = & - n_1\cdot h_1 + u_H	\\\notag
\frac{d\overline\xi}{dt} = & - \overline x_1\cdot \overline{H}\cdot \overline\xi + { G\cdot \overline{i_1}}
\end{align}

The reduced bihormonal--glucose model \eqref{red:glucose} is a reduced version of the bihormonal--glucose model \eqref{glucose}--\eqref{glycogen} presented in \ref{Chapter:complete_model}. It is obtained following the transformations to address the lack of local structural and practical identifiability, which are described in Supplementary material S.8 and Supplementary material S.9. The reduced model accounts for 6 states and 10 parameters (1 state and 6 parameters less than the complete model).

The reduced bihormonal--glucose model \eqref{red:glucose} is validated using experimental data from pigs where insulin and glucagon boluses were administered in the IP cavity. The details are presented in the next two sections. Furthermore, the local practical and structural parameter identifiability of the reduced bihormonal--glucose model \eqref{red:glucose} were analyzed as explained in Section \ref{sec:PLreduced} and  Section \ref{sec:SVDreduced}.

\subsection{Profile Likelihood of the reduced bihormonal--glucose model}\label{sec:PLreduced}

\begin{figure}[!]\centering
\includegraphics[height=21cm]
{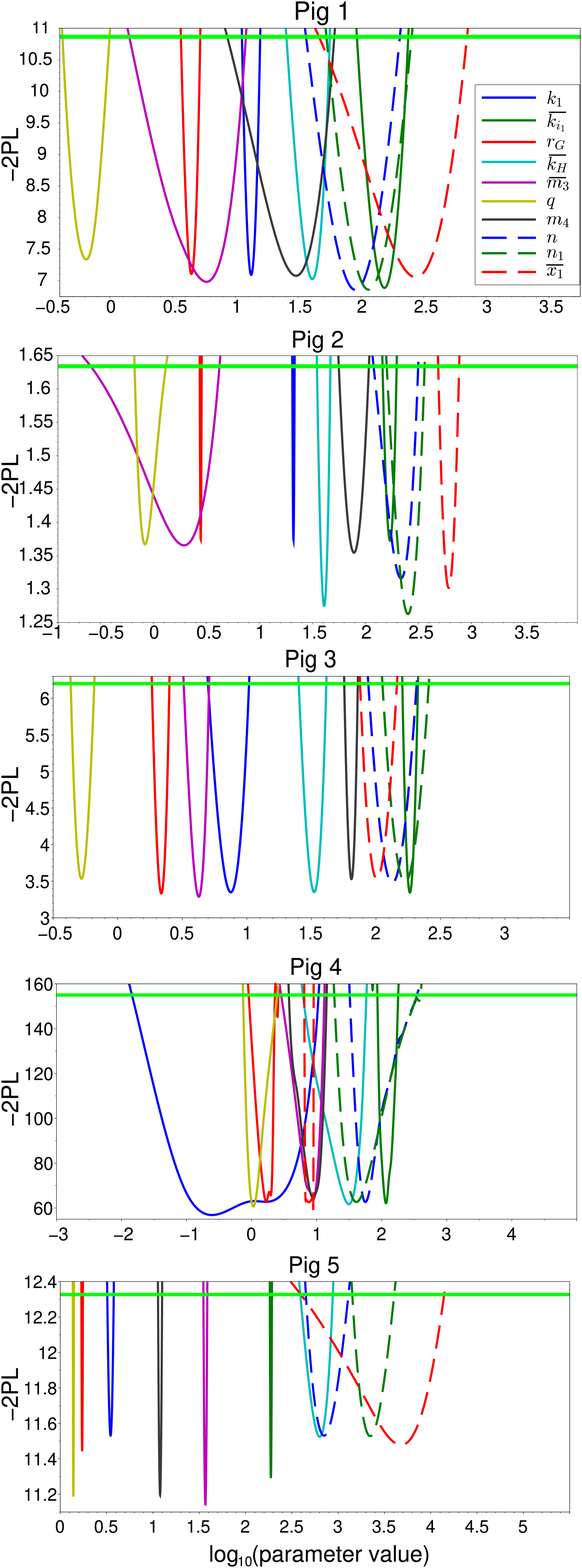}
\caption{Profile likelihoods of the reduced bihormonal--glucose model \eqref{red:glucose}. The parameters have profile likelihoods with single minima and that exceed the confidence threshold twice, which suggests that they are locally practically identifiable. The green horizontal lines indicate the confidence thresholds of 99\% for Pig 1, 98\% for Pigs 2,3, and 4, and 96\% for Pig 5.}\label{fig:PLreduced}
\end{figure}

Profile likelihoods \cite{kreutz2013profile} were computed for the reduced bihormonal--glucose model \eqref{red:glucose} using the method described in Supplementary material S.5. The results are depicted in Fig. \ref{fig:PLreduced}. All profile likelihoods have a single minima and exceed the confidence threshold twice. Therefore, no parameters of the reduced bihormonal--glucose model \eqref{red:glucose} with lack of local practical identifiability were found.

\subsection{Singular Value Decomposition of the reduced bihormonal--glucose model}\label{sec:SVDreduced}

\begin{figure}[!]\centering
\includegraphics[width=4.25cm]{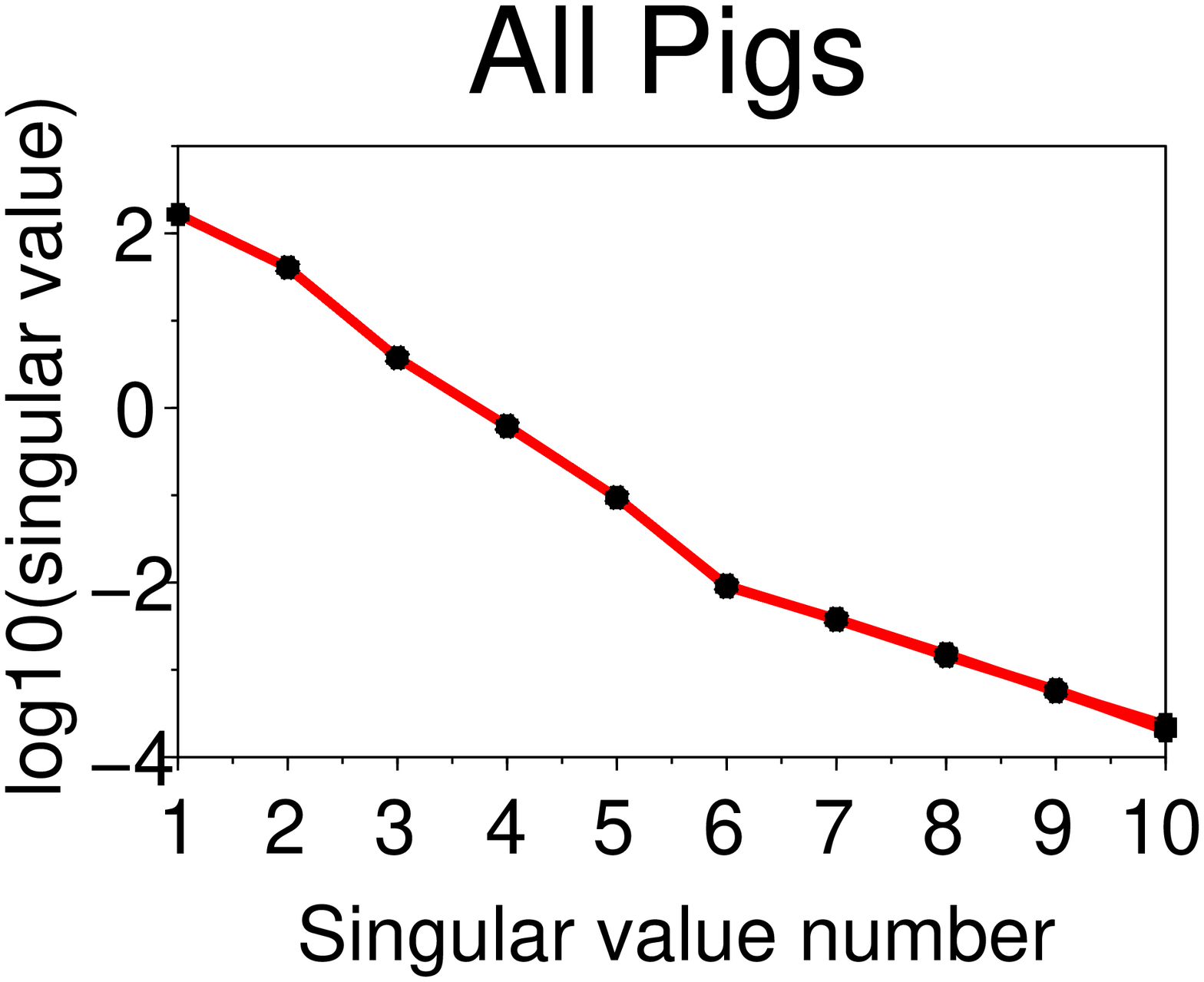}\\
\includegraphics[width=4.25cm]{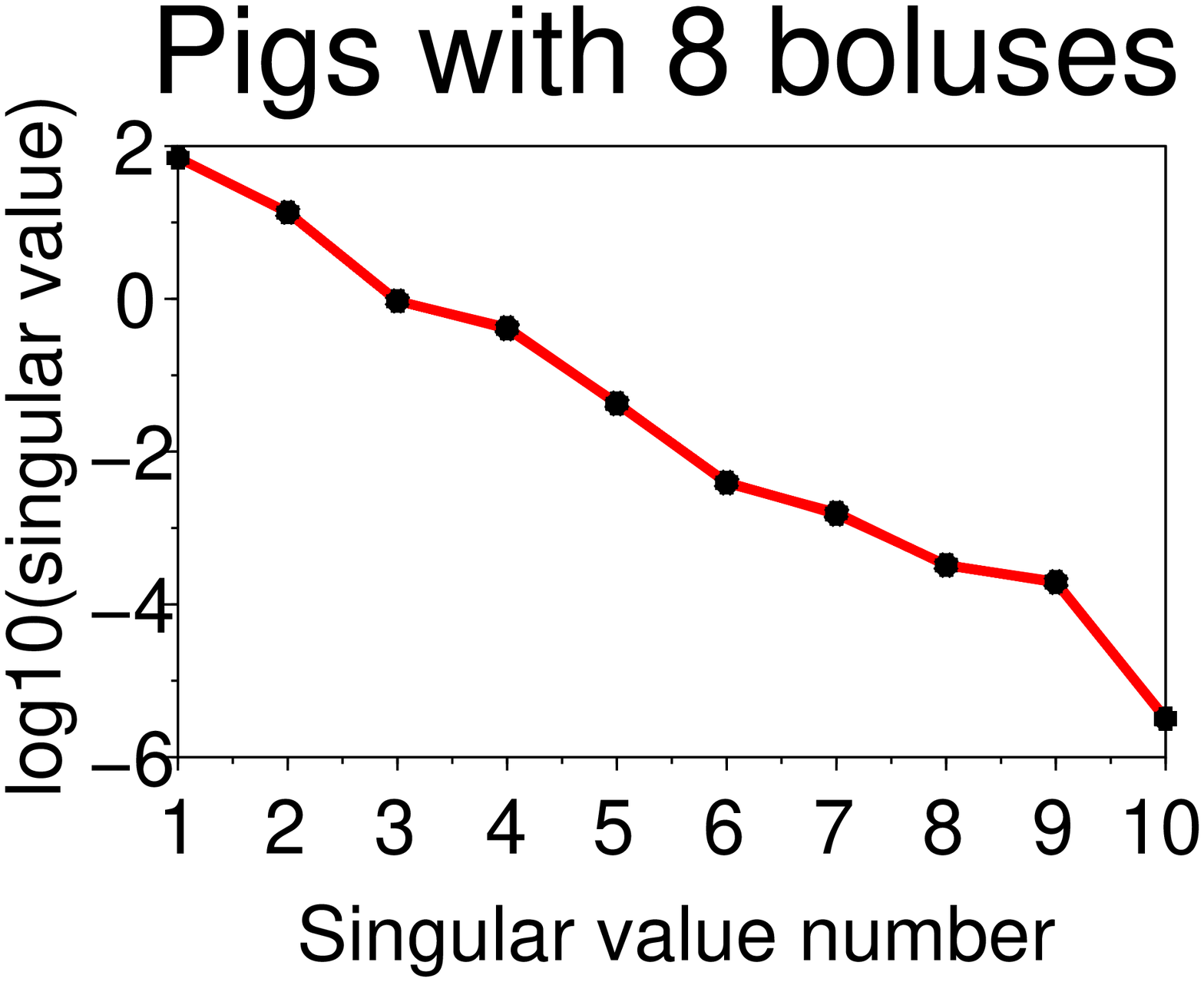}
\includegraphics[width=4.25cm]{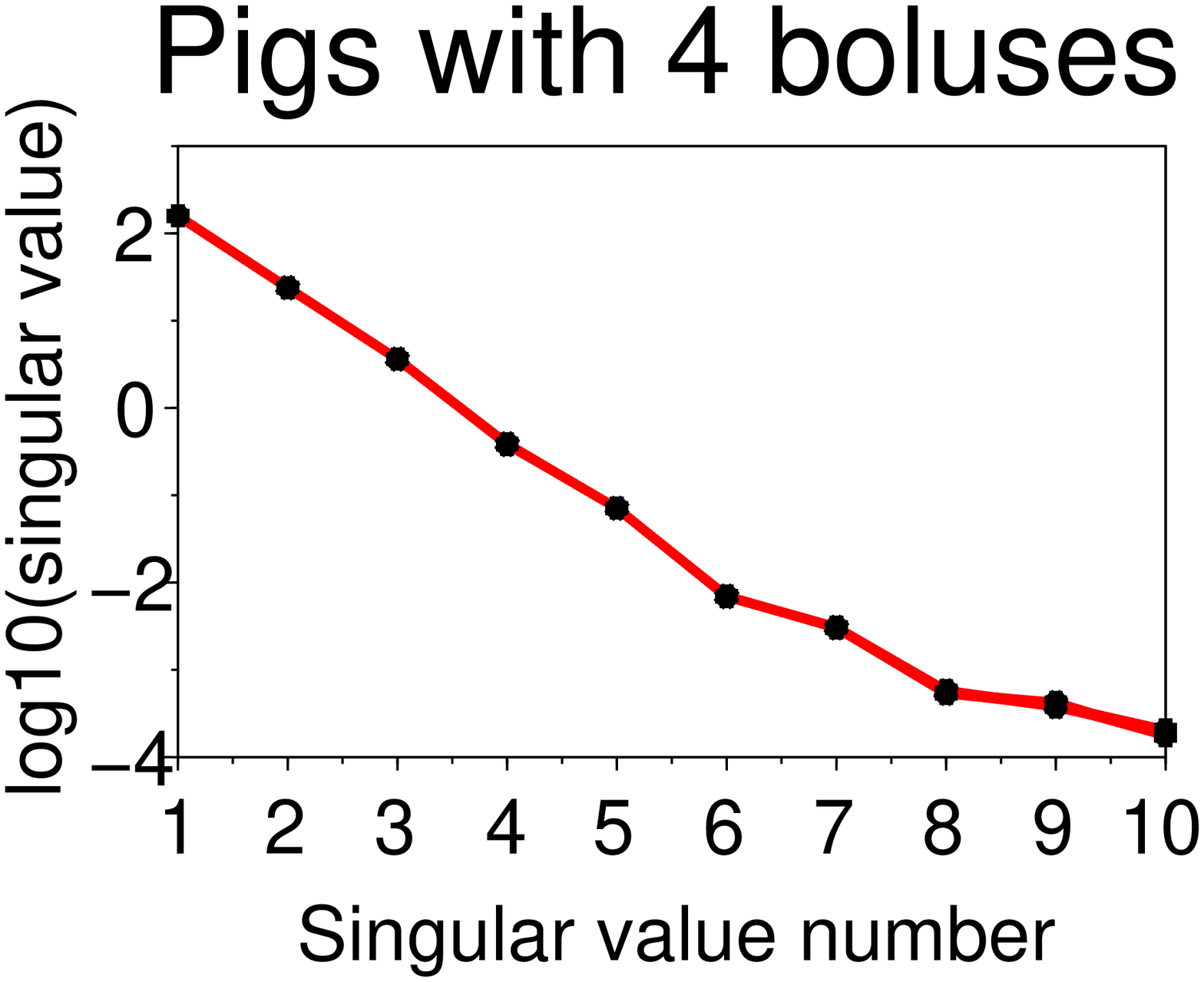}
\caption{Singular Value Decomposition of the reduced bihormonal--glucose model \eqref{red:glucose}. No singular value equal to zero or very small was found, suggesting that parameters are locally structurally identifiable. }\label{fig:SVDreduced}
\end{figure}

The Singular Value Decomposition method \cite{stigter2015fast},\cite{staal2019glucose} described in Supplementary material S.6 was performed for the reduced bihormonal--glucose model \eqref{red:glucose} combining cases. The result are depicted in Fig. \ref{fig:SVDreduced}. All singular values has order superior to $10^{-6}$, while for the bihormonal--glucose model \eqref{glucose}--\eqref{glycogen} there are singular values with order $10^{-9}$ or less. In conclusion, no parameter with lack of local structural identifiability was found for the reduced bihormonal--glucose model \eqref{red:glucose}.

\section{Data Collection and Methods}\label{section:data}

\subsection{Data Collection}

The reduced bihormonal--glucose model \eqref{red:glucose} was calibrated using data from experiments with pigs. Each experiment was carried out in about 8 hours. Blood glucose levels were measured at least every 5 minutes from IV blood samples collected in syringes and analyzed on a Radiometer
ABL 725 blood gas analyzer (Radiometer Medical ApS, Brønshøj, Denmark) \cite{dirnena2019intraperitoneal, dirnena2018intraperitoneal, am2020intraperitoneal}.

Glucose was intravenously infused at different rates during all the experiment. In order to simulate food intake, glucose infusion was increased and then decreased as a step function after the first half of the experiments (glucose infusion is depicted in Fig. \ref{fig:calibration54}). In the experiments for Pig 1 and Pig 4, there were several increments in glucose infusion to avoid hypoglycemia.

Insulin and glucagon boluses were pumped into the IP cavity. Porcine insulin and glucagon endogenous production were neglected for modeling, since they were suppressed by a combination of octreotide and pasireotide during the experiments.

\subsection{Parameter Estimation: Minimization Method}\label{subsec:minimization_method}

 Model calibration is performed in order to personalize the reduced bihormonal--glucose model \eqref{red:glucose} for each subject. Parameter estimation for each experiment was carried out using the Nelder--Mead algorithm to minimize the sum of square errors between the model and experimental data. The \emph{fminsearch} tool was used in Scilab to obtain parameter values that minimize the cost function
\begin{align*}
F(\boldsymbol\theta)=&\sum_{t\in T_{BGA}} \big[BGA(t)-G(t,\boldsymbol\theta)\big]^2,
\end{align*}where $\boldsymbol\theta$ is the vector of parameters to be estimated, $BGA(t)$ blood glucose measurements, $T_{BGA}$ the set of time--points at which glucose was measured, and $G(t,\boldsymbol\theta)$ the glucose state of the reduced bihormonal--glucose model \eqref{red:glucose} with the parameters in $\boldsymbol\theta$.

After the calibration, the Mean Square Error was computed for each experimental case:
\begin{align*}\hat\sigma^2=\frac{\sum_{t\in T_{BGA}}\big[BGA-G(t,\boldsymbol\theta)\big]^2}{n},\end{align*}where $n$ is the sample size.

To compare the accuracy of the model estimation, the BIC (Bayesian information criterion) value was calculated for each experiment case. BIC criterion considers the complexity of the model, i.e. the number of parameters and the sample size. The BIC value is defined as
$$BIC = n\cdot\log(\hat\sigma^2) + p_\# \cdot\log(n)$$
where $p_\#$ is the number of parameters. BIC can take negative values. The better the performance of the model, the lower the BIC value \cite{burnham2004multimodel}.

The parameters estimated are in Table \ref{table:parameters_reduced}. For the measurable state (i.e. glucose state), the initial condition can be taken as the last data point registered. For the non--measurable states (dimensionless states), the initial conditions are fixed and not estimated with the minimization method, because they can be readjusted by a change of variable without affecting the output of the system (the measurable state).

A second algorithm based on the quasi--Newton method was used to calibrate the parameters of the reduced bihormonal--glucose model \eqref{red:glucose}. The results are presented in Table S.6 of Supplementary Material S.10. The estimates obtained with the Nelder--Mead algorithm and the quasi--Newton method are close. The mean of the relative difference between the parameters is $0.11$.

\section{Model Calibration Results}\label{Section:reduced_simulations}

\begin{figure*}\centering
\includegraphics[width=17cm]{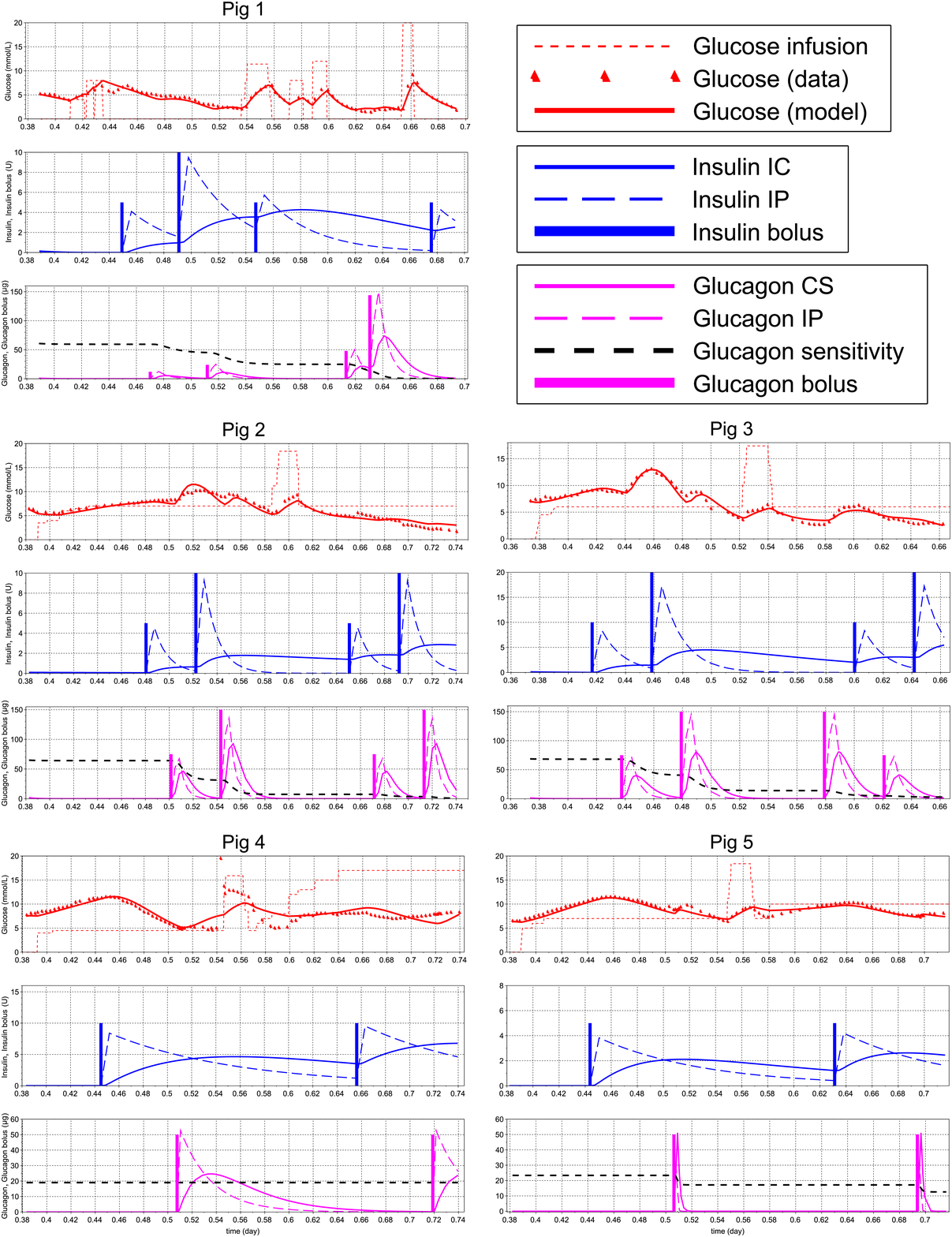}
\caption{Reduced bihormonal--glucose model \eqref{red:glucose} compared to experimental data. During the experiments, several intraperitoneal insulin and glucagon boluses were used to identify the nonlinear dynamics of both hormones and the glucose response. The parameters estimated are described in Table \ref{table:parameters_reduced}.}\label{fig:calibration54}
\end{figure*}

 The reduced bihormonal--glucose model \eqref{red:glucose} is personalized according to the experimental data of each subject. For this purpose, its parameters were estimated for each experimental case, as explained in Section \ref{section:data}. The numerical results are plotted in Fig. \ref{fig:calibration54}.

The BIC values for 4 out of the 5 experimental cases are -47.97 in average (see Table \ref{table:parameters_reduced}), which suggest that data are accurately approximated with a model of adequate complexity. The case with the largest BIC value (Pig 4) had an abnormal late and sharp response to the first glucagon bolus.

\begin{table}[!]\scriptsize
\begin{center}
\caption{Parameters estimated with the Nelder--Mead algorithm for the reduced bihormonal--glucose model \eqref{red:glucose}, parameter Coefficient of Variation (CV), Mean Square Error (MSE), and BIC values of the model approximation. The better the model performance to approximate the data, the lower the MSE and BIC values.}\label{table:parameters_reduced}
\begin{tabular}{cccccccc}
\hline
Parameter			&	Pig 1	&	Pig 2	&	Pig 3	&	Pig 4				&	Pig 5	&CV	\\
\hline
$k_1$				&	13.79	&	21.56	&	7.41	&	0.98				&	3.92	&0.87	\\
			
$\overline{k_{i_1}}$	&	171.68	&	168.68	&	181.79	&	115.36				&	185.46	&0.17	\\			

$\overline{k_H}$		&	38.50	&	44.02	&	33.20	&	28.63				&	655.46	&1.73	\\
	
$r_G$				&	4.73	&	2.79	&	2.21	&	1.77				&	1.73	&0.47	\\

$\overline{m_3}$	&	4.83	&	2.33	&	4.38	&	8.74				&	37.89	&1.28	\\			

$m_4$				&	27.84	&	85.05	&	64.96	&	9.49				&	12.30	&0.84	\\			

$q$					&	0.48	&	0.81	&	0.52	&	0.96				&	1.37	&0.44	\\
 				
${n}$				&	110.34	&	177.30	&	142.24	&	37.85				&	709.33	&1.15	\\				
$n_1$				&	138.95	&	200.09	&	177.44	&	38.52				&	2152.87	&1.67	\\				

$\overline{x_1}$ 		&	237.39	&	196.85	&	102.37	&	0.0014				&	3392.80	&1.86	\\
\hline
MSE	&	0.30	&	0.54	&	0.36	&	3.45	&	0.25\\				
\hline
BIC					&	-57.40	&	-8.15	&	-38.47	&	174.15				&	-87.85		\\	
\hline
\end{tabular}
\end{center}
\end{table}

The variability of parameters between individuals is due to the inter--subject variability of glucose dynamics \cite{heise2004lower, toffanin2019multiple}. To determine the range of variation for the parameters will require many experimental cases (which might only be possible with clinical cases) to make the estimations statistical meaningful.

The Mean Square Error between the data and the model was computed for the complete and reduced model. This show that the error obtained approximating the data with the reduced model is lower than with the complete model (compare Table \ref{table:parameters_reduced} and Table S.3 of Supplementary material S.4.

 Furthermore, the BIC criteria was used to compare the performance of the complete and reduced models to approximate the data of the full experiments (i.e. about 8 hours). BIC values are lower for the reduced model than for the complete model (see Table \ref{table:parameters_reduced} and Table S.3 of Supplementary material S.4.), indicating that the reduced model has a better performance for approximating the experimental data. 
This is shown with the number of cases had at the moment when the work was done. To consider more experimental cases or longer experiments may be done in future works.

\section{Discussion}\label{Section:discussion}

The bihormonal--glucose model \eqref{glucose}--\eqref{glycogen} and the reduced bihormonal--glucose model \eqref{red:glucose} have been presented in this work. Both models can be personalized to represent glucose nonlinear dynamics when intraperitoneal insulin and glucagon boluses are administered. Furthermore, these models can represent the variants in glucose infusions as performed during the experiments, leading to conclude that the models can approximate the dynamics of glucose even when there is food intake.

The local identifiability of the complete model was addressed. 
After local structural identification analysis, the reduction of the bihormonal--glucose model was accomplished through re--parametrization and state transformations. 
Moreover, a structural identification analysis was performed for the reduced bihormonal--glucose model \eqref{red:glucose} and no parameter was found to lack of local structural identifiability.

To address practical identifiability of the complete model, two insulin states active in glucose removal were reduced to a single state, and it has been shown that this does not imply less accuracy in the approximation of glucose measurements. Moreover, parameter profile likelihoods were computed for the reduced bihormonal--glucose model \eqref{red:glucose} and the results suggest that parameters are locally practically identifiable.

In this way, the number of parameters to be estimated was reduced for the bihormonal--glucose model \eqref{glucose}--\eqref{glycogen}, obtaining the reduced bihormonal--glucose model \eqref{red:glucose} which can represent experimental data with the same or better accuracy.

The reduced model was inferred taking into account that insulin and glucagon are not measured (to approximate insulin and glucagon measurements, the complete model may be used). Indeed, the reason for reducing is to have a simpler model which can be quickly calibrated after 2--3 hours of animal experiments where only blood glucose measurements are available, so that the model can be tested in MPC. 

In this work, two different methods were used to address parameter estimation: a minimization routine based on the Nelder--Mead algorithm and the quasi--Newton method. Other approaches can be tested to identify the models with prior knowledge on model parameters, for instance, using Bayesian methodologies \cite{cobelli1999minimal, visentin2016one} and Markov chain Monte Carlo techniques \cite{magni2004insulin, haidar2013stochastic}.

New technologies concerning the treatment of T1DM have been proposed to be incorporated in AP, for instance, to analyze body sounds to detect meal ingestion \cite{kolle2019feasibility, setti2019pilot}. But so far for this work, it is considered that only blood glucose measures will be available and that the information about glucose input is given.
Besides, there is low probability that the necessary data to make practically identifiable the parameters of the bihormonal--glucose model could be obtained from a non--invasive device, specially for the parameters related to the concentration of hormones within the intraperitoneal cavity.

Also, glucagon sensitivity might be related to the glycogen stored in the liver and the generation of glucose from amino acids (gluconeogenesis) \cite{chee2007closed}. As future work, the formulation of the glucagon sensitivity state can be revised, for instance to analyze if in larger periods of time there is an increment in this state and whether there is saturation. In these experiments, the repletion of glycogen in the liver could not be observed, probably due to the short duration of the experiments. However, it is expected that hepatic glycogen repletion will enhance the effects of glucagon doses \cite{blauw2016pharmacokinetics}.

On the other hand, the main risks observed in the use of IP insulin infusions are skin infections and interruption of the insulin supply attributable to catheter obstruction \cite{zisser2011clinical, jones2017artificial, aam2018effect, chakrabarty2019new}. The addition of a second hormone (glucagon) does not necessarily imply a second port, but rather a single catheter with two channels. Therefore, increasing one to two intraperitoneal hormones may not significantly increase the  risk of infection. Even if the efficiency of a single port (or a two--lumen catheter) to infuse dual hormones intraperitoneally is still to be established, the primary purpose of this work is to propose a suitable and relatively simple model to develop a controller that allows progress in animal experiments and future clinical trials.

Although the risk of these complications has been reduced with experience \cite{chakrabarty2019new}, the use of the intraperitoneal route is something that must be balanced with its benefit \cite{gehr2019continuous}. In--silico comparisons have shown that with intraperitoneal insulin infusions glucose levels can be controlled within the normal range by giving smaller glucose excursions after meals compared to subcutaneous insulin infusions, keeping blood glucose levels lower and preventing hypoglycemia even without the need for boluses before meals. \cite{fougner2007silico, jones2017artificial}. 

Furthermore, in animal experiments it has been observed that, compared to subcutaneous administration of hormones, intraperitoneal boluses induce faster effects on glucose levels while reducing the concentration of hormones in the circulatory system \cite{dirnena2018intraperitoneal, dirnena2019intraperitoneal, aam2020intraperitonealphd, am2020intraperitoneal}.

 Finally, normalizing glucose levels can eradicate the long--term adverse effects of diabetes that affect many patients after decades of disease.

\section{Conclusion}\label{Section:conclusion}
In summary, a low--order nonlinear bihormonal--glucose model accounting for intraperitoneal insulin and glucagon infusions is introduced in this work. The innovations of the model are the use of power--law kinetics for representing intraperitoneal insulin absorption and a separate glucagon sensitivity state. The model was reduced addressing its practical and structural lack of parameter identifiability, given glucose estimations from animal experiments. The parameters of the reduced model were found to exhibit local practical as well as structural identifiability. Both the complete and the reduced model can fit data from animal experiments, where insulin and glucagon boluses were introduced in the IP cavity, which is completely novel to the best of the knowledge of the authors.

\section*{Ethical approval}
 The animal experiments were approved by the Norwegian Food Safety Authority (FOTS number 12948) and were in accordance with “The Norwegian Regulation on Animal Experimentation” and
“Directive 2010/63/EU on the protection of animals used for scientific purposes”.

\section*{Acknowledgment}

The conceptualization of the models and numerical simulations were done by Dr. C. Lopez-Zazueta. 
Prof. Ø. Stavdahl, Dr. A. L. Fougner, and Dr. C. Lopez-Zazueta contributed to the discussion about the models and the parameter identification analysis. 
Prof. S. M. Carlsen, Dr. S. C. Christiansen, Dr. M. K. \AA m, Dr. A. L. Fougner, Mr. P. B{\"o}sch, and Dr. C. Lopez-Zazueta contributed in the discussion of the protocol for the experiments.
The experiments were performed by Dr. M. K. \AA m, Mr. P. B{\"o}sch, and Mrs. O. Lyng.
Dr. M. K. \AA m collected the data of the experiments with insulin and glucagon intraperitoneal boluses.

  \bibliographystyle{ieeetr}
  \bibliography{final_version_arxiv}


\appendix

{\vspace{1cm}{\hspace{-.5cm}\huge\bf{Supplementary Material}}}

\gdef\thesection{S.1}\section{Insulin--Glucose Submodel}\label{insulin_model}
The following model accounts for glucose dynamics of a subject with T1DM when insulin boluses are administered in the IP cavity: 
\begin{align}
\gdef\theequation{S.1}\label{insulin_submodel}
\frac{dG}{dt}  = & - [k_1 + k_I\cdot (I {+I_b}) { + k_{i_1}\cdot i_1}]\cdot G   + r_G\cdot Ra_G	\\\notag
\frac{dI}{dt} = & - m_1\cdot I + m_2\cdot i_1^p \\\notag
\frac{di_1}{dt} = & { - m_3\cdot i_1^q}+ m_4\cdot i_2 \\\notag
\frac{di_2}{dt} = & - m_4\cdot  i_2 + u_I.\notag
\end{align}
States and parameters are described in Table I
. In this case, where insulin measures were available, the state $I$ and the constant $I_b$ have units mU/L.

The extended insulin--glucose model \eqref{insulin_submodel} is a generalized version of model (2) in \cite{lopez2019simple}. The extended model allows two insulin compartments, $I$ and $i_1$, to be active in glucose--removal. It also integrates a nonlinear insulin consumption / transport term for the intermediate insulin compartment (the power $q$ for the state $i_1$). 

The extended insulin--glucose model \eqref{insulin_submodel} was calibrated with data from experiments with pigs, where insulin boluses were injected in the intraperitoneal cavity \cite{dirnena2019intraperitoneal}. The parameters estimated are in Table \ref{table:model_insulin}. The approximation was carried out as in \cite{lopez2019simple}{\color{black}, i.e. parameter estimation was carried out using the Nelder--Mead algorithm to minimize the sum of square errors between the model and data. The \emph{fminsearch} tool was used in Scilab to obtain parameter values that minimize the cost function
\begin{align*}
F(\boldsymbol\theta)=&\sum_{t\in T_{BGA}} \big[BGA(t)-G(t,\boldsymbol\theta)\big]^2\\
&+\sum_{t\in T_I}[IM(t)-I(t,\boldsymbol\theta)]^2,
\end{align*}where $\boldsymbol\theta$ is the vector of parameters to be estimated, $BGA(t)$ blood glucose measurements, $T_{BGA}$ the set of time--points at which glucose was measured, $G(t,\boldsymbol\theta)$ the glucose state in model \eqref{insulin_submodel} with the parameters in $\boldsymbol\theta$, $IM(t)$ blood insulin measurements, $T_I$ the set of time--points at which insulin was measured, $I(t,\boldsymbol\theta)$ the insulin state in model \eqref{insulin_submodel} with the parameters in $\boldsymbol\theta$.}

\begin{figure}[!]\centering
\hspace{-.4cm}\includegraphics[width=10cm]{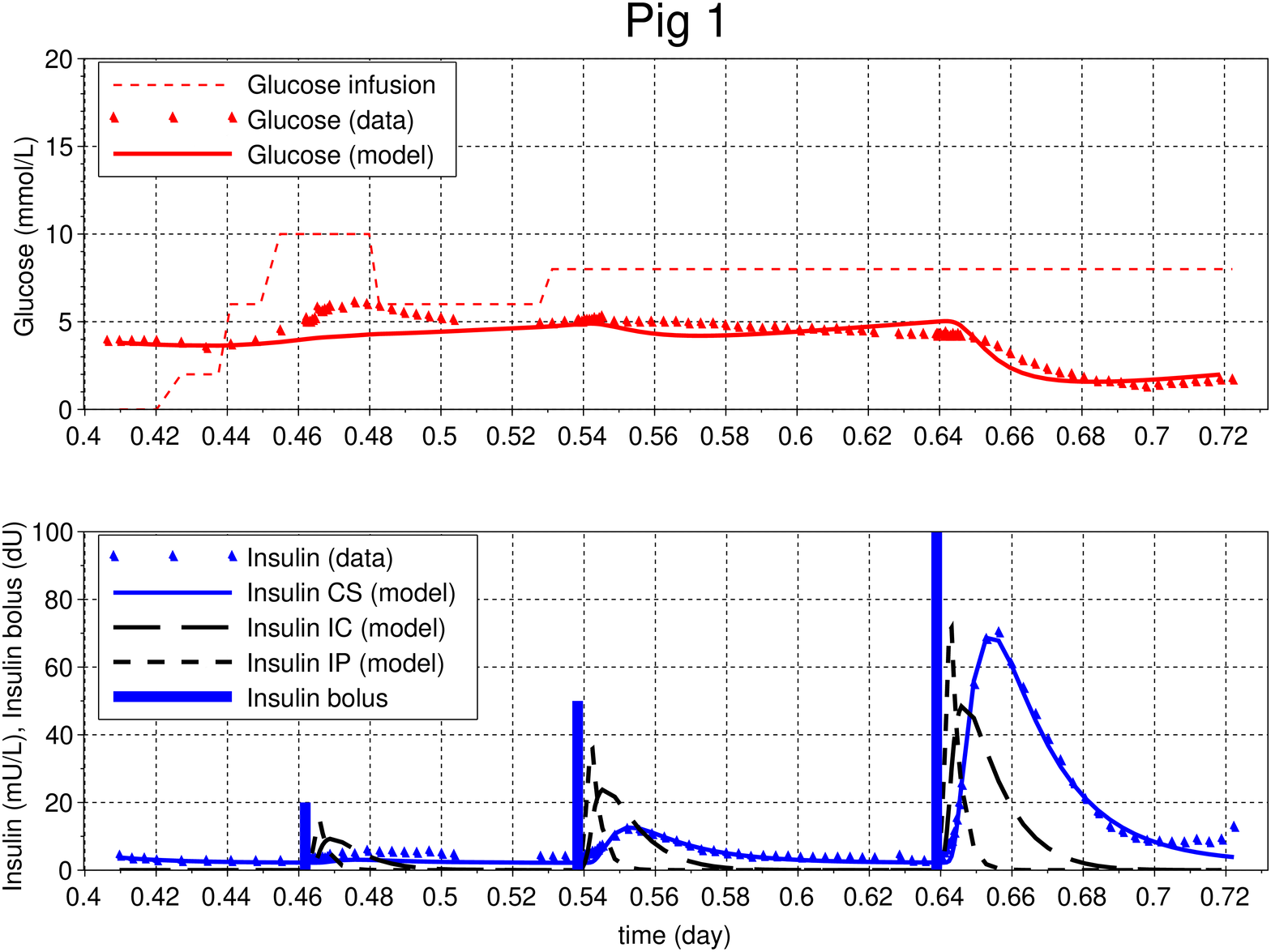}\\

\hspace{-.4cm}\includegraphics[width=10cm]{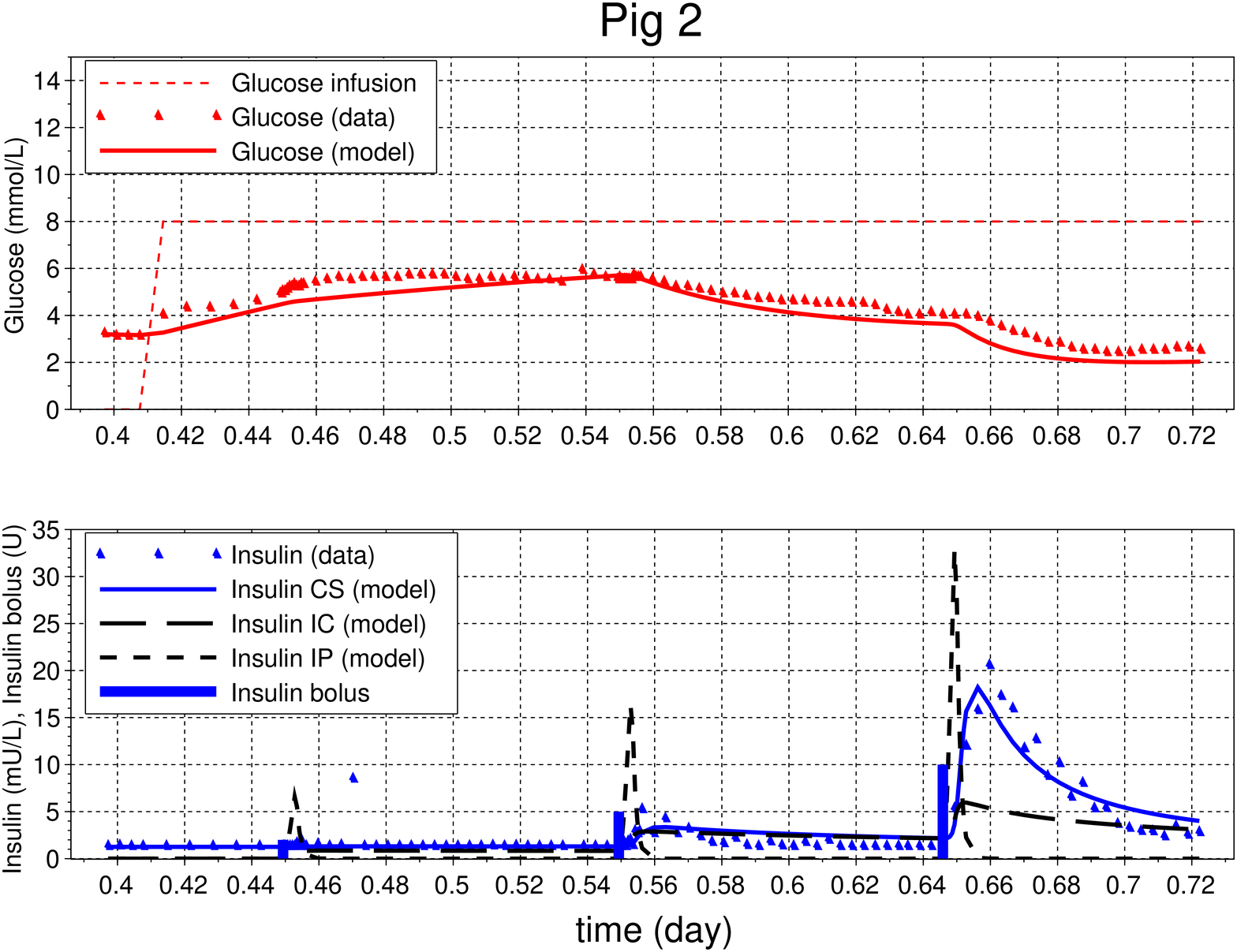}\\

\hspace{-.4cm}\includegraphics[width=10cm]{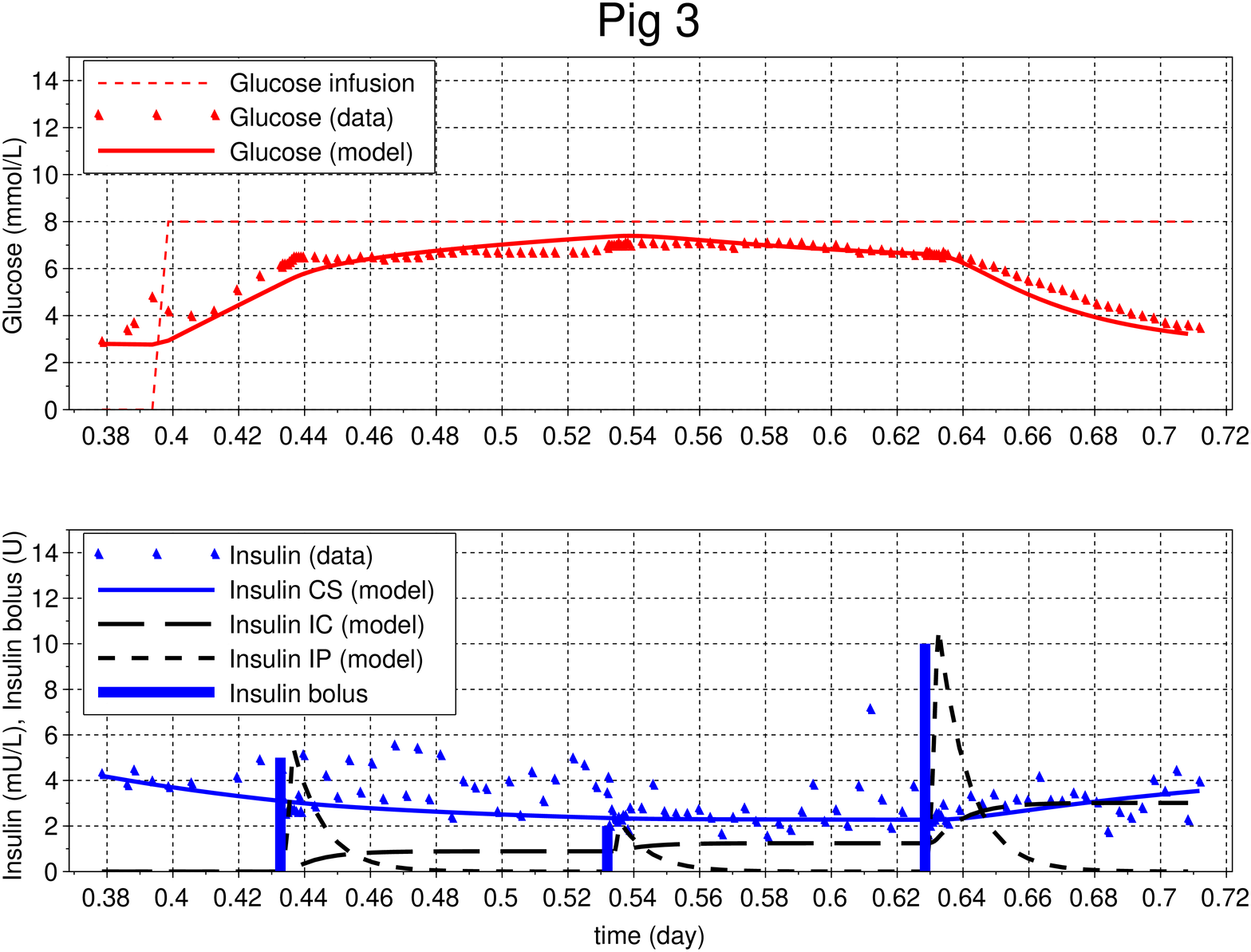}
\caption{Extended insulin--glucose model \eqref{insulin_submodel} compared to experimental data. Insulin boluses were introduced in the intraperitoneal cavity for these experiments. The parameters estimated are described in Table \ref{table:model_insulin}.}\label{figure:insulin}
\end{figure}

Results are depicted in Figure \ref{figure:insulin}. The extended insulin--glucose model \eqref{insulin_submodel} turns out to have lower BIC values (see Table \ref{table:model_insulin}) than model (2) in \cite{lopez2019simple} (this last only considers one insulin--state active in glucose removal).

\gdef\thefigure{S.1}

\gdef\thetable{S.1}
\begin{table*}
\begin{center}
\caption{Parameters estimated for the extended insulin--glucose model \eqref{insulin_submodel}. The minimization method to approximate data is the same described in \cite{lopez2019simple}. BIC values of the extended model are lower than model (2) in \cite{lopez2019simple}. MSE is the Mean Square Error of the model approximation, {\color{black}Var the parameter variance} and CV the parameter Coefficient of Variation (i.e. the standard deviation divided by the mean)}\label{table:model_insulin}
\begin{tabular}{ccccccc}
\hline
Parameter	&	Pig 1				&	Pig 2				&	Pig 3				&	Units		&{\color{black}Var}	&  CV	\\
\hline
$k_1$		&	0.0062				&	0.55				&	0.74				&	1/d		& {\color{black}0.10}& 		0.88	\\

$k_I$		&	0.82				&	0.51				&	0.000002			&	$\frac{\text{L}}{\text{mU.d}}$	&{\color{black}0.11}&	0.93\\

$k_{i_1}$	&	0.09				&	0.10				&	0.08				&	1/d	& {\color{black}0.00}	&		0.11	\\

$r_G$		&	0.61				&	0.98				&	1.84				&	$\frac{\text{h}}{\text{L.d}}$ &{\color{black}0.27}	& 	0.55\\

$m_1$		&	58.02				&	497.52				&	9.43				&	1/d		&{\color{black}48194.79}	&	1.43	\\

$m_2$		&	0.05				&	0.0004				&	0.00199				&	$\frac{\text{mU}}{\text{L.d}}$ & {\color{black}0.00}&	1.61\\

$m_3$		&	146.65				&	0.0002				&	0.00246				&	1/d &{\color{black}4779.07}	&	1.73	\\

$m_4$		&	381.56				&	900.18				&	107.70			&	1/d 	&{\color{black}107998.95}	& 	0.87		\\

$p$ 		&	2.56				&	2.99				&	1.65				&	-		&{\color{black}0.31}		&	0.29\\

$q$ 		&	0.90				&	3.02				&	1.54				&	-	&{\color{black}0.79}		&	0.60\\
\hline
MSE	&	1.73	&	1.03	&	0.82	&	\\
\hline
BIC			&	168.01				&	60.11				&	5.69				&						\\
\hline
\end{tabular}
\end{center}
\end{table*}

\emph{Note.} Glucose infusion was constant for these experiments, except at the beginning of the experiment with Pig 1. For Pig 1, since it was the first case, at the beginning of the experiment the value of glucose infusion was changed until reaching a value that would stabilize pig glucose levels (before insulin infusion). The data of glucose infusion was put in the model as reported from the experiment, but it is not known why the first part of the model cannot match the data as for the rest of the experiment.

\gdef\thesection{S.2}
\section{Glucagon--Glucose Submodel}\label{glucagon_model}

The following model was conceived to approximate experimental data where glucagon boluses were introduced in the IP cavity and subcutaneously (for details about the animal trials, see \cite{am2020intraperitoneal}):
\begin{align}
\gdef\theequation{S.2}\label{glucagon_subsystem}
\frac{dG}{dt}  = & - k_1 \cdot G + k_H\cdot (H+H_b) \cdot \xi  \\\notag
\frac{dH}{dt} = & - n\cdot H+ {n_4}\cdot h_2+ {n_2 \cdot h_1} 	\\\notag
\frac{dh_1}{dt} = & - n_1\cdot h_1 +  u_H	\\\notag
\frac{dh_2}{dt} = & - n_3\cdot h_2 +  u_{H,SC}	\\\notag
\frac{d\xi}{dt} = & - x_1\cdot H\cdot \xi,  \notag
\end{align}
where $h_2$ is the subcutaneous glucagon state and $u_{H,SC}$ the subcutaneous glucagon input. The states of the glucagon--glucose model \eqref{glucagon_subsystem} are described in Table I
.  In this case, where glucagon measures were available, the state $H$ and the constant $H_b$ have units pmol/L.

In contrast to insulin dynamics, the transport of glucagon between compartments did not exhibit a marked non--linear behavior. However, the response of glucose to glucagon boluses was variable and nonlinear. For these reason, a new auxiliary state $\xi$ was included to represent this nonlinearity in glucose dynamics. 

The glucagon sensitivity $\xi$ was observed to decrease during the time of the experiments. But in longer periods of times it can be hypothesized that this sensitivity is restored, by the storage of glucose in the liver in  presence of insulin for example. However, for the experiments with only glucagon boluses it was assumed that there was no insulin in the system, since the hormone endogenous secretion was suppressed with drugs. 

Furthermore, from the experimental measurements in \cite{dirnena2019intraperitoneal} and \cite{am2020intraperitoneal}, it was observed that glucagon transport to the circulatory system has faster dynamics than insulin. Blood insulin concentration requires about 30 minutes to reach a maximum after one bolus, while blood glucagon takes approximately 15 minutes. Therefore, instead of using three compartments as for the insulin transport, only two states are assumed for glucagon. Indeed, dynamic systems involving fast and slow rates can be reduced through singular perturbation techniques, eliminating the states with fast dynamics \cite{verhulst2007singular, lopez2019dynamical}.

Experimental data from pigs where intraperitoneal and subcutaneous glucagon boluses were administered \cite{am2020intraperitoneal} was used to calibrated  the glucagon--glucose model \eqref{glucagon_subsystem}. The parameters estimated are in Table \ref{table:glucagon}. 
{\color{black}Parameter estimation was carried out using the Nelder--Mead algorithm to minimize the sum of square errors between the model and data. The \emph{fminsearch} tool was used in Scilab to obtain parameter values that minimize the cost function}
\begin{align*}
F_{GN}(\boldsymbol\theta)=&\sum_{t\in T_G} \Big[ BGA(t)-G(t,\boldsymbol\theta) \Big]^2\\
&+\sum_{t\in T_{GN}} \Big[ GNM(t)-H(t,\boldsymbol\theta) \Big]^2,
\end{align*}
{\color{black}where $BGA(t)$ denotes blood glucose measurements, $T_G$ the set of time--points at which glucose was measured, $G(t,\boldsymbol\theta)$ the glucose state in model \eqref{glucagon_subsystem} with the parameters in $\boldsymbol\theta$, $GNM(t)$ blood glucagon measurements, $T_{GN}$ the set of time--points at which glucagon was measured and $H(t,\boldsymbol\theta)$ the glucagon state in model \eqref{glucagon_subsystem} with the parameters in $\boldsymbol\theta$. }

\gdef\thefigure{S.2}
\begin{figure}[!]\centering
\includegraphics[width=8.5cm]{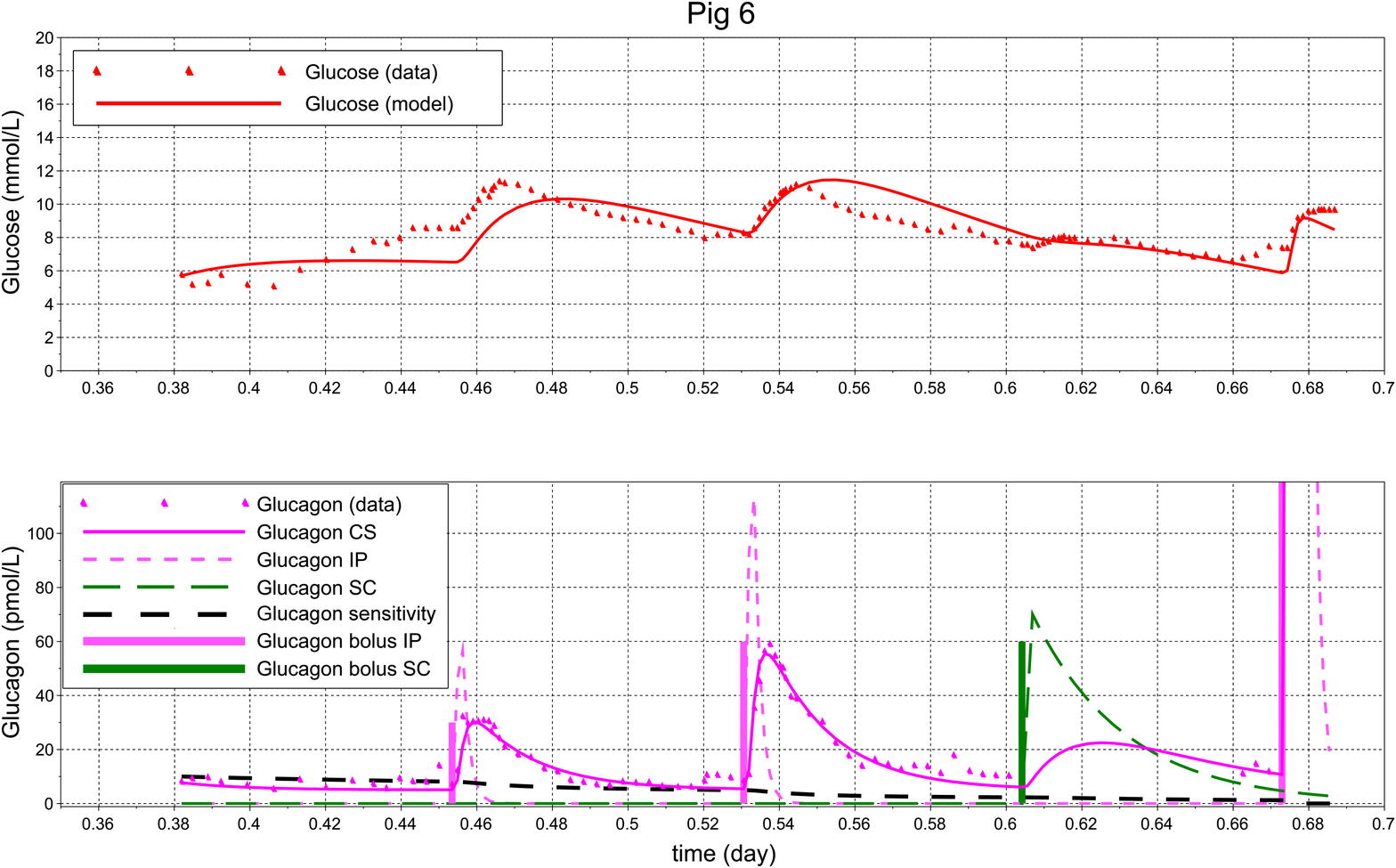}\vspace{.2cm}
\includegraphics[width=8.5cm]{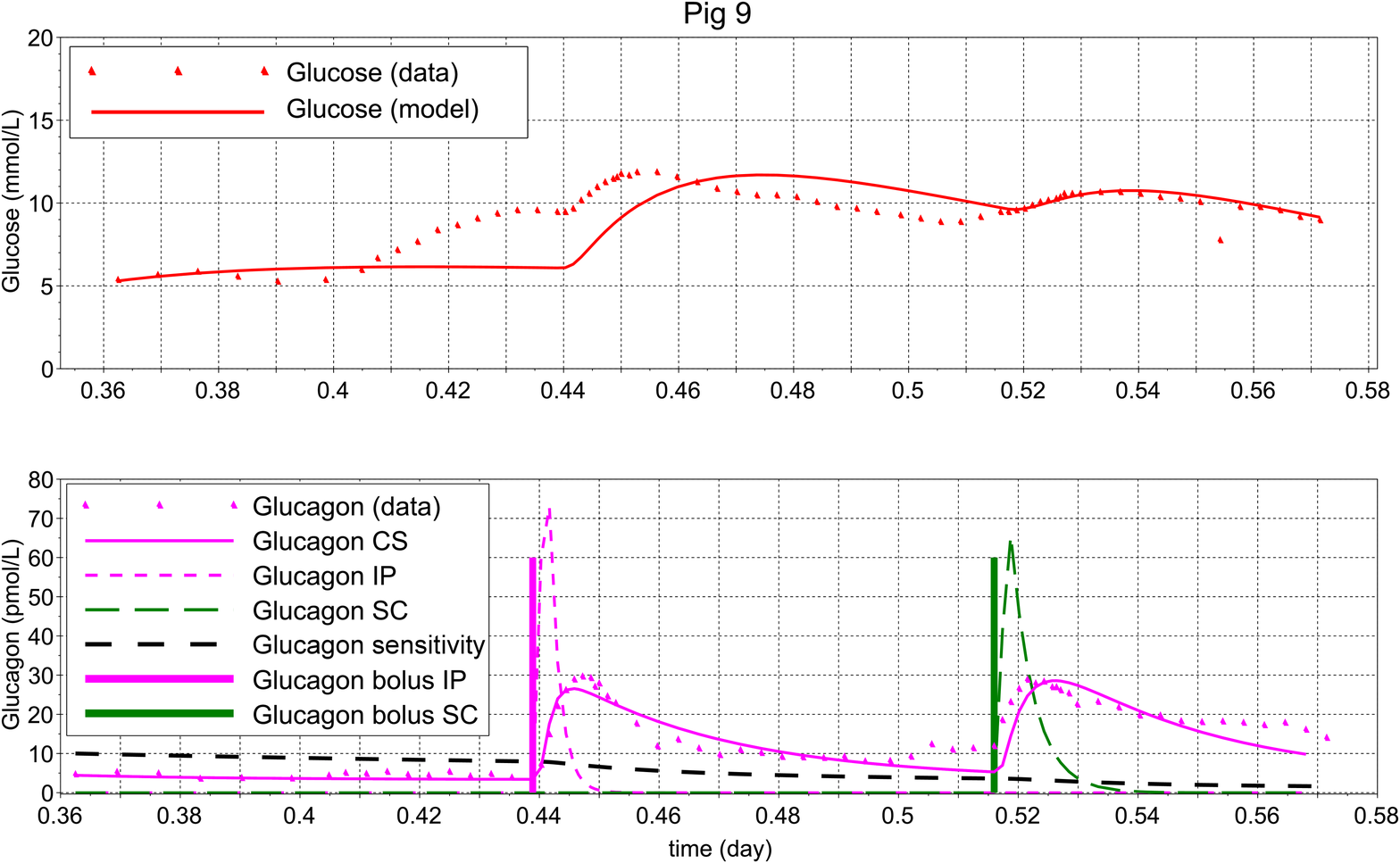}\vspace{.2cm}
\includegraphics[width=8.5cm]{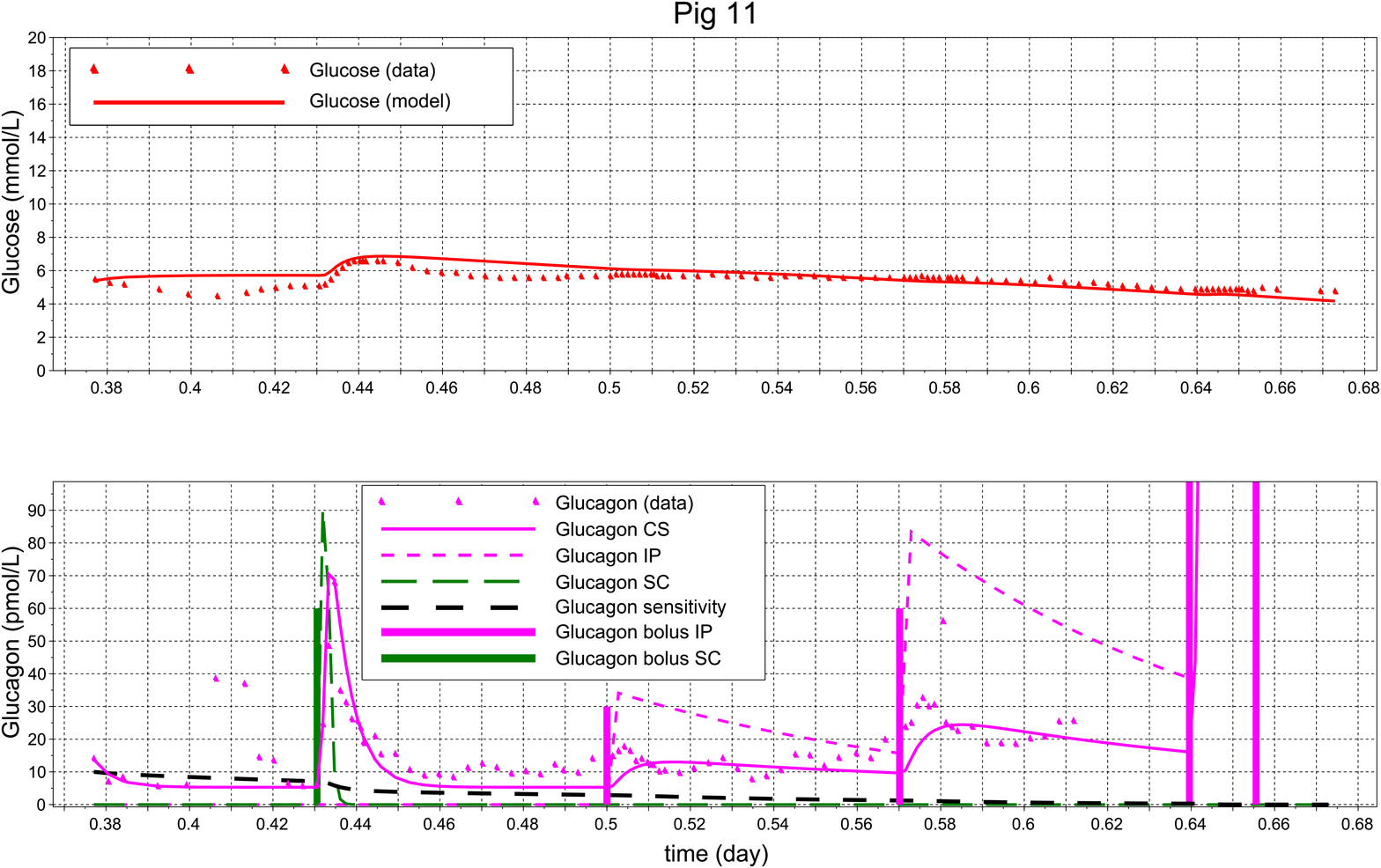}
\caption{Glucagon--glucose model \eqref{glucagon_subsystem} compared with experimental data. The parameters to fit the data are described in Table \ref{table:glucagon}. At the beginning of experiments for Pig 6 and Pig 9, there was an unexpected increase in glucose levels (notice that in these cases there was no glucose infusion). For this reason, the model fits less the data from the beginning of these experiments.}\label{figure:glucagon}
\end{figure}

\gdef\thetable{S.2}
\begin{table*}
\begin{center}
\caption{Parameters for the glucagon--glucose model \eqref{glucagon_subsystem}. The parameters were obtained minimizing the cost function $F_{GN}$. BIC values were computed to evaluate the performance of the model to fit data. MSE is the Mean Square Error of the model approximation, {\color{black}Var the parameter variance} and CV the parameter Coefficient of Variation (i.e. the standard deviation divided by the mean)}\label{table:glucagon}
\begin{tabular}{cccccccc}
\hline
Parameter	&	Pig 6	&	Pig 9	&	Pig 11	&	Units & {\color{black}Var} & CV\\
\hline
$k_1$		&	11.54	&	11.43	&	3.75	&	1/d	& {\color{black}13.30}&0.50	\\

$k_H$		&	16.84	&	23.68	&	5.32	&	$\frac{\text{mmol}}{\text{pmol.d}}$& {\color{black}57.40}&	0.61\\

$n$			&	59.41	&	36.05	&	212.79	&	1/d& {\color{black}6145.35}	&0.93\\

$n_1$		&	561.17	&	568.98	&	11.54	&	1/d& {\color{black}68099.28}&	0.84\\

$n_2$		&	4438.52	&	1912.86	&	467.40	&	$\frac{pmol}{L.d}$& {\color{black}2693123.01}&	0.88\\

$n_3$		&	41.00	&	257.62	&	1314.97	&	1/d& {\color{black}309768.09}	&1.27\\

$n_4$		&	279.83	&	1019.12	&	16016.44&	$\frac{\text{pmol}}{\text{L.d}}$& {\color{black}52567450.13}&	1.54\\

$x_1$ 		&	0.53	&	0.80	&	1.064	&	$\frac{\text{L}}{\text{pmol.d}}$& {\color{black}0.05}	&0.33\\
\hline
MSE	&	4.45	&	7.06	&	48.64	&	\\
\hline
BIC			&	321.08	&	308.20	&	860.78	&		\\
\hline
\end{tabular}
\end{center}
\end{table*}

 Variability in the parameters is related to the inter--subject variability of glucose homeostasis \cite{heise2004lower, toffanin2019multiple}.
BIC value considers the logarithm of the Mean Square Error, thus, when approximations are less accurate, this exponentially grows.

On the other hand, a hypothetical endogenous production of glucose can be added in the glucagon--glucose model \eqref{glucagon_subsystem} to improve data fitting in the first part of the experiments, where glucose increase was not expected due to the hormone suppression drugs administered. 

\gdef\thesection{S.3}
\section{Notes on the Bihormonal--Glucose Model}

The deduction of the insulin--glucose and glucagon--glucose submodels are described in \cite{lopez2019simple}, Sections \ref{insulin_model} and \ref{glucagon_model} of this work. Both submodels rely on basic features of glucose homeostasis: insulin--dependent and insulin--independent consumption of glucose, storage of glucose (in the form of glycogen) in the presence of insulin and its release in the presence of glucagon, inducing an increment in glucose levels, as well as the transport of hormones between compartments.

The submodels were deduced from simple linear models \cite{candas1994adaptive, yipintsoi1973mathematical, lopez2019simple}, heuristically adapting them to approximate experimental data obtained with intraperitoneal hormone infusions.

 All the model parameters, the insulin basal value $I_b$ and the glucagon basal value $ H_b$ are nonnegative. The inputs $Ra_G$, $u_I$ and $u_H$ are nonnegative functions.

 Insulin and glucagon basal values were added in the model because experimental data indicates basal values larger than the lowest thresholds of the insulin and glucagon analysis kits (1.25 mU/L for insulin and  2.3 pmol/L for glucagon), even in the presence of drugs suppressing insulin and glucagon production (see Figures \ref{figure:insulin} and \ref{figure:glucagon}). Thus, these basal values have to be considered as hormone levels that remain almost constant in a time interval where there is no perturbation from exogenous hormone injections.

Note that insulin and glucagon have different functions in glucose homeostasis, therefore their equations are also different.

The term $x_2\cdot G\cdot I$ is included to represent the storage of glucose in the form of glycogen. This storage is caused by insulin and it occurs mainly in the liver and skeletal muscles \cite{chee2007closed, wasserman2009four, blauw2016pharmacokinetics, cinar2018artificial, dirnena2019intraperitoneal}). This term allows to increment glucagon sensitivity when glucose and insulin concentrations increase, that it is when the storage of glucose in form of glycogen occurs. Then, glycogen reserves decrease when glucagon induces the breakdown of glycogen (glycogenolysis).  Basically, if there are no glycogen reserves, glucose response to glucagon will be minimized.

Notice that the term $x_2\cdot G\cdot I$ in the glucagon sensitivity equation does not contradict the conclusion of \cite{hinshaw2015glucagon} of that the difference in glucagon dose response based on glycemia was not significant at euglycemia or hypoglycemia, because the response depends on the glycogen stored and not on glucose levels.

 In comparison with the glucagon--glucose model \eqref{glucagon_subsystem}, for the bihormonal--glucose model (1)--(7) 
 one term allowing the increase of glucagon sensitivity in presence of insulin and glucose is considered. Moreover, the compartments and the terms related to the subcutaneous glucagon infusion were omitted, because for the experiments with both insulin and glucagon boluses all the injections were made in the IP cavity.

For the sake of simplicity, diffusion of glucose and hormones in a compartment is considered to be homogeneous. On the other hand, when a state is not considered as dimensionless, it is representing the concentration of a substance (glucose, insulin or glucagon) in a compartment. Dimensionless states are considered as auxiliary states at each compartment that allow to describe hormone transport from the IP compartment to the circulatory system and their effects on glucose levels.

It is assumed that the bihormonal--glucose model is a reduced system, since the model is a simplified version of insulin, glucagon and glucose dynamics and a power--law is used. To compute mass conservation in a reduced model is more complex than in a closed system or a non--reduced compartmental model, which accounts for all the compartments, all transfers between them or out of the system and dilution terms \cite{jacquez1993qualitative}. Typically, the mass conservation stated in a complete system is hidden in the reduced one and to recover it would require knowing the parameters of the complete model.

To determine mass conservation in the  bihormonal--glucose model (1)--(7) 
would require the statement of a larger, more complex and detailed model, that accounts for all the processes to which these hormones are submitted. However, such a complete model could not be validated, assuring the hypothesis at each stage to be correct, with the data available, because there are no taken measurements from internal compartments such as the different quadrants in the intraperitoneal cavity \cite{aam2018effect}, the portal vein, etc.   
Furthermore, from glucose measurements the amount of insulin and glucagon can only be estimated. Determining with accurate precision the concentration of these hormones to validate mass conservation would require considering more factors, for instance, the quantity of hormone receptors.

More detailed models have been proposed with the purpose of describing hormone concentration through the body (e.g. UVA/Padova model \cite{man2014uva}). But they are usually of large dimension, which carries out other difficulties as lack of identifiability and the non--practical implementation of these in a controller  \cite{cobelli2009diabetes, messori2018individualized}. Here, it is presented an example of the use of Power--law to describe in a reduced way the nonlinear transport of insulin from the IP cavity to the blood stream.

The bihormonal--glucose model (1)--(7) 
was calibrated with blood glucose measurements from experiments with pigs, where several insulin and glucagon boluses where introduced in the IP cavity. The parameters estimated are in Table \ref{table:bihormonal}. The experimental framework is described in Section IV
. For each experiment, 16 model parameters were estimated using the minimization method explained in Section IV--B
. Numerical simulations show the model (1)--(7) 
can approximate the data (see  \ref{app_calibration}). 

Although the bihormonal--glucose model (1)--(7) 
has few parameters that can be estimated to approximate glucose measurements, one objective of this work is to verify whether its parameters can be identified after alternating intraperitoneal insulin and glucagon boluses in an experiment. Moreover, in the initial experiments blood samples were analyzed to find the insulin and glucagon concentrations, but in future experiments such measurements will not be available, so the model might be reduced in order to fit a smaller set of data. For this purpose, an identification analysis is presented in the following two Sections.

\gdef\thesection{S.4}
\section{Calibration of the Bihormonal--Glucose Model}\label{app_calibration}

\begin{figure*}\centering
\includegraphics[width=7cm]{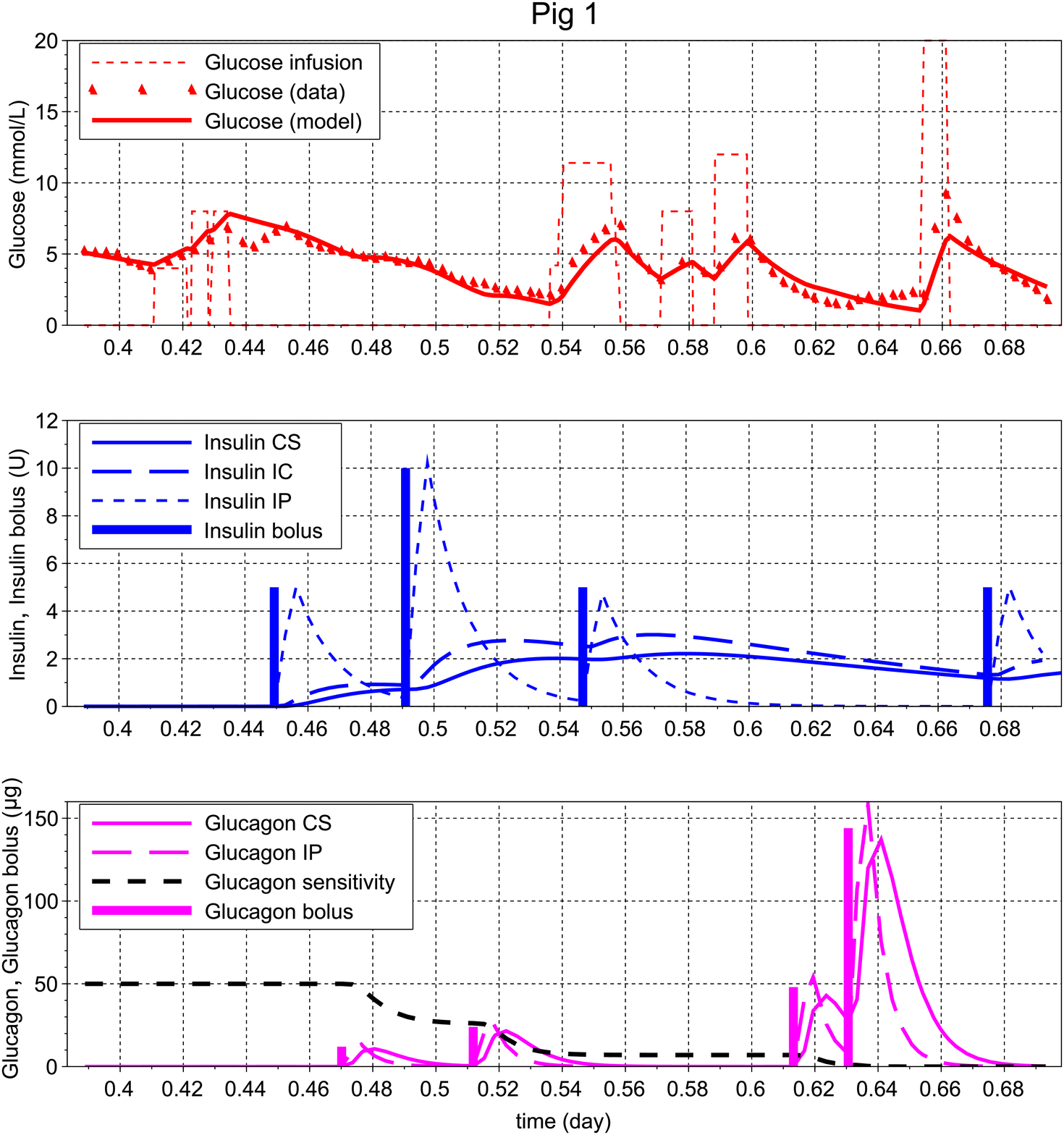}\includegraphics[width=7cm]{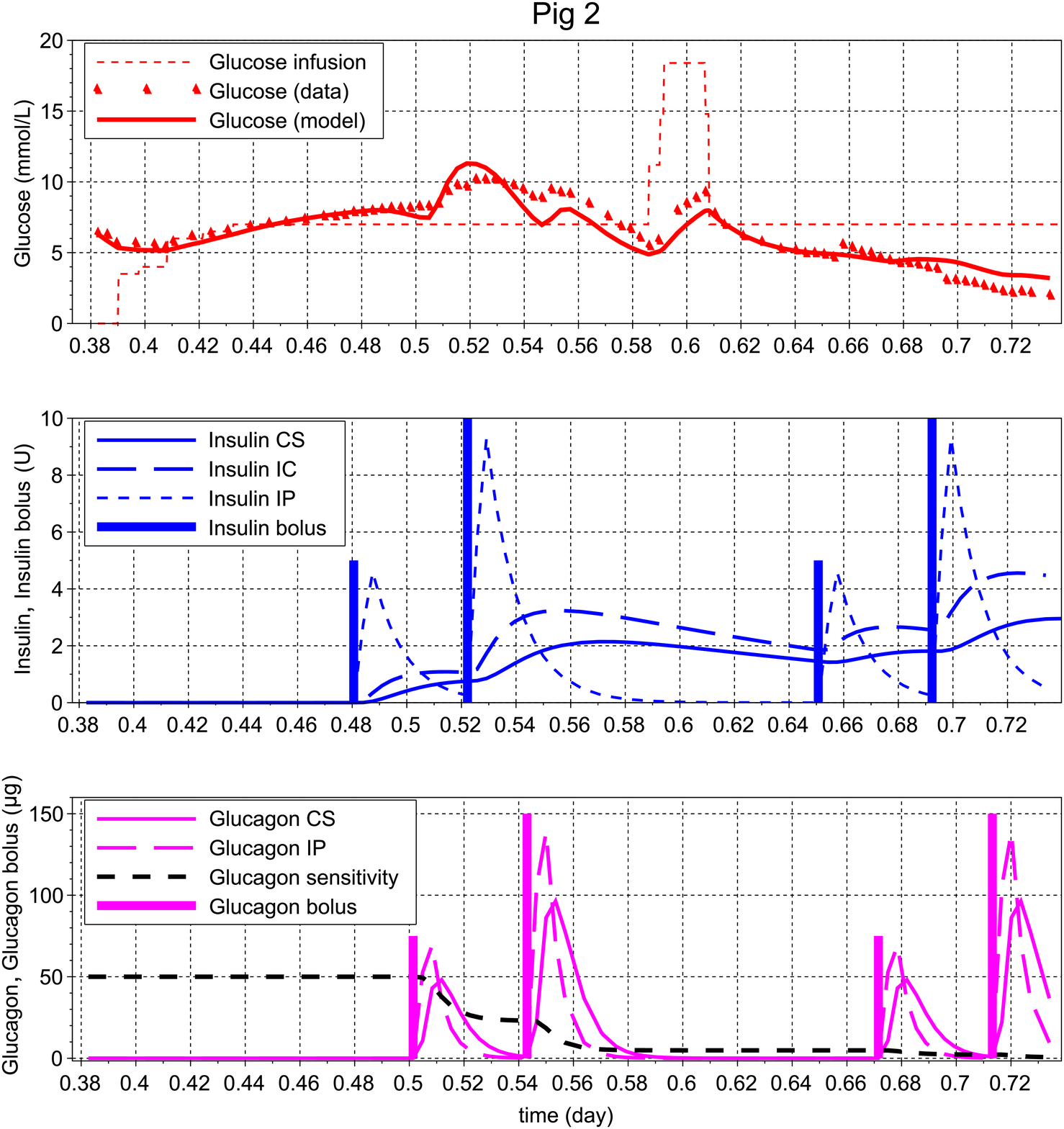}\\
\includegraphics[width=7cm]{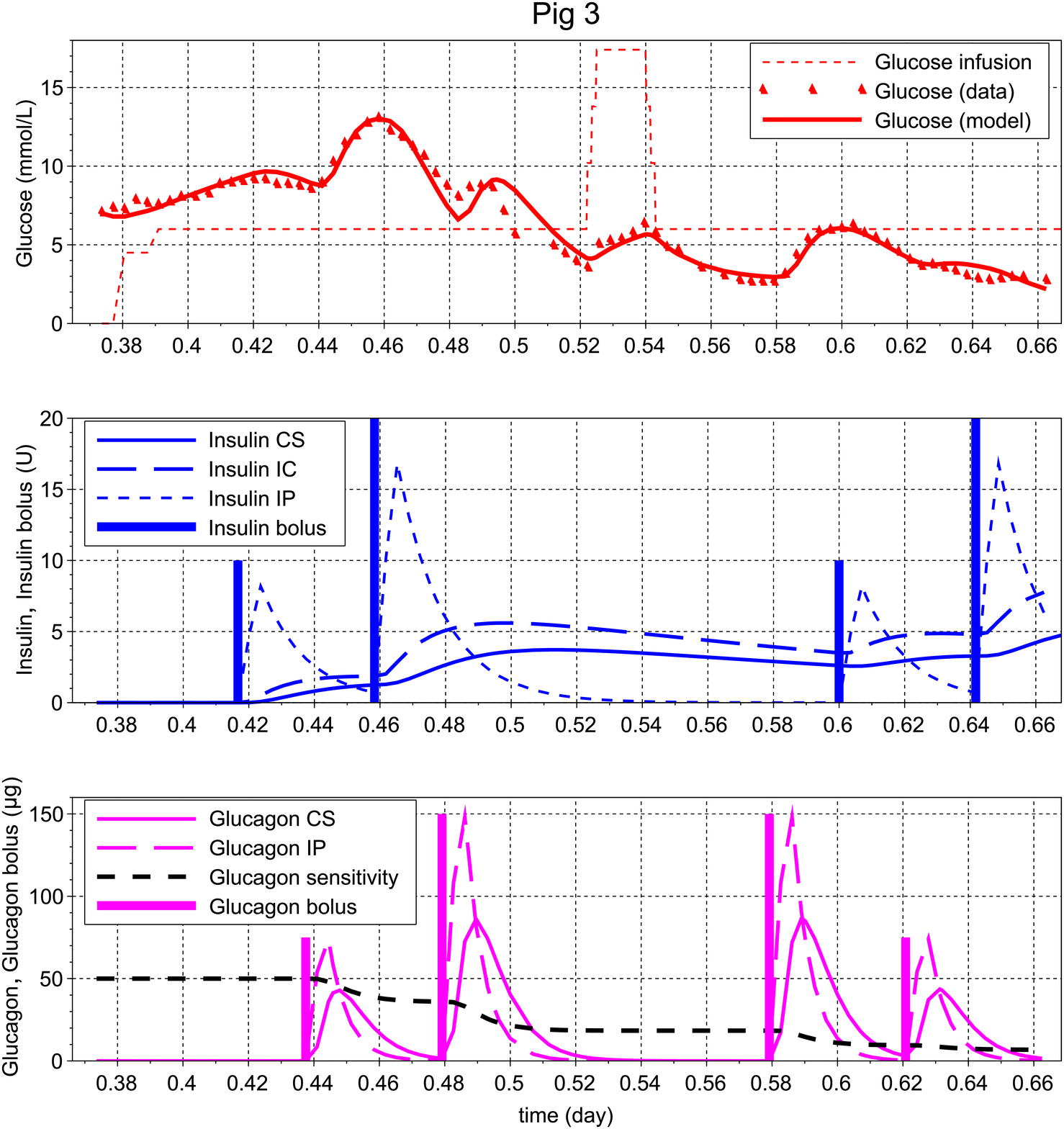}\\
\includegraphics[width=7cm]{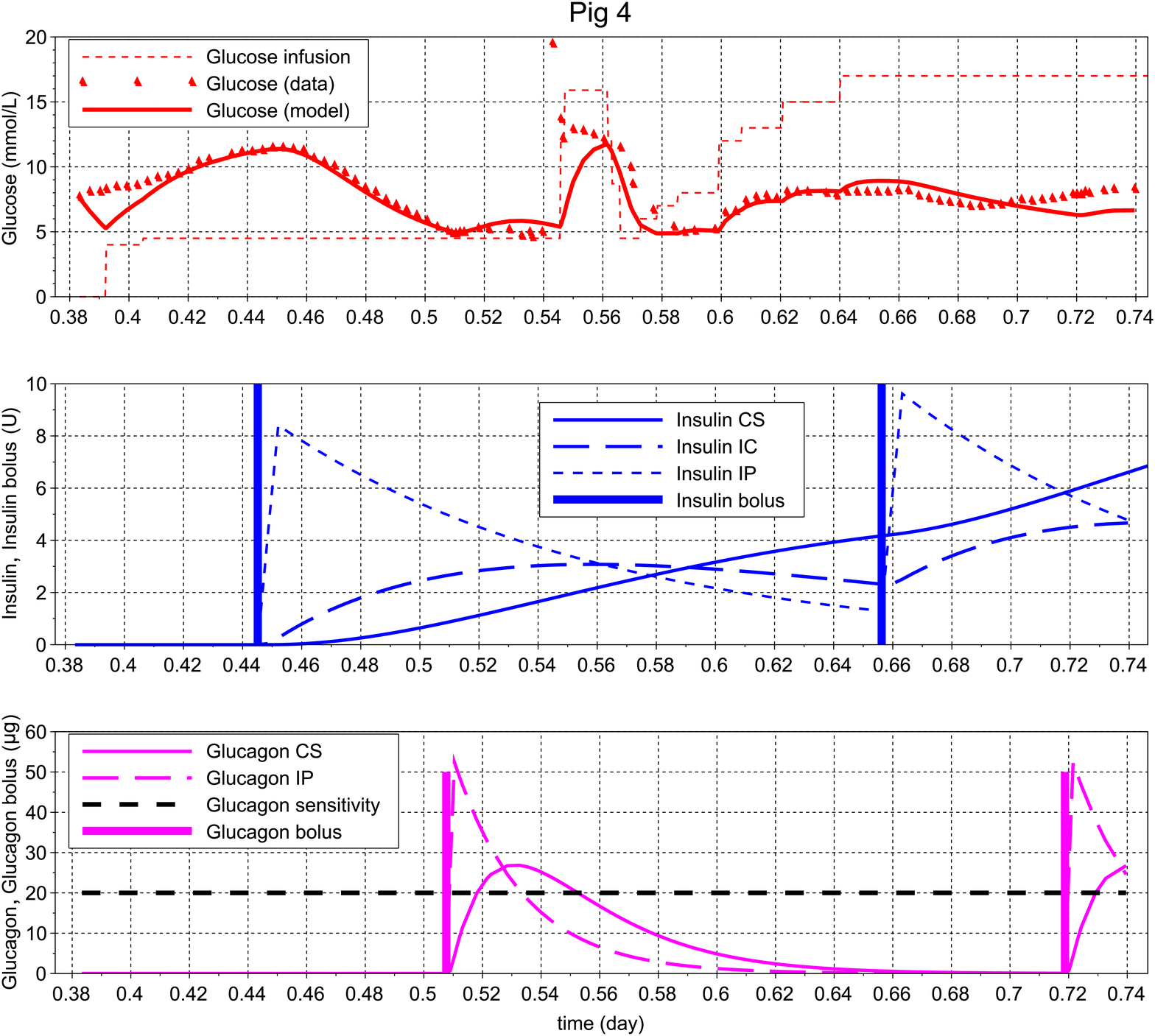}
\includegraphics[width=7.5cm]{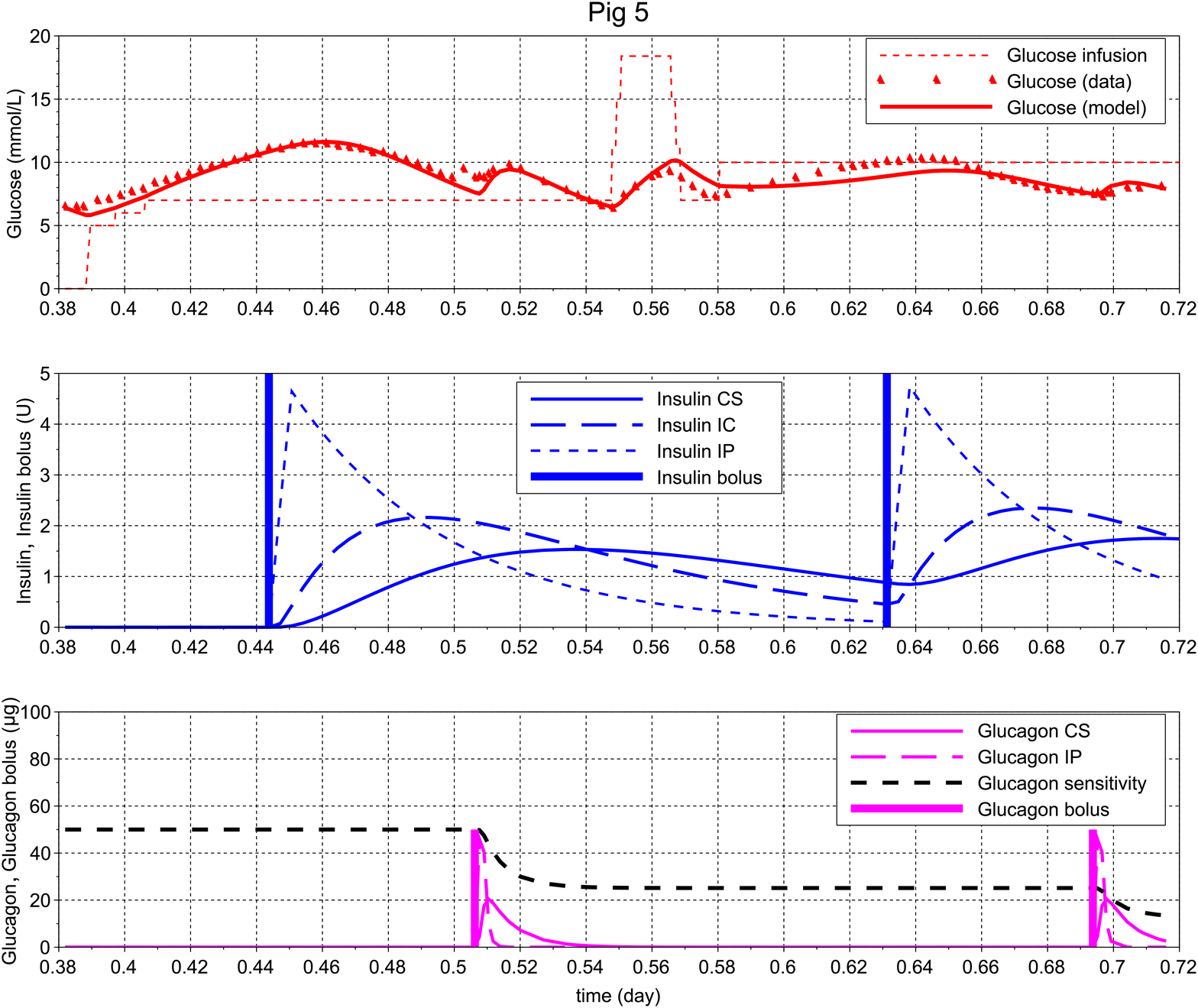}
\caption{Bihormonal--glucose model (1)--(7) 
compared with experimental data. All boluses of insulin and glucagon were introduced in the intraperitoneal cavity. The model parameters used for the numerical simulation are described in Table \ref{table:bihormonal}.}\label{fig:complete_3_4_5}
\end{figure*}

The results after calibrating the bihormonal--glucose model (1)--(7) 
 with experimental glucose measurements are presented in Figure \ref{fig:complete_3_4_5} and Table \ref{table:bihormonal}. The data was approximated using the minimization method described in Section V
 .

\gdef\thefigure{S.3}

\gdef\thetable{S.3}
\begin{table*}[!]
\begin{center}
\caption{Parameters estimated for the bihormonal--glucose model (1)--(7)
. The values were obtained minimizing the quadratic error between the model and the data. Results are depicted in Figure \ref{fig:complete_3_4_5}. MSE is the Mean Square Error of the model approximation and CV the parameter Coefficient of Variation (i.e. the standard deviation divided by the mean).}\label{table:bihormonal}
\begin{tabular}{ccccccc}
\hline
Parameter		&	Pig 1				&	Pig 2				&	Pig 3				&	Pig 4				&	Pig 5 & CV	\\	
\hline
$k_1$		&	8.24				&	22.29				&	8.21				&	41.54				&	14.93	&0.73\\

$k_I$		&	0.0001				&	1.05				&	0.023				&	22.21				&	62.66	&1.58\\

$k_{i_1}$	&	3.87				&	2.36				&	2.62				&	17.06				&	0.0638&1.30	\\

$k_H$		&	0.083				&	0.28				&	0.21				&	0.23				&	0.40	&0.48\\

$r_G$		&	3.73				&	2.80				&	2.41				&	11.20				&	3.33&0.78	\\

$m_1$		&	94.56				&	86.58				&	94.00				&	2.03				&	24.52	&0.73\\

$m_2$		&	13.97				&	19.31				&	211.92				&	11.80				&	12.71&1.64\\	

$m_3$		&	19.63				&	17.17				&	16.80				&	12.28				&	23.59	&0.23\\

$m_4$		&	76.92				&	83.49				&	70.65				&	9.17				&	20.75&0.66\\	

$p$			&	0.92				&	0.91				&	0.93				&	1.18				&	0.82&	0.14\\
		
$q$			&	0.57				&	0.60				&	0.61				&	0.77				&	1.10&0.30\\	
 			
$n$			&	149.15				&	175.38				&	156.11				&	49.42				&	119.63&0.38\\	

$n_1$		&	177.95				&	192.26				&	171.29				&	41.92				&	1002.72&1.22	\\

$n_2$		&	217.95				&	194.52				&	213.71				&	149.85				&	294.52&	0.24\\
$x_1$ 		&	1.42				&	1.04				&	0.32				&	0.00001				&	2.02&	0.85\\
$x_2$		&	0.00001				&	0.00001				&	0.000001			&	0.00105				&	0.00214	&1.48\\
\hline
MSE	&	0.63	&	0.69	&	0.41	&	3.92	&	0.35\\		
\hline
BIC			&	30.79				&	38.92				&	-2.00				&	234.33				&	-27.09	\\
\hline
\end{tabular}
\end{center}
\end{table*}

\gdef\thesection{S.4.1}
\section{Model Prediction}\label{appendix:prediction}

In this section, the bihormonal--glucose model (1)--(7) 
was calibrated using only the data from the first 3 hours and a part, which is a time interval that comprises 2 intraperitoneal boluses of each hormone (insulin and glucagon). Then, using the parameters obtained from the calibration, the approximation of the data on the second part of the experiment and in the entire experiment was evaluated with the BIC value. 

\begin{figure}[!]
\includegraphics[width=10cm]{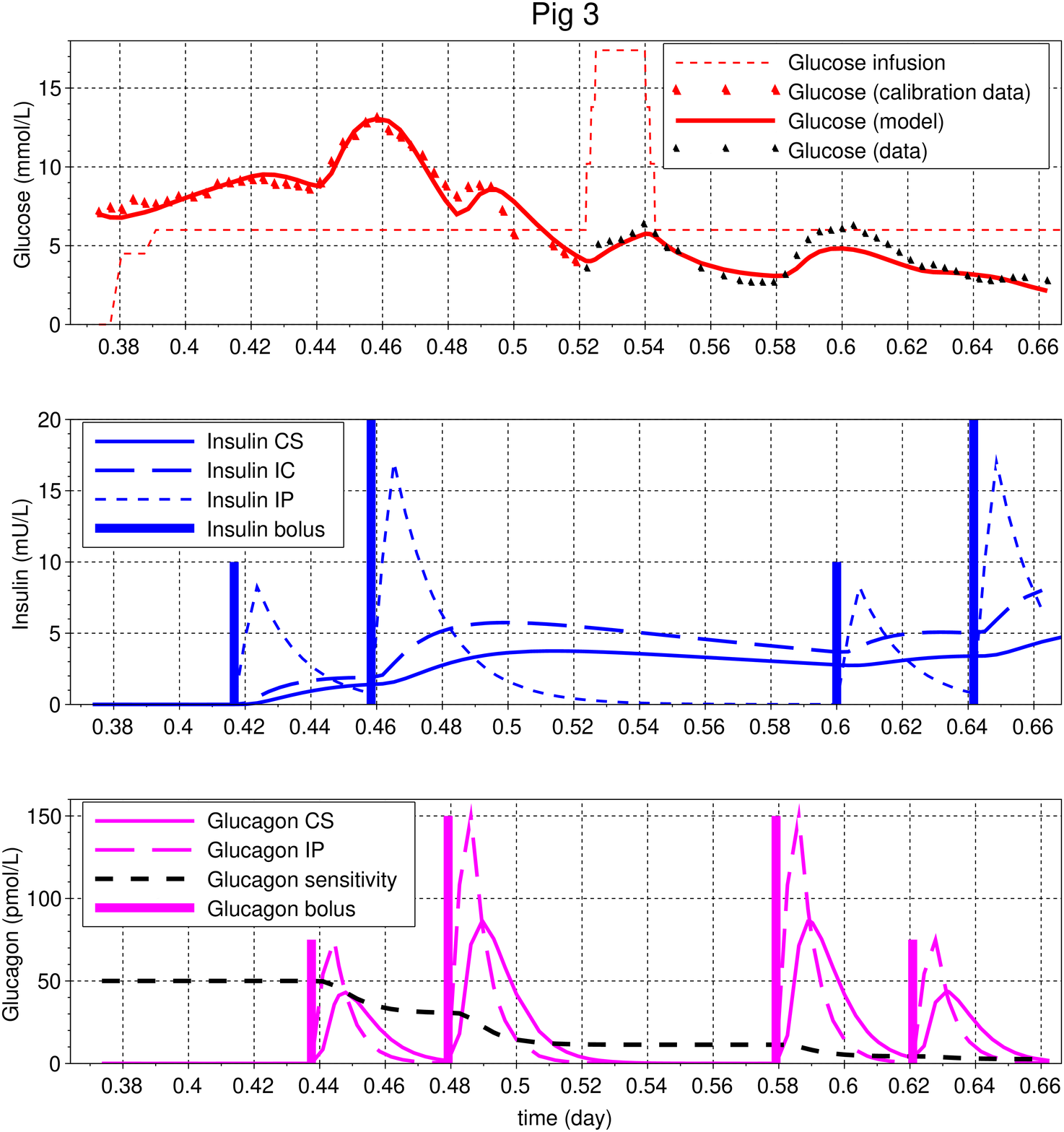}
\caption{Approximation of glucose data when the bihormonal--glucose model (1)--(7) 
is only calibrated with the data from the first part of the experiment (red triangles). The parameters used are described in Table \ref{table:prediction3}.}\label{fig:prediction3}
\end{figure}

\begin{figure}[!]
\includegraphics[width=10cm]{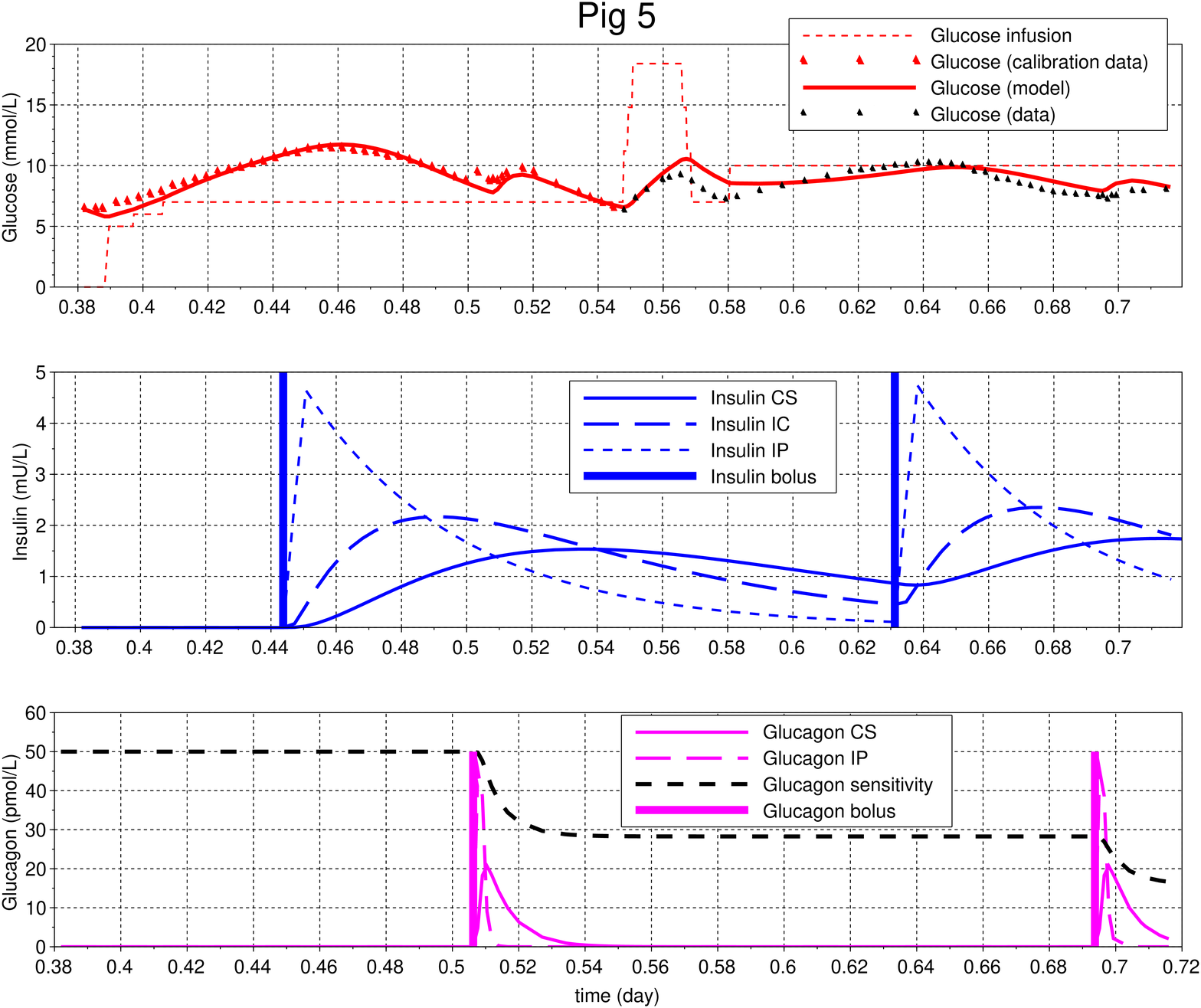}
\caption{Approximation of glucose data when the bihormonal--glucose model (1)--(7) 
is only calibrated with the data from the first part of the experiment (red triangles). The parameters used are described in Table \ref{table:prediction3}.}\label{fig:prediction5}
\end{figure}

The BIC value comparing the data from all the experiment and the approximation given by the model is negative, which indicates a good performance. Moreover, when evaluating the performance of the model to approximate the data on the first part of the experiment (where it was calibrated) and the second part there are obtained similar BIC values. This suggests that, when the model is calibrated with the first glucose measurements, it is then able to predict glucose levels of the next 4 hours with a similar accuracy. See Figures \ref{fig:prediction3}--\ref{fig:prediction5} for a graphical representation and Tables \ref{table:prediction3}--\ref{table:prediction5} where BIC values obtained are set.

\gdef\thefigure{S.4}

\gdef\thetable{S.4}
\begin{table*}[!]
\begin{center}
\caption{Parameters estimated using the Nelder--Mead algorithm to fit all the data and the first part of the experiment data (see Figure \ref{fig:prediction3}). The Error percent corresponds to the difference between the two parameters values obtained.}\label{table:prediction3}
\begin{tabular}{cccc}
\hline
&&Pig 3&\\
Parameter		&	Fit all data		&	Fit first part &	Error percent		
\\
\hline
$k_1$			&	8.21	&	8.50	&	3.42		
\\
$k_I$			&	0.023	&	0.025	&	{10.11}	
\\
$k_{i_1}$		&	2.62	&	2.41	&	7.96		
\\
$k_H$			&	0.212	&	0.202	&	4.65	
\\
$r_G$			&	2.41	&	2.37	&	1.41		
\\
$m_1$			&	94.00	&	96.72	&	2.89		
\\
$m_2$			&	211.92	&	233.44	&	{10.16}	
\\
$m_3$			&	16.80&	16.68	&	0.76		
\\
$m_4$			&	70.65	&	68.02	&	3.71		
\\
$p$				&	0.932	&	0.829	&	{10.99}	
\\
$q$				&	0.613	&	0.605	&	1.22		
\\
$n$				&	156.11	&	152.23	&	2.49		
\\
$n_1$			&	171.29	&	167.97	&	1.94		
\\
$n_2$			&	213.71	&	223.28	&	4.47		
\\
$x_1$ 			&	0.320	&	0.435	&	{35.75}	
\\
$x_2$			&	0.0000007	&	0.0000008	&	{14.29}	
\\
\hline
BIC total			&	-2.00		&	-3.13	&	--			
\\
BIC first part		&	41.28		&	21.49	&	--			
\\
BIC second part	&	-13.60		&	23.11	&	--			
\\
\hline
\end{tabular}
\end{center}
\end{table*}

\gdef\thefigure{S.5}

\gdef\thetable{S.5}
\begin{table*}[!]
\begin{center}
\caption{Parameters estimated using the Nelder--Mead algorithm to fit all the data and the first part of the experiment data (see Figure \ref{fig:prediction5}). The Error percent corresponds to the difference between the two parameters values obtained.}\label{table:prediction5}
\begin{tabular}{cccc}
\hline
&&Pig 5&\\
Parameter		&		Fit all data	&	Fit first part &	Error percent	
\\
\hline
$k_1$			&	14.93	&	15.29	&	2.43		
\\
$k_I$			&	62.66&	61.28 &	2.19		
\\
$k_{i_1}$		&	0.064	&	0.065&	1.86		
\\
$k_H$			&	0.401	&	0.348	&	{13.11}	
\\
$r_G$			&	3.33 &	3.42&	{2.72}	
\\
$m_1$			&	24.52	&	25.07	&	2.22		
\\
$m_2$			&	12.71	&	12.57	&	1.07		
\\
$m_3$			&	23.59&	23.91	&	1.34		
\\
$m_4$			&	20.75	&	20.78	&	0.14		
\\
$p$				&	0.822&	0.821		&	0.078		
\\
$q$				&	1.097	&	1.104	&	0.62		
\\
$n$				&	119.63		&	133.30	&	{11.43}
\\
$n_1$			&	1002.72	&	1045.28	&	{4.24}	
\\
$n_2$			&	294.52	&	295.86	&	0.46		
\\
$x_1$	 		&	2.02	&	1.93	&	{4.15}	
\\
$x_2$			&	0.002139	&	0.002133	&	0.28		
\\
\hline
BIC total			&	-27.10	&	-13.51	&	--			
\\
BIC first part		&	-7.05	&	-16.07	&	--			
\\
BIC second part	&	25.42	&	38.62	&	--			
\\
\hline
\end{tabular}
\end{center}
\end{table*}

\gdef\thesection{S.5}
\section{Analysis of Local Practical\\ Identifiability: Profile Likelihoods}\label{sec:PL}

Profile likelihoods are used in Systems Biology for addressing parameter identifiability and assessing confidence intervals for the parameters. The profile likelihood of a single parameter represents statistically the goodness of data fit for each value that the parameter can take in its confidence interval \cite{kreutz2013profile}. 

In this work, profile likelihoods were computed for the bihormonal--glucose model (1)--(7) 
parameters to analyze their practical identifiability. For this, it is assumed that glucose measurements have white Gaussian noise $\varepsilon\sim N(0,\sigma^2)$. Then, minus two times the log--likelihood (-2LL) is computed as the least square estimation
\begin{equation}\notag
\chi^2=\arg \min_{\boldsymbol\theta}\sum_i\big[ BGA(t_i) - G\big(t_i,u,\boldsymbol\theta\big) \big]^2/\sigma^2,
\end{equation}
where $G$ represents blood glucose concentration from (1)
, $BGA(t_i)$ blood glucose measurements, $u$ the system input ($Ra_G$, $u_I$ and $u_G$), $\boldsymbol\theta$ the vector of parameters and $\sigma^2$ the variance of the measurements. Indeed, in typical Systems Biology applications likelihood is equivalent to the least--squares criterion \cite{raue2009structural, kreutz2013profile}. Moreover, model glucose sensor errors have been proposed in the simplest way as white Gaussian noise processes \cite{facchinetti2010modeling}. 

The profile likelihood of a parameter $\boldsymbol\theta_j$ evaluated at $p$ corresponds to \cite{kreutz2013profile}
\begin{equation}\notag
PL_j(p)=\max_{\boldsymbol\theta\in\{\boldsymbol\theta\vert\boldsymbol\theta_j=p\}}LL(y\vert\boldsymbol\theta).
\end{equation}

 Then, the profile likelihood 
\begin{equation}
-2PL_j(p)=\arg \min_{\{\boldsymbol\theta:\boldsymbol\theta_j=p\}}\sum_i\frac{\big[BGA(t_i)-G(t_i,u,\boldsymbol\theta)\big]^2}{\sigma^2}
\end{equation}must be computed for each value $p$ in an interval of possible values of $\boldsymbol\theta_j$.

The cost function -2PL can possess multiple local minima, which complicates the reckoning of profile likelihoods \cite{schon2011system} and makes their numerical computation demanding. Here, profile likelihoods were computed numerically using the minimization tool \emph{fminsearch} in Scilab 6.0.2, which computes the minimum of the cost function with the Nelder--Mead algorithm (www.scilab.org). 
The initial parameters used for the algorithm were those found to approximate data with the bihormonal--glucose model (1)--(7) 
(see \ref{app_calibration}). 
However, the result of the minimization algorithm relies on the initial values given for the parameters and the number of iterations made. For these reasons, the numerical results here presented may be considered as local results.

Finally, in order to obtain a confidence threshold $\Delta(\alpha)$, the $\alpha$--quantiles of the -2PL distribution were computed with the \emph{perctl} function in Scilab. Then, the threshold with $\alpha$--$\%$ confidence is the $\alpha$--th percentile of the -2PL distribution.

\gdef\thefigure{S.6}
\begin{figure}[!]\centering
\includegraphics[width=9cm]{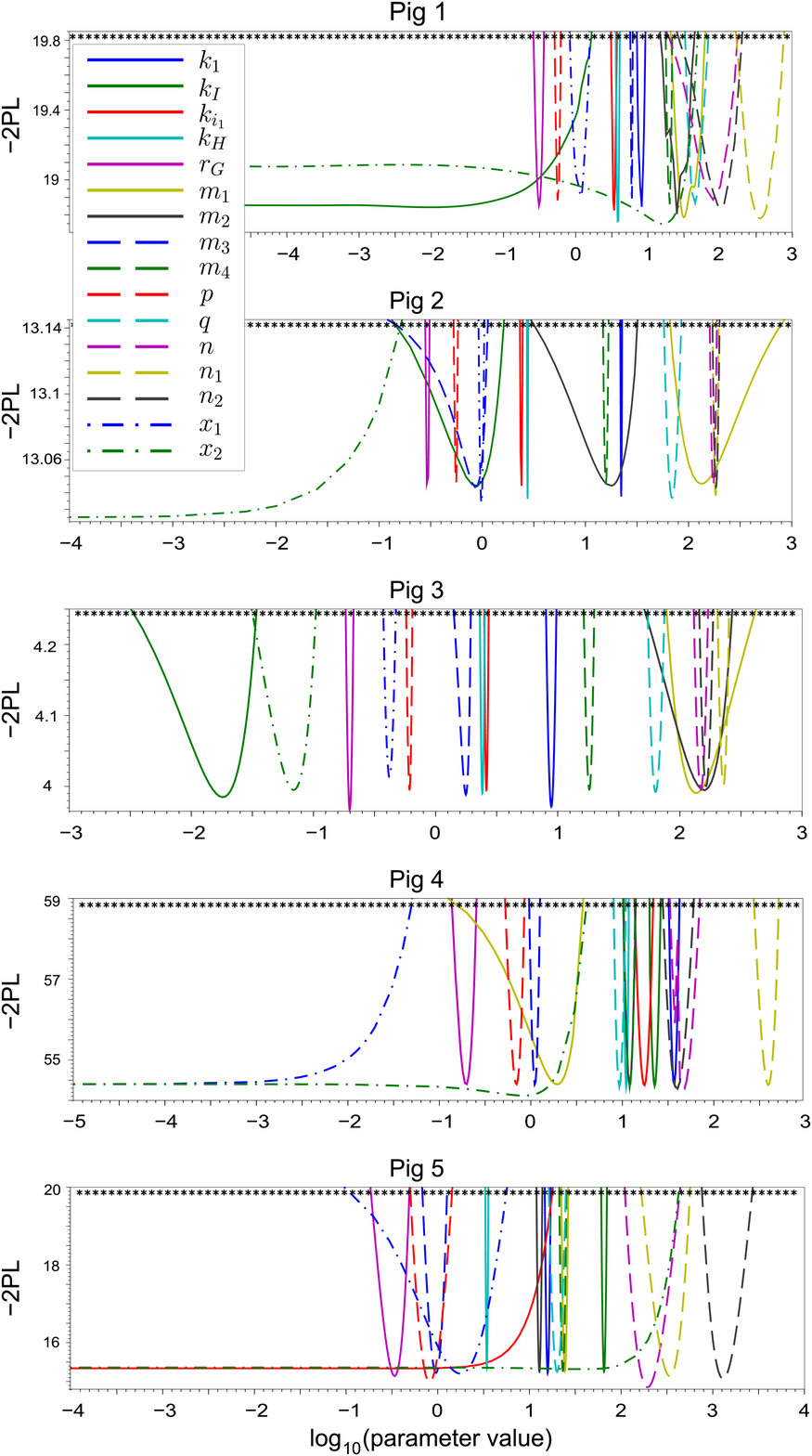}
\caption{Parameter profile likelihoods of the bihormonal--glucose model (1)--(7)
. Parameters $k_I$, $k_{i_1}$, $x_1$ and $x_2$ appear to be practically non--identifiable, because their profile likelihoods flatten towards small values in some cases. The rest of the parameters are locally identifiable, since they have profile likelihoods with single minima and that exceed the confidence threshold twice, delimitating confidence intervals of finite size. Asterisks denote the threshold $\Delta(\alpha)$ with confidence 98\% for pigs 1, 2 and 4, and 99\% for pigs 3 and 5.}\label{fig:PL} 
\end{figure}

\gdef\thesection{S.5.1}
\section{Results}

Profile likelihoods for the bihormonal--glucose model (1)--(7) 
 are depicted in Figure \ref{fig:PL}. According to this, the parameters $k_I$, $k_{i_1}$, $x_1$ and $x_2$ are practically non--identifiable, because in some cases they have profile likelihoods that flatten towards small values. The rest of the parameters are locally identifiable, since their profiles have unique minima and exceed the confidence thresholds giving confidence intervals of finite size \cite{raue2009structural, kreutz2013profile}.

The lack of practical identifiability 
of $k_I$ and $k_{i_1}$ can be explained as the infeasibility of identifying independently the rates of insulin--dependent glucose removal. This leads to conclude that, even assuming that glucose consumption can occur in different compartments, the rates at which this occurs in each compartment cannot be distinguished just from the blood glucose measurements of the experiments.

On the other hand, $x_1$ appears to be non--identifiable only for one experiment where glucose response to glucagon boluses was apparently small and glucagon sensitivity remains almost constant (see Figure \ref{fig:complete_3_4_5}). The non--identifiability of $x_2$ may be consequence of the apparent lack of glucagon sensitivity restoration, since the glucose response to glucagon boluses was almost always decreasing in these specific animal trials (see simulations and data in  \ref{app_calibration}). The possible reason for this unexpected behavior is discussed by \cite{am2020intraperitoneal}.

\gdef\thesection{S.6}
\section{Local Structural Identifiability: \\Sensitivity Analysis}\label{sec:SVD}

In this section, the local structural identifiability of the bihormonal--glucose model (1)--(7) 
 is addressed. The purpose is to check if the  lack of practical identifiability 
 detected in the previous Section is related to the model structure or just to the experimental design.
Practical identifiability considers the amount and quality of measured data that was used for parameter calibration. Structural identifiability is related to the model structure independently of experimental data \cite{raue2009structural}.

The method applied to assess the local structural identifiability and the possible correlations between model parameters is based on the SVD of the sensitivity matrix and it has been performed following the works in \cite{stigter2015fast} and \cite{staal2019glucose}.
 
 Given $(N+1)$ equidistant time--points $[t_0, t_1,\dots,t_N]$, the first step is determining the parameter--to--output sensitivity matrix
 \begin{equation}\notag
 S(\boldsymbol\theta)=\begin{pmatrix}
 \frac{\partial Y_1(t_0)}{\partial\boldsymbol\theta_1} & \dots & \frac{\partial Y_1(t_0)}{\partial\boldsymbol\theta_m} \\
 \vdots & \ddots & \vdots \\
 \frac{\partial Y_n(t_0)}{\partial\boldsymbol\theta_1} & \dots & \frac{\partial Y_n(t_0)}{\partial\boldsymbol\theta_m} \\
  \vdots & \ddots & \vdots \\
 \frac{\partial Y_1(t_N)}{\partial\boldsymbol\theta_1} & \dots & \frac{\partial Y_1(t_N)}{\partial\boldsymbol\theta_m} \\
 \vdots & \ddots & \vdots \\
 \frac{\partial Y_n(t_N)}{\partial\boldsymbol\theta_1} & \dots & \frac{\partial Y_n(t_N)}{\partial\boldsymbol\theta_m} \\
 \end{pmatrix},
 \end{equation}
 where $Y_1$, $Y_2$,...,$Y_n$ are the output states and $\boldsymbol\theta_1$, $\boldsymbol\theta_2$,..., $\boldsymbol\theta_m$ the model parameters. 
Having the state and output equations
 \begin{align}\notag
 \frac{dX}{dt}&=f(t,X,u,\boldsymbol\theta),\\\notag
 Y&=\bold H(X,\boldsymbol\theta),
 \end{align}the parameter--to--output sensitivities can be derived from the associated parametric state--to--output sensitivity system \cite{stigter2015fast}:
\begin{align}\notag
\frac{dX_{\boldsymbol\theta}}{dt}&=\frac{\partial f}{\partial X}X_{\boldsymbol\theta}+\frac{\partial f}{\partial\boldsymbol\theta}\\\notag
Y_{\boldsymbol\theta}&=\frac{\partial \bold H}{\partial X}X_{\boldsymbol\theta}+\frac{\partial\bold H}{\partial\boldsymbol\theta},
\end{align}where $X_{\boldsymbol\theta}$ and $Y_{\boldsymbol\theta}$ denote the partial derivatives of the state and output vector.

In this work, the output equation consists only of the blood glucose state. Then, the output function is $\bold H(X)=G$. Given the simplicity of the bihormonal--glucose model (1)--(7)
, the partial derivatives $\frac{df}{dX}$ and $\frac{df}{d\boldsymbol\theta}$ can be computed analytically.  Then, the equation $\frac{dX_{\boldsymbol\theta}}{dt}$ was numerically solved.

The SVD of the parameter--to--output sensitivity matrix must be performed for different vectors of parameters $\boldsymbol\theta^i$. Parameters $\boldsymbol\theta^i$ must be chosen from realistic parameter values \cite{stigter2015fast}. For this, parameter vectors $\boldsymbol\theta^i$ were chosen considering the confidence intervals obtained in Section \ref{sec:PL}. For obtaining a uniform selection of parameters, values that fit glucose measurements (see  \ref{app_calibration}) and their translations of 1\% of the respective confidence interval were considered for the vectors $\boldsymbol\theta^i$. That is to say, for each parameter three values are considered: 1) the value in Table \ref{table:bihormonal}, 2) the value of Table \ref{table:bihormonal} plus 1/100 times the length of its confidence interval and 3) the value of Table \ref{table:bihormonal} minus 1/100 times the length of its confidence interval. Then, from the vector of parameter values reported in Table \ref{table:bihormonal}, new vectors are generated by changing one entry with the values translated 2) and 3) for each parameter.

For each parameter vector $\boldsymbol\theta^i$, the SVD of the parameter--to--output sensitivity matrix was computed with the tool \emph{svd} in Scilab. The distances between the points $t_i$ were set as 0.004 in the scale of days (d), equivalent to 5m46s approximately (during the experiments blood samples were taken at least every 5 minutes). All the equidistant points $t_0, t_1,\dots, t_N$ in the full experiment time--interval (about 8 hours) were taken.

\begin{figure}[!]\centering
\includegraphics[width=7.2cm]{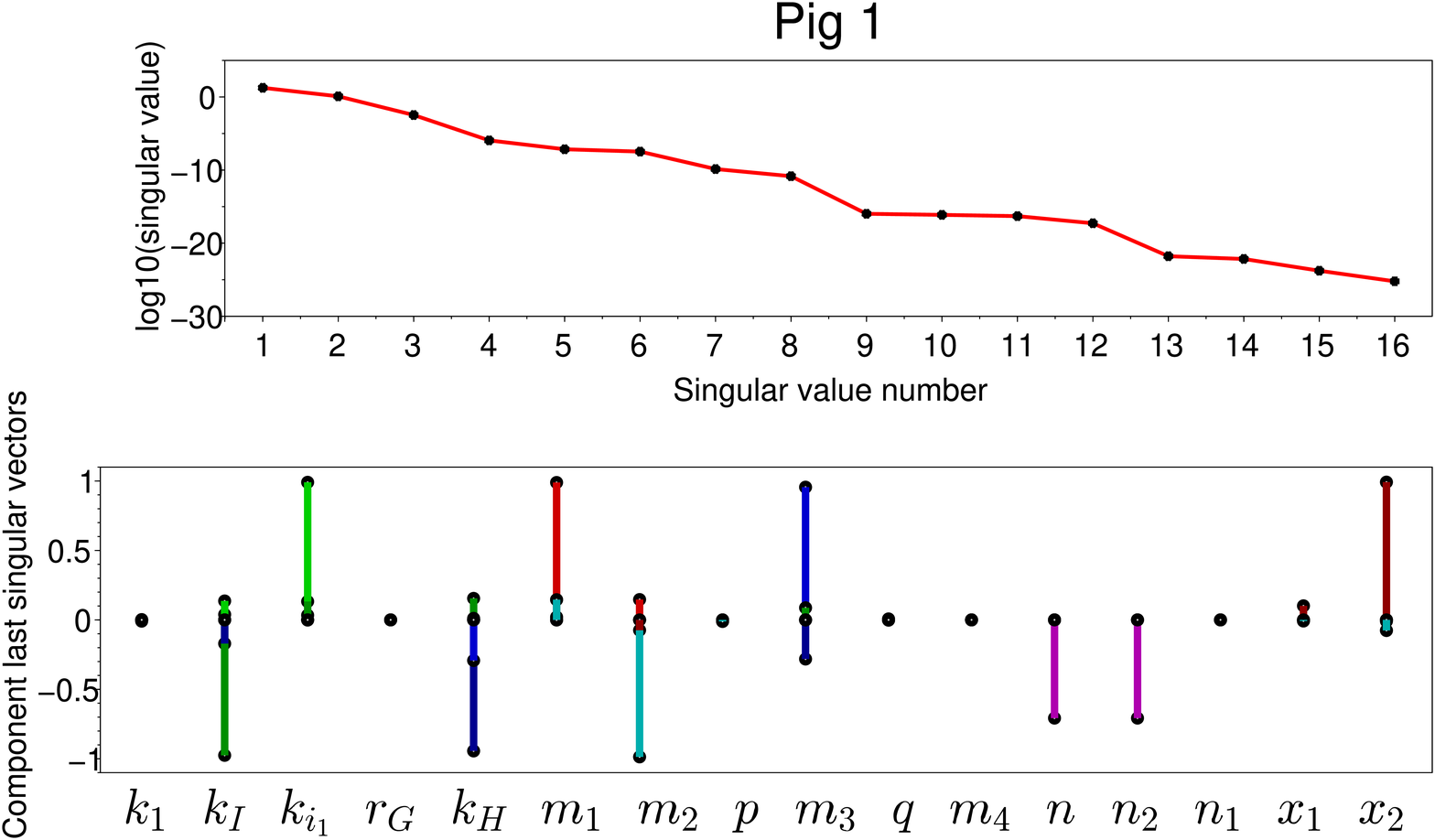}\\[.25cm]
\includegraphics[width=7.2cm]{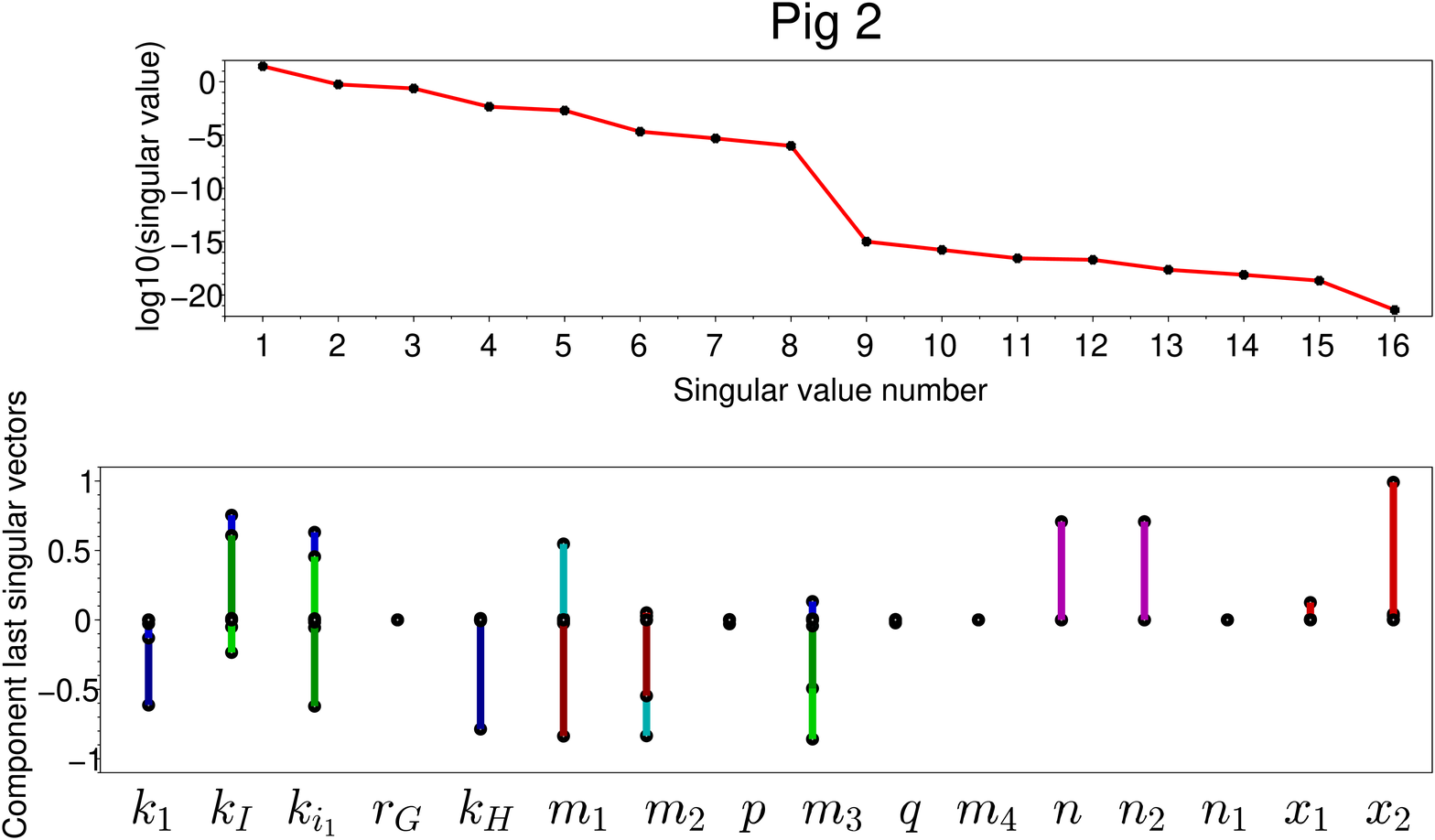}\\[.25cm]
\includegraphics[width=7.2cm]{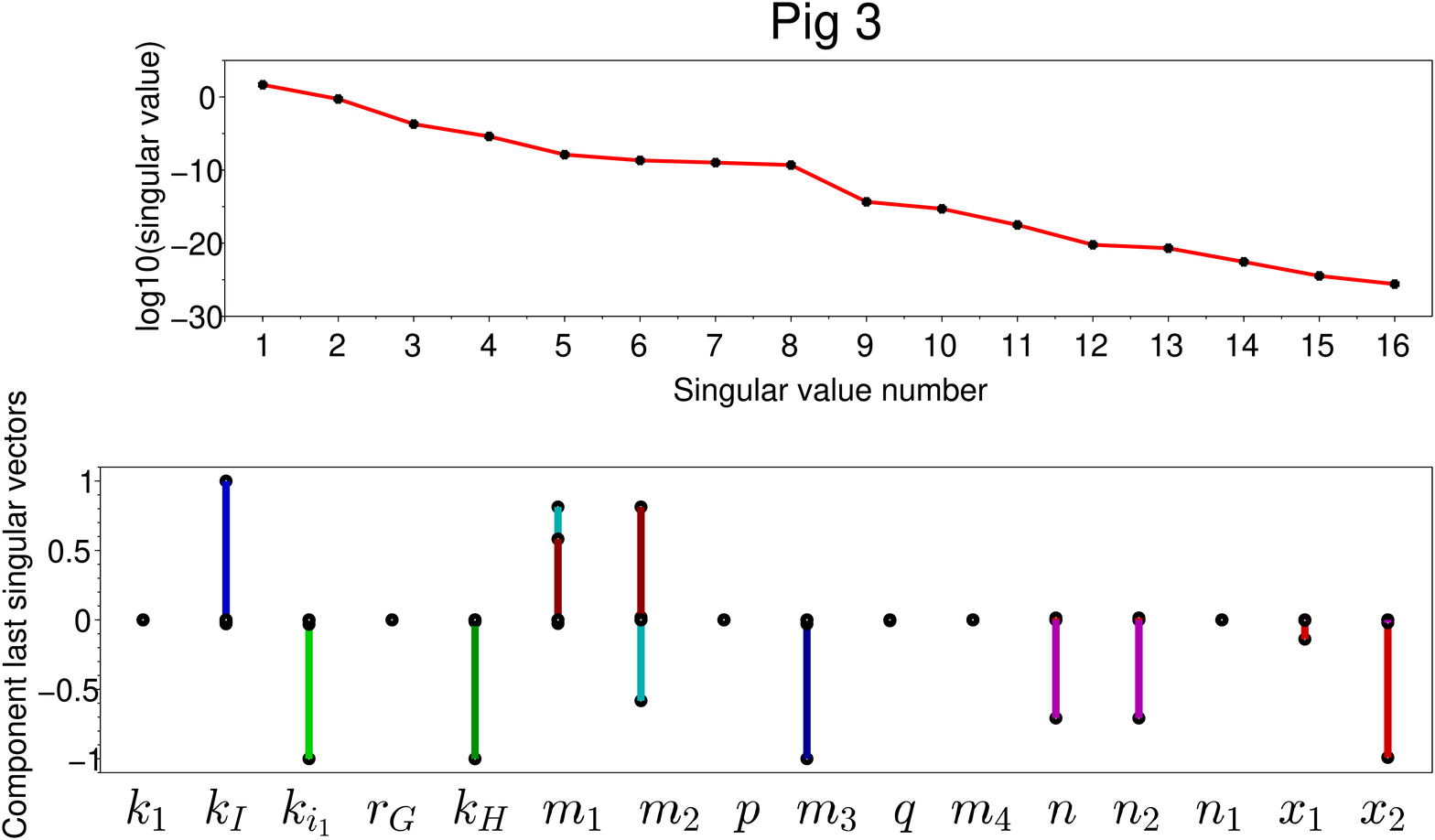}\\[.25cm]
\includegraphics[width=7.2cm]{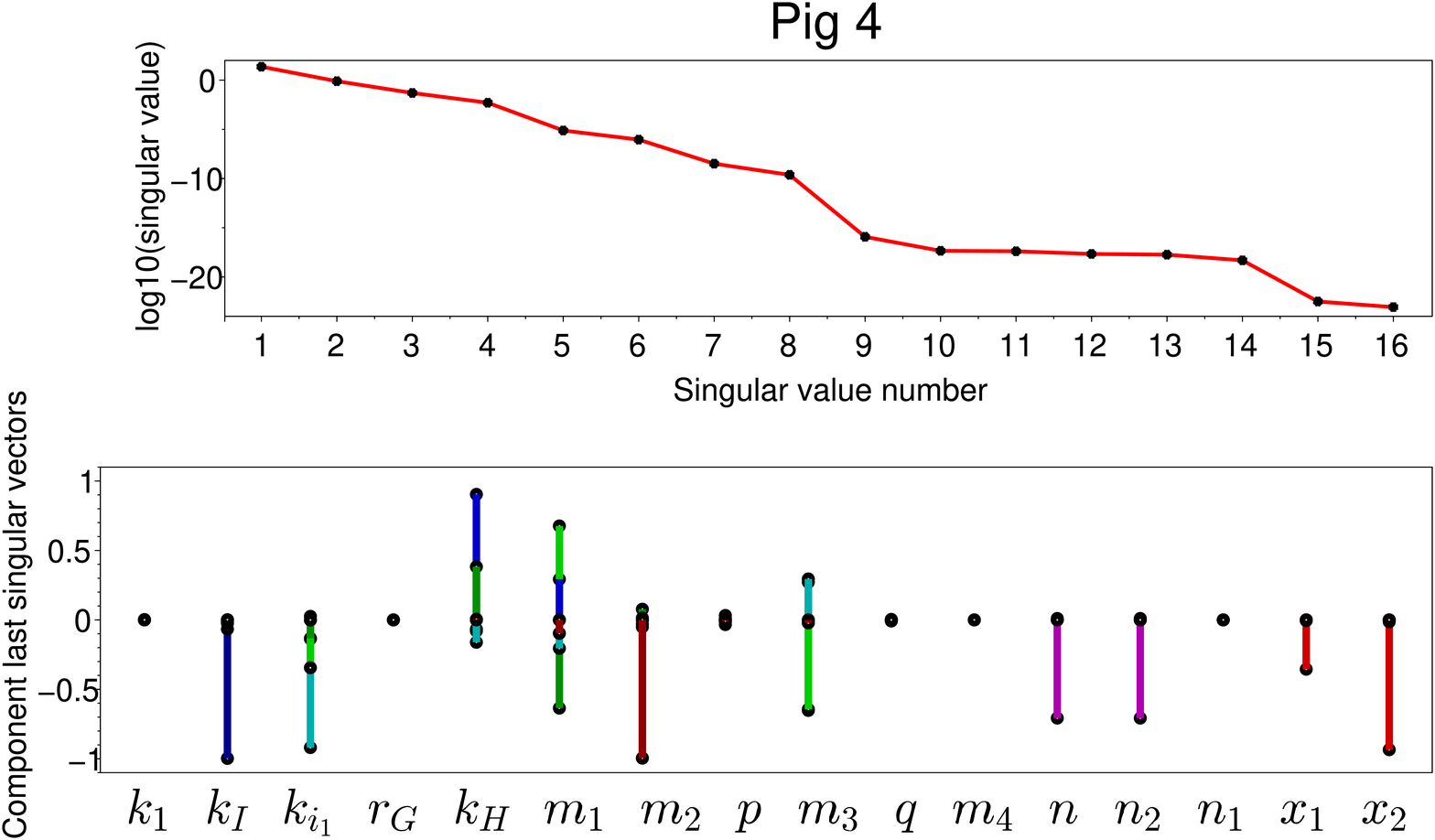}\\[.25cm]
\includegraphics[width=7.2cm]{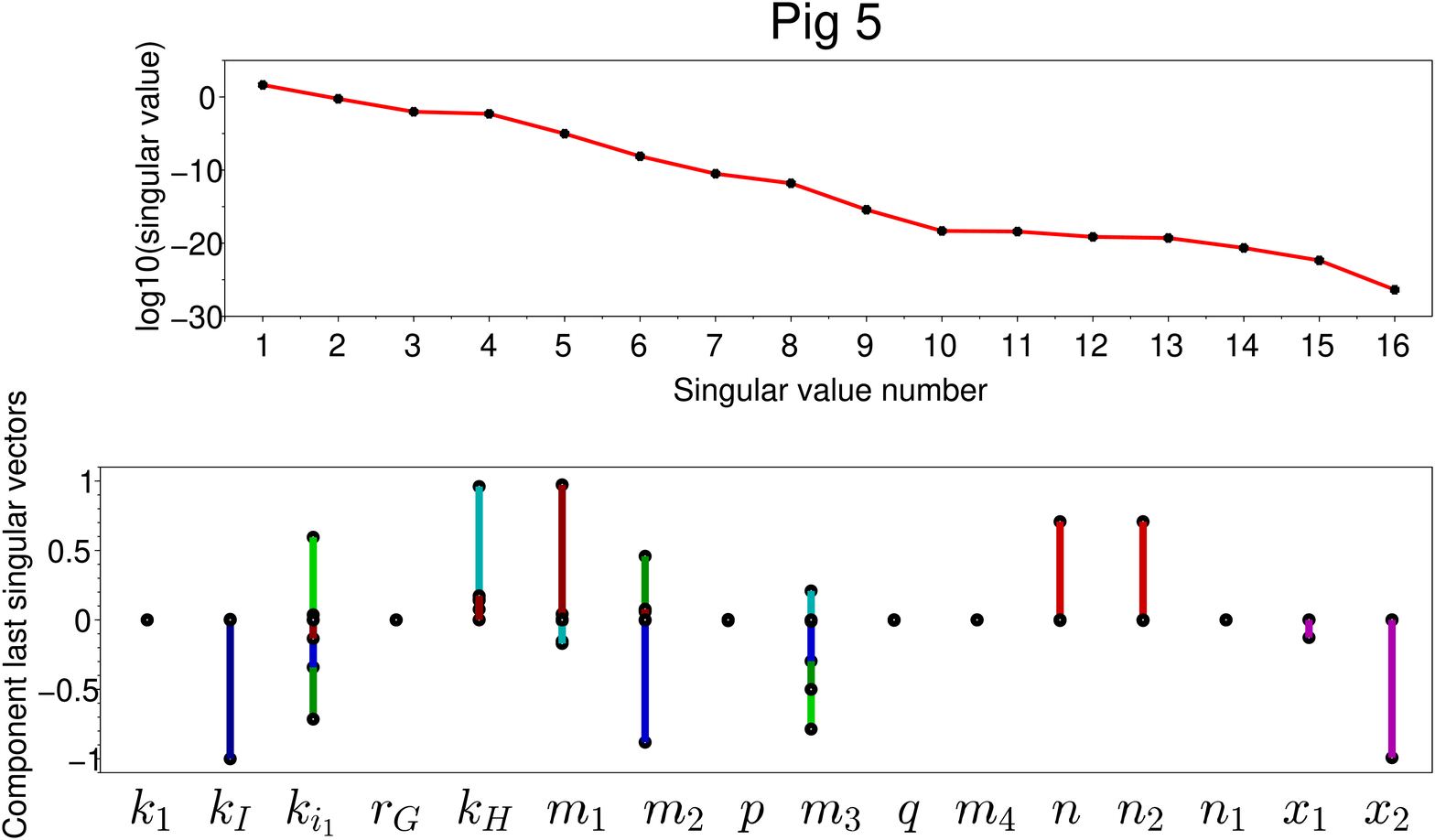}
\caption{Singular value decomposition (SVD) of the sensitivity matrix of the bihormonal--glucose model (1)--(7)
. Each color corresponds to one singular vector. Then, at each experimental case, correlated parameters are those with bars of the same color.}\label{fig:SVD}
\end{figure}

The singular vectors related to the small singular values were plotted to observe the non--identifiable parameters and their correlations \cite{stigter2015fast}. The results are depicted in Figure \ref{fig:SVD}.

A locally structurally identifiable re--parametrization can be computed in case of having a non--trivial  null--space of the Jacobi matrix \cite{stigter2015fast}
\begin{align*}
\frac{dJ}{d\boldsymbol\theta}(X(0),\boldsymbol\theta)=
\begin{pmatrix}
\frac{\partial\bold H}{\partial\boldsymbol\theta_1}&\dots&\frac{\partial\bold H}{\partial\boldsymbol\theta_m}\\
\frac{\partial\mathcal L_{f_0}\bold H}{\partial\boldsymbol\theta_1}&\dots&\frac{\partial\mathcal L_{f_0}\bold H}{\partial\boldsymbol\theta_m}\\
\frac{\partial\mathcal L_{f_0}\mathcal L_{f_1}\bold H}{\partial\boldsymbol\theta_1}&\dots&\frac{\partial\mathcal L_{f_0}\mathcal L_{f_1}\bold H}{\partial\boldsymbol\theta_m}\\
\vdots&\ddots&\vdots\\
\frac{\partial\mathcal L_{f_{j_0}}\dots\mathcal L_{f_{j_1}}\bold H}{\partial\boldsymbol\theta_1}&\dots&\frac{\partial\mathcal L_{f_{j_0}}\dots\mathcal L_{f_{j_1}}\bold H}{\partial\boldsymbol\theta_m}\\
\vdots&\ddots&\vdots
\end{pmatrix},
\end{align*}where $f_0$ represents the system with no control and $f_i$ the control vectors of the \emph{affine control system}
\begin{align*}
\frac{dX}{dt}&=f_0(X,\boldsymbol\theta)+\sum_{i=1}^r f_i(X,\boldsymbol\theta)\cdot u_i\\
 Y&=\bold H(X,\boldsymbol\theta),
\end{align*}and $\mathcal L_f \bold H$ denotes the Lie derivative of function $\bold H$ in the direction of the vector function $f$, i.e. \cite{walter1996identifiability, stigter2015fast}
\begin{align*}
\mathcal L_f \bold H=\frac{d\bold H}{dX}\cdot f.
\end{align*}

\gdef\thefigure{S.7}

\gdef\thesection{S.7}
\section{Singular Value Decomposition of the\\
 Bihormonal--Glucose Model}\label{SMfull}

\begin{figure}[!]\centering
\includegraphics[width=9cm]{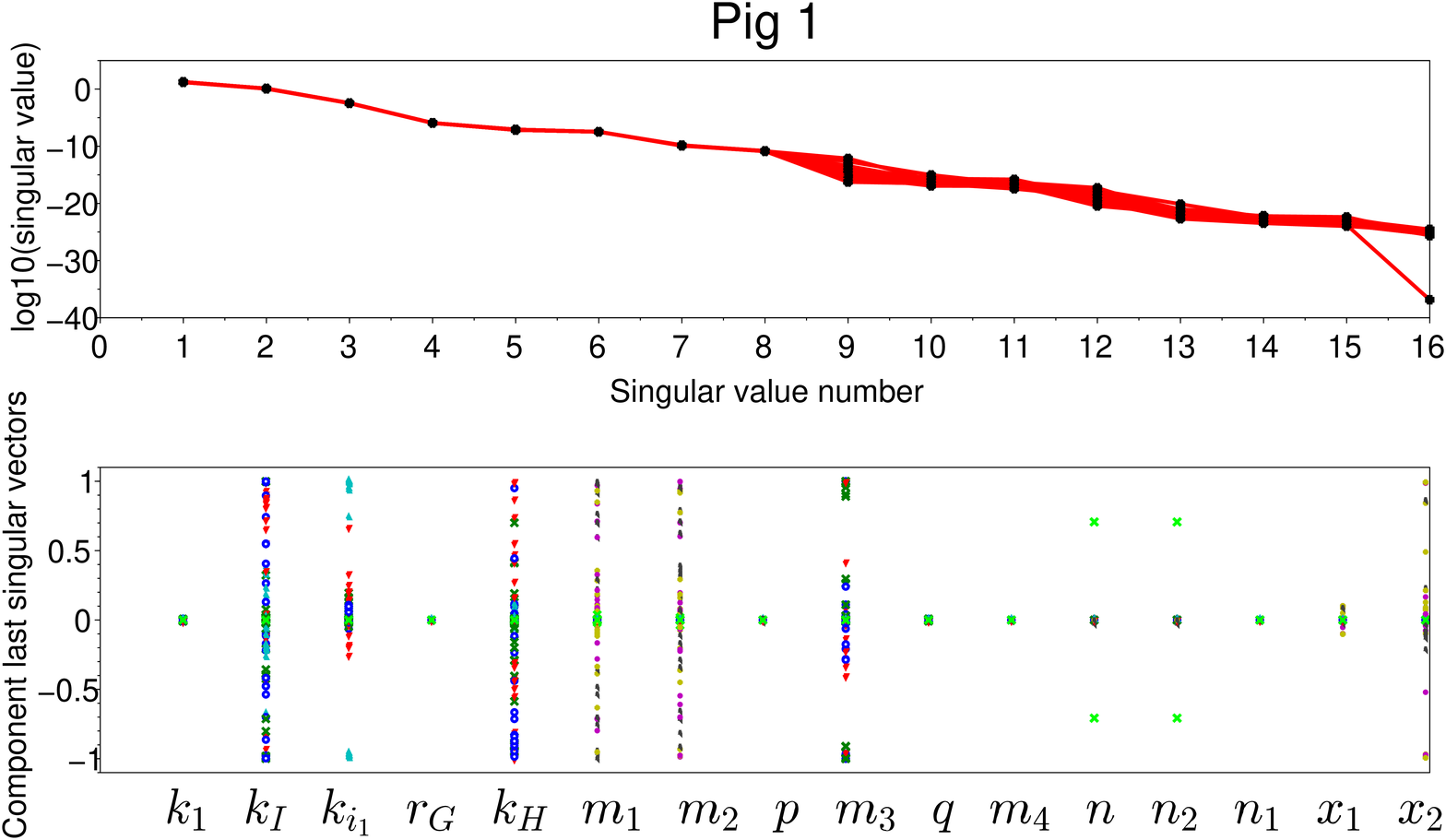}\\
\includegraphics[width=9cm]{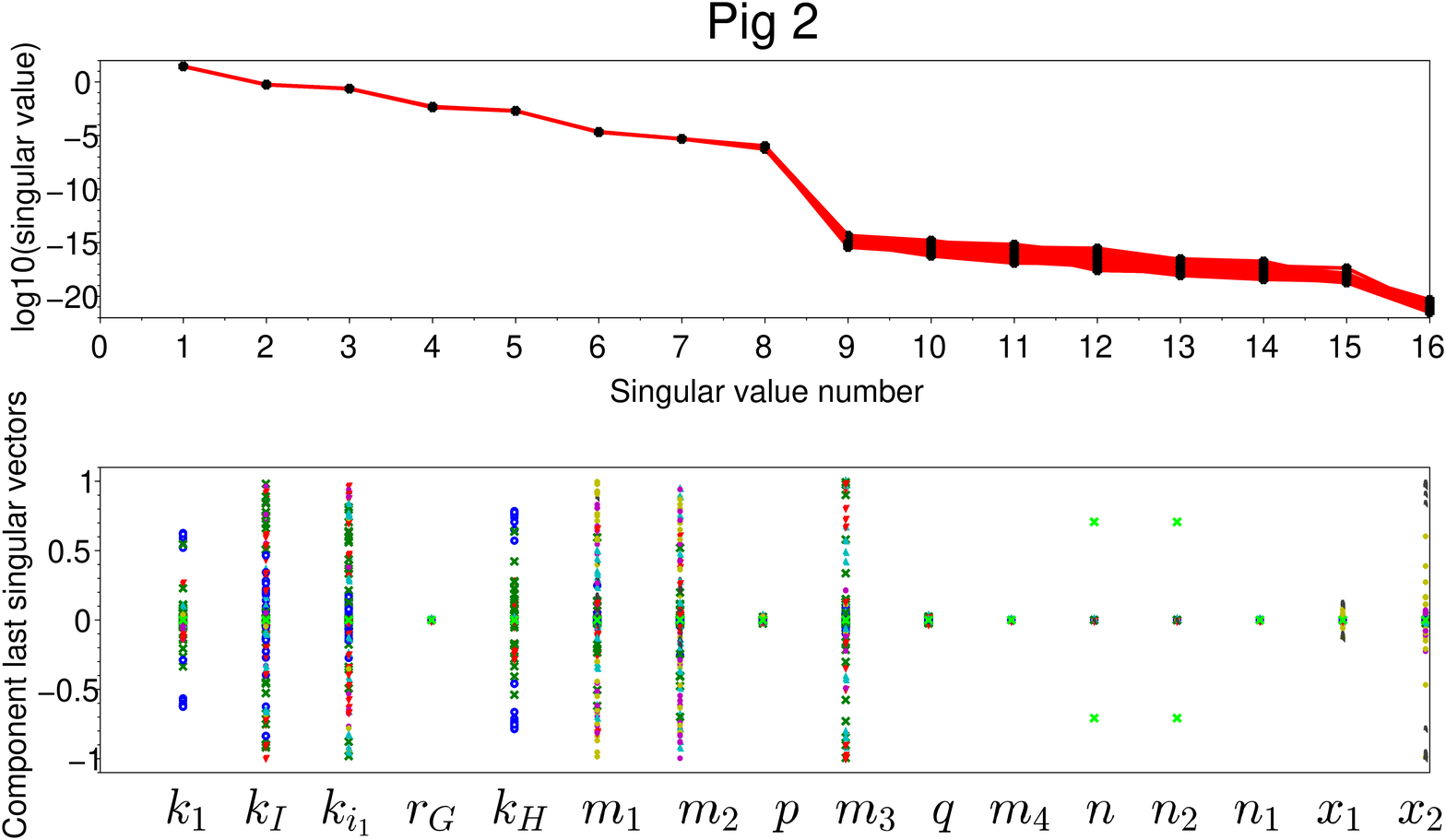}\\
\includegraphics[width=9cm]{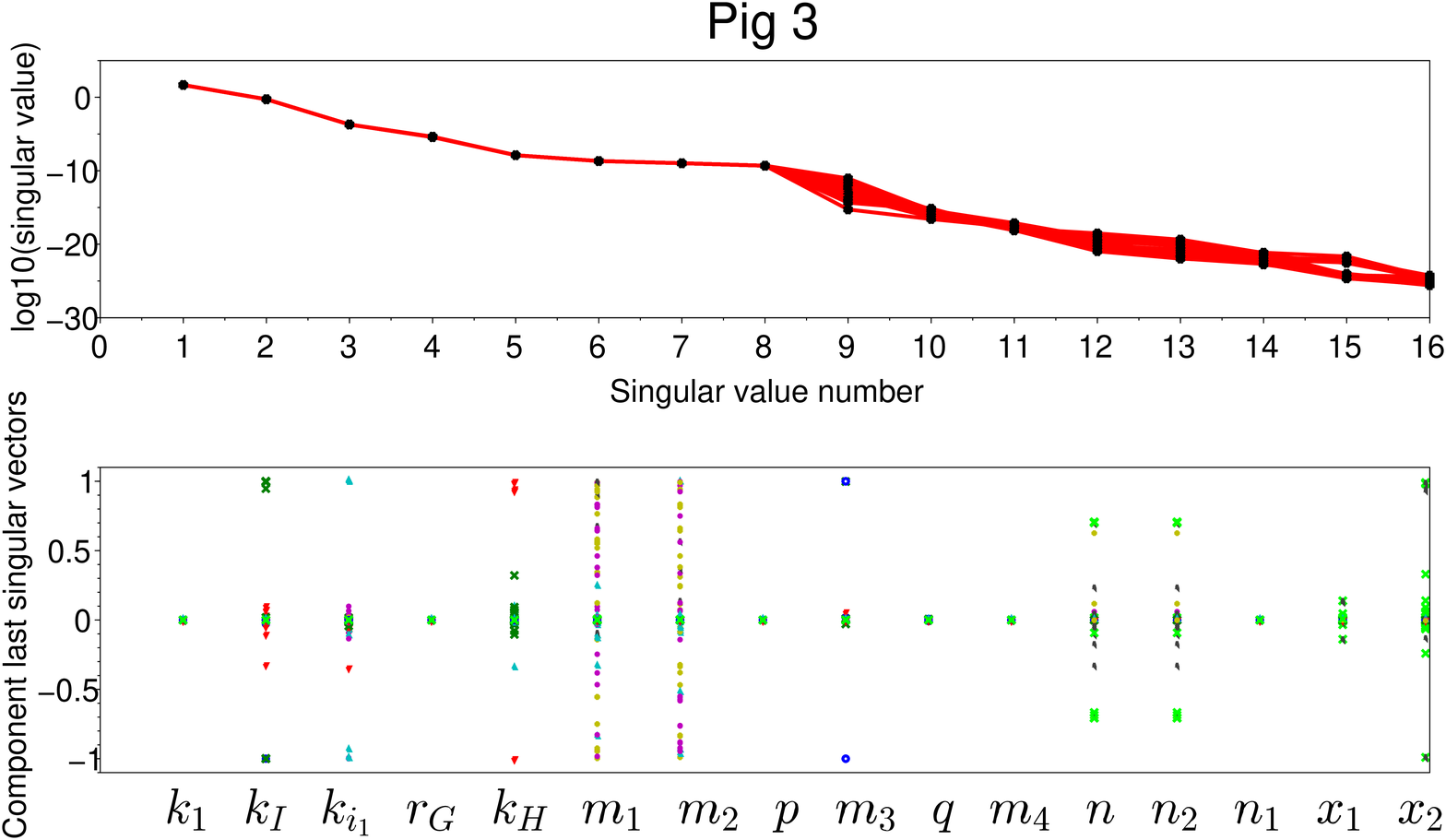}\\
\includegraphics[width=9cm]{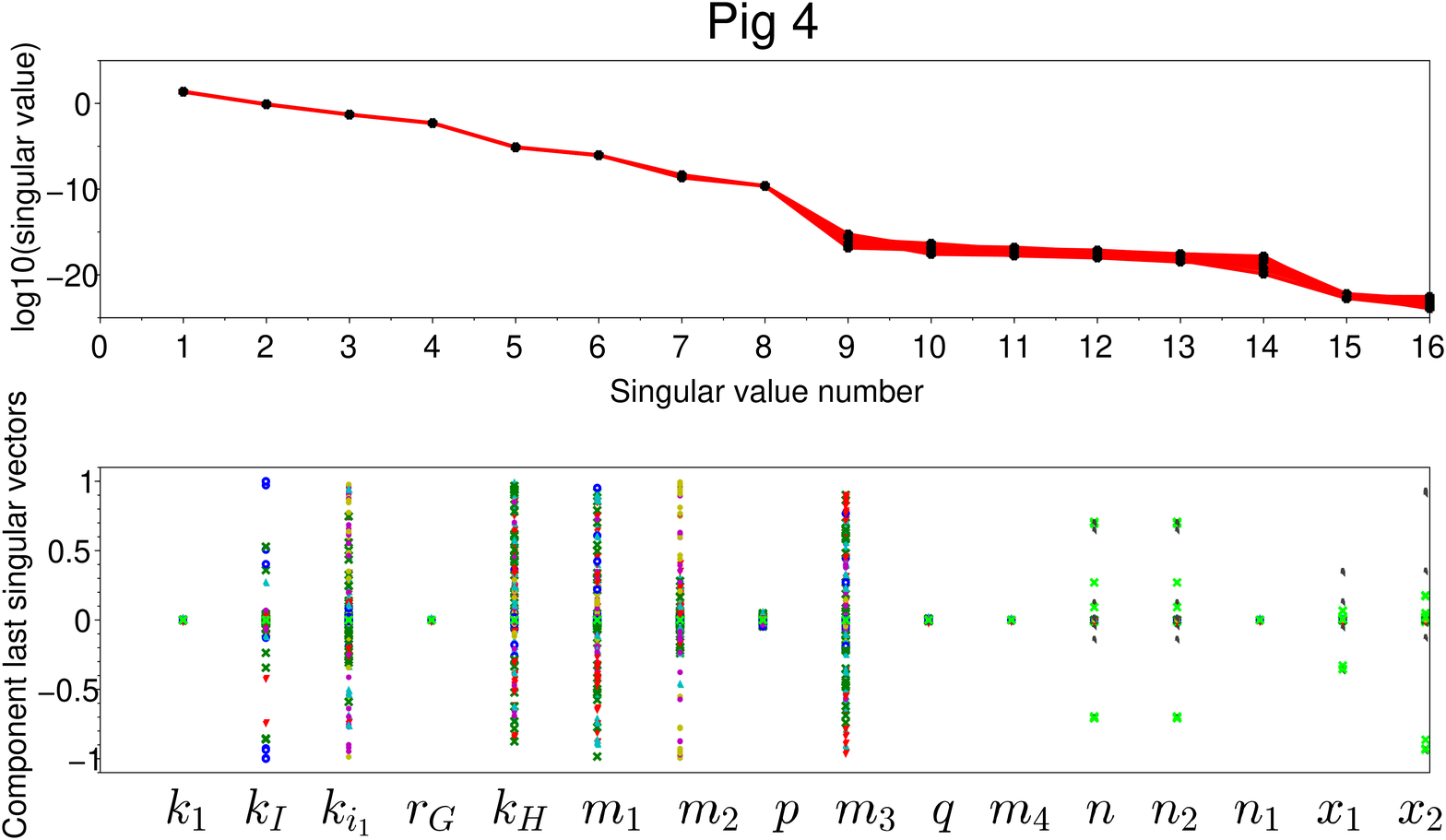}\\
\includegraphics[width=9cm]{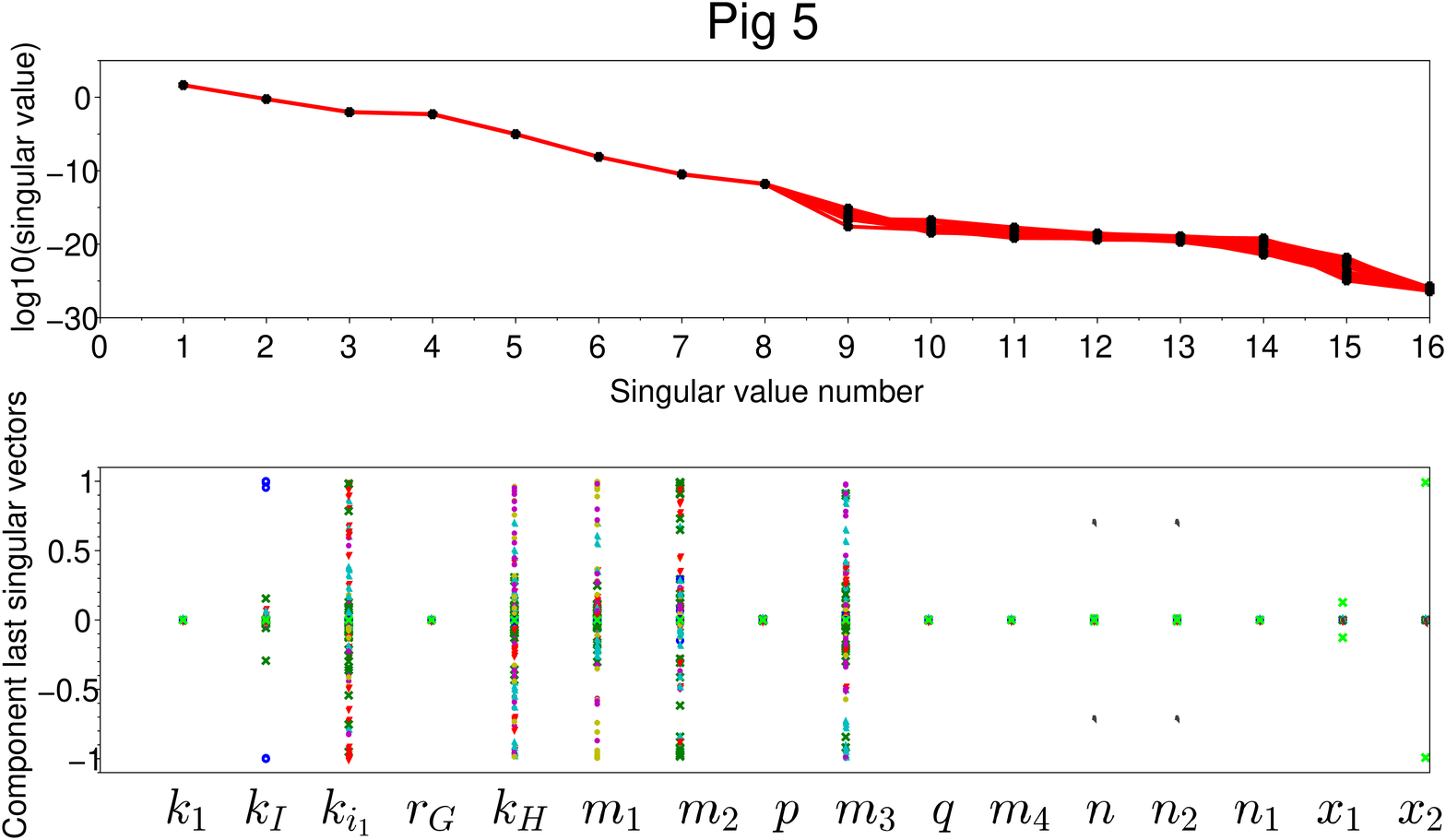}
\caption{Singular value decomposition (SVD) of the sensitivity matrix of the bihormonal--glucose model (1)--(7) 
 for each subject. Each symbol / color corresponds to one singular vector. Then, at each experimental case, correlated parameters are those with the same symbol / color.}\label{fig:SVD_full}
\end{figure}

In Figure \ref{fig:SVD_full} are depicted all the singular vectors related to the lowest singular values of the sensitivity matrix  of the bihormonal--glucose model (1)--(7)
. For each case, there are 8 singular values considered as very small and, according to each set of parameter $\boldsymbol\theta^i$, the singular vectors related to these values vary. This fact makes particularly difficult to identify the parameters that are correlated. However, in  \ref{SMcomb}, more clear results are obtained using the SVD of a bigger sensitivity matrix that combines several cases. 

{\color{black}In order to guarantee that a locally structurally identifiable re--parametrization can be found, a non--trivial null--space has been exhibited for the Jacobi matrix of Lie derivatives that only includes the parameters with lack of structural identifiability \cite{stigter2015fast}. In this case, the re--parametrization proposed can be carried-out without changing the output, which only consists of the glucose state, because the hormone states are not compared with experimental data and they are considered as dimensionless. }

\gdef\thefigure{S.8}

\gdef\thesection{S.7.1}
\section{Jacobi Rank Test}\label{app_Jacobi}
{\color{black}
For the bihormonal--glucose system (1)--(7)
, the \emph{affine control system} can be written defining the vector with no control as
\begin{align*}
&f_0(X,\boldsymbol\theta):=\\
&\begin{pmatrix}
-[k_1 + k_I\cdot ( I + I_b ) + k_{i_1}\cdot i_1]\cdot G + k_H\cdot( H +H_b )\cdot\xi	\\
- m_1\cdot I + m_2\cdot i_1^p 	\\
- m_3\cdot i_1^q+ m_4\cdot i_2 	\\	
- m_4\cdot  i_2 	\\	
- n\cdot H + n_2 \cdot h_1 	\\
- n_1\cdot h_1 \\
- x_1\cdot H\cdot \xi + x_2\cdot G\cdot I
\end{pmatrix},
\end{align*}
the glucose input by 
\begin{align*}
f_1(X,\boldsymbol\theta)&:=(r_G,0,0,0,0,0,0)^T\\
u_1&:=Ra_G,
\end{align*}
the insulin input by
\begin{align*}
f_2(X,\boldsymbol\theta)&:=(0,0,0,1,0,0,0)^T\\
u_1&:=u_I,
\end{align*}
and the glucagon input by
\begin{align*}
f_3(X,\boldsymbol\theta)&:=(0,0,0,0,0,1,0)^T\\
u_3&:=u_H.
\end{align*}

Furthermore, since $\bold H=G$, the Lie derivatives are
\begin{align*}
&\mathcal L_{f_0}\bold H=\\
&-[k_1 + k_I\cdot ( I + I_b ) + k_{i_1}\cdot i_1]\cdot G + k_H\cdot( H +H_b )\cdot\xi,	\\
&\mathcal L_{f_1}\bold H=r_G,\\
&\mathcal L_{f_2}\bold H=0,\\
&\mathcal L_{f_3}\bold H=0,\\
&\mathcal L_{f_0}\mathcal L_{f_0}\bold H=\\
&(-[k_1 + k_I\cdot ( I + I_b ) + k_{i_1}\cdot i_1]\cdot G + k_H\cdot( H +H_b )\cdot\xi)\\
&\cdot(-[k_1 + k_I\cdot ( I + I_b ) + k_{i_1}\cdot i_1])\\
&+(- m_1\cdot I + m_2\cdot i_1^p)\cdot(-k_IG)\\ 
&+ (- m_3\cdot i_1^q+ m_4\cdot i_2)\cdot(-k_{i_1}G)\\
&+ (- n\cdot H + n_2 \cdot h_1)\cdot(k_H\xi)\\ 
&+ (- x_1\cdot H\cdot \xi + x_2\cdot G\cdot I)\cdot(k_H\cdot( H +H_b) ),	\\
&\mathcal L_{f_1}\mathcal L_{f_0}\bold H=r_G\cdot(-[k_1 + k_I\cdot ( I + I_b ) + k_{i_1}\cdot i_1])\\
&\mathcal L_{f_2}\mathcal L_{f_0}\bold H=0\\
&\mathcal L_{f_3}\mathcal L_{f_0}\bold H=0.
\end{align*}

Finally, the Jacobi matrix containing the partial derivatives with respect to $m_1,m_2,n,n_2,x_1$ and $x_2$ is
\begin{align*}
\frac{dJ}{d\boldsymbol\theta}(X(0),\boldsymbol\theta)=&\begin{pmatrix}
0&0&0&0&0&0\\
0&0&0&0&0&0\\
0&0&0&0&0&0\\
0&0&0&0&0&0\\
0&0&0&0&0&0\\
j^{00}_{m_1}&j^{00}_{m_2}&j^{00}_{n}&j^{00}_{n_2}&j^{00}_{x_1}&j^{00}_{x_2}\\
0&0&0&0&0&0\\
0&0&0&0&0&0\\
0&0&0&0&0&0\\
\vdots&\vdots&\vdots&\vdots&\vdots&\vdots
\end{pmatrix}
\end{align*}where
\begin{align*}
j^{00}_{m_1}=&\frac{\partial \mathcal L_{f_0}\mathcal L_{f_0}\bold H}{\partial m_1}=I(0)k_IG(0),\\
j^{00}_{m_2}=&\frac{\partial \mathcal L_{f_0}\mathcal L_{f_0}\bold H}{\partial m_2}=-i_1^p(0)k_IG(0),\\
j^{00}_{n}=&\frac{\partial \mathcal L_{f_0}\mathcal L_{f_0}\bold H}{\partial n}=-H(0)k_H\xi(0),\\
j^{00}_{n_2}=&\frac{\partial \mathcal L_{f_0}\mathcal L_{f_0}\bold H}{\partial n_2}=h_1(0)k_H\xi(0),\\
j^{00}_{x_1}=&\frac{\partial \mathcal L_{f_0}\mathcal L_{f_0}\bold H}{\partial x_1}=-H(0)\xi(0)(k_H(H(0)+H_b)),\\
j^{00}_{x_2}=&\frac{\partial \mathcal L_{f_0}\mathcal L_{f_0}\bold H}{\partial x_2}=G(0)I(0)(k_H(H(0)+H_b)).
\end{align*}

If none of the initial conditions ${I(0)},{i_1(0)},{H(0)},{h_1(0)},\xi(0),{G(0)},$ is zero, the matrix $\frac{dJ}{d\boldsymbol\theta}(X(0),\boldsymbol\theta)$ above has a nontrivial null--space since 
\begin{align*}
&\Big(\frac{1}{I(0)},\frac{1}{i_1^p(0)},0,0,0,0\Big)^T\in\mathcal N\Big(\frac{dJ}{d\boldsymbol\theta}(X(0),\boldsymbol\theta)\Big),\\
&\Big(0,0,\frac{1}{H(0)},\frac{1}{h_1(0)},0,0\Big)^T\in\mathcal N\Big(\frac{dJ}{d\boldsymbol\theta}(X(0),\boldsymbol\theta)\Big),\\
&\Big(0,0,0,0,\frac{1}{H(0)\cdot\xi(0)},\frac{1}{G(0)\cdot I(0)}\Big)^T\in\mathcal N\Big(\frac{dJ}{d\boldsymbol\theta}(X(0),\boldsymbol\theta)\Big).
\end{align*}
}

\gdef\thesection{S.7.2}
\section{Results}
 {\color{black}
The results of the SVD of the sensitivity matrix for the bihormonal--glucose model (1)--(7) 
are depicted in Figure \ref{fig:SVD}. 33 vectors of parameters $\boldsymbol\theta_i$ were used, but for simplicity only some cases (that represent the total information obtained) are plotted. See Figure \ref{fig:SVD_full} in  \ref{SMfull} for the singular vectors obtained with each set of parameters.

 In each experimental case, the last 8 singular values are separated from the rest of singular values with a difference of at least 4 decades on the logarithmic scale. Then, 8 singular vectors related to these singular values were plotted in order to visualize the possible parameter correlations. The results are variable from one case to another, but the singular vectors that appears more frequent in different cases are:
$$(m_1,m_2),\quad (n,n_2) \quad\text{and} \quad (x_1,x_2).$$

Parameters $k_1$, $k_I$, $k_{i_1}$, $k_H$ and $m_3$ have nonzero entries in the singular vectors, but their correlations were variable and not repeated in several different cases. For this reason, no further deductions about the correlation of these parameters were done.
}

\gdef\thesection{S.7.3}
\section{Singular Value Decomposition of the Bihormonal--Glucose Model Combining Cases}\label{SMcomb}

{\color{black}
For the bihormonal--glucose model (1)--(7)
, singular value decomposition (SVD) of sensitivity matrices including several experimental case were performed with the same method described in Section \ref{sec:SVD}. Vectors corresponding to the lowest singular values were plotted to reckon the parameters with lack of structural identifiability. The lowest singular values considered were those with value of order $10^{-9}$ or less and separated from the rest of singular values by a gap of order $10^3$ at least. 
}

\gdef\thefigure{S.9}
\begin{figure}[!]\centering
\hspace{-.25cm}\includegraphics[width=10cm]{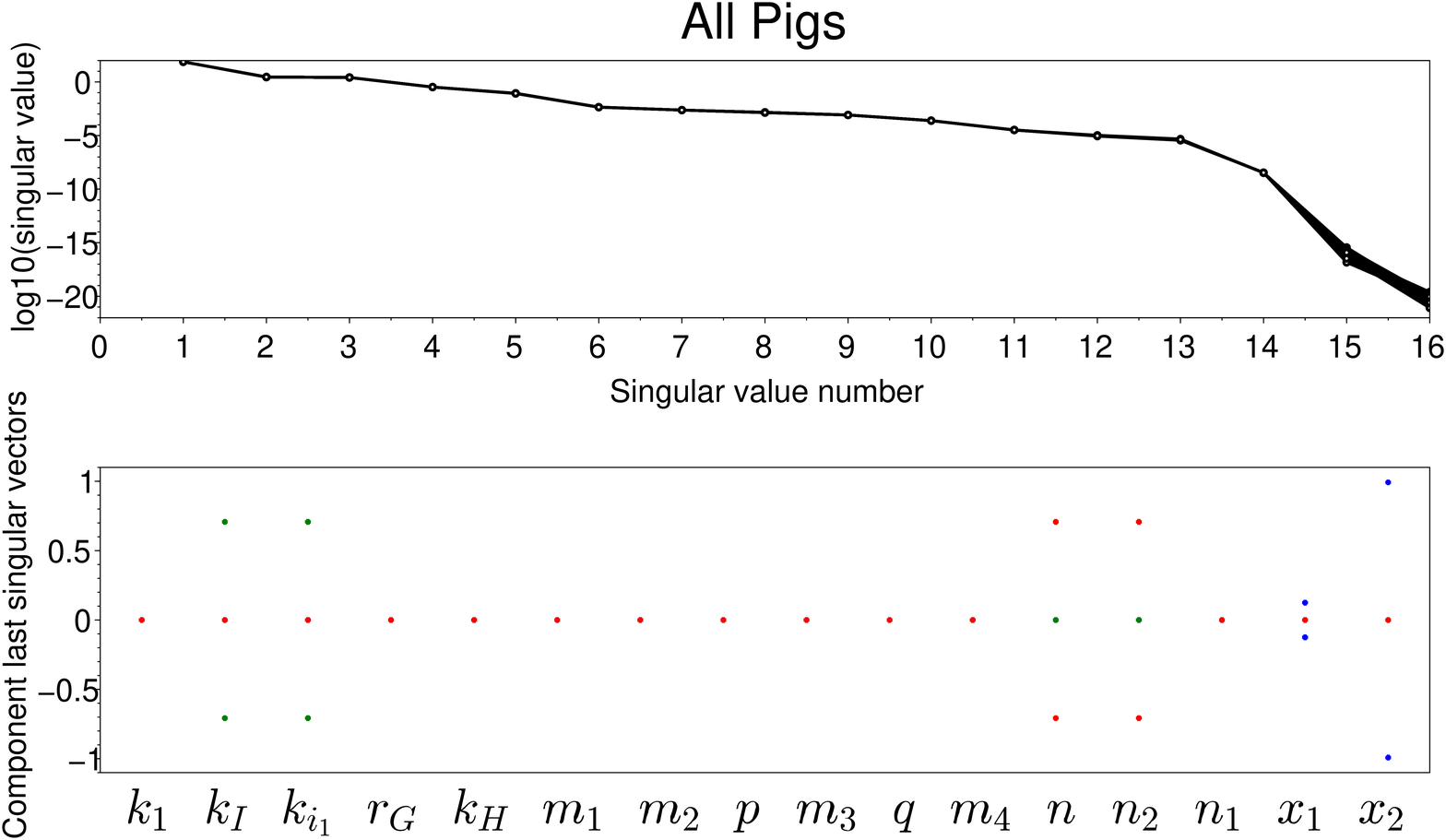}\\

\hspace{-.25cm}\includegraphics[width=10cm]{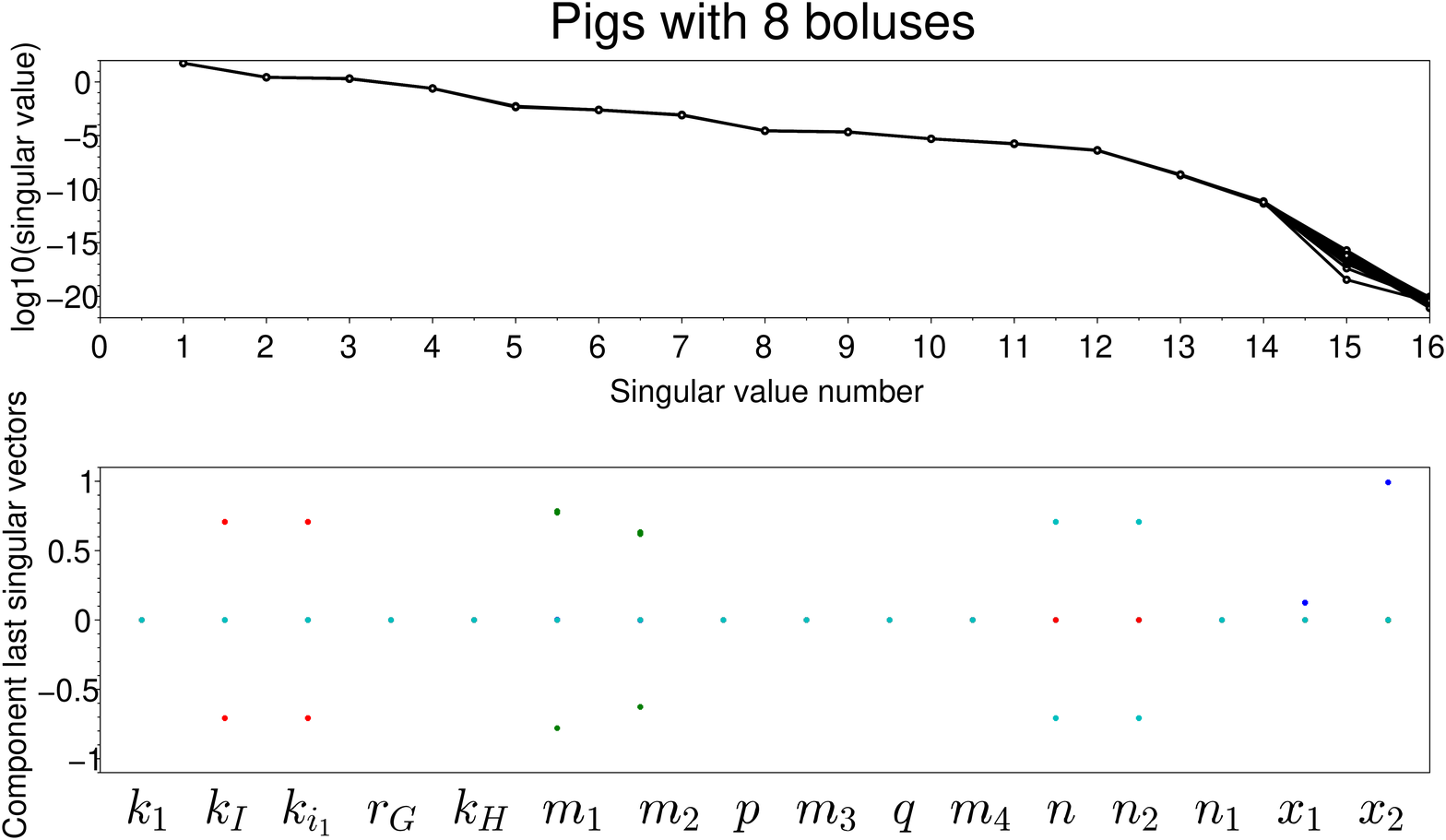}\\

\hspace{-.25cm}\includegraphics[width=10cm]{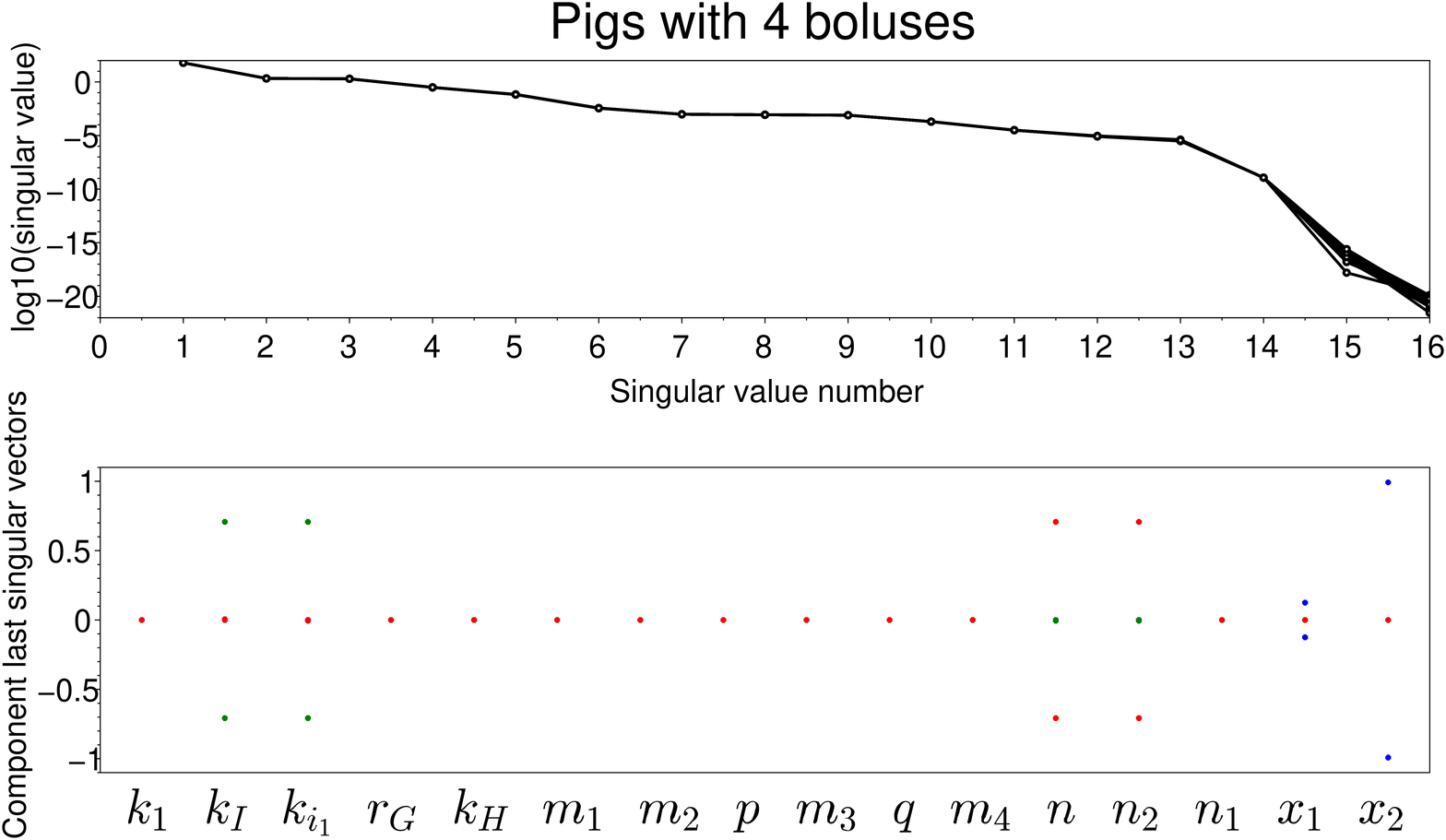}\\
\caption{Singular value decomposition (SVD) of the sensitivity matrix of the bihormonal--glucose model (1)--(7) 
 combining cases. All pigs refers to Pigs 1,2,3,4 and 5. Pigs with 8 boluses are Pigs 1, 2 and 3. Pigs with 4 boluses are Pigs 1,2 and 3 with their first 4 boluses and Pigs 4 and 5 with all their boluses. Each color corresponds to one singular vector. Then, at each case, correlated parameters are those with the same color.}\label{fig:SVD_comb}
\end{figure}

\subsection{Jacobi Rank Test}\label{app_Jacobi_comb}

{\color{black}
The Jacobi matrix containing the partial derivatives with respect to $k_I, k_{i_1}, m_1,m_2,n,n_2,x_1$ and $x_2$ is
\begin{align*}
\frac{d\hat J}{d\boldsymbol\theta}(X(0),\boldsymbol\theta)=&
\begin{pmatrix}
0&0&0&0&0&0&0&0\\
j^{0}_{k_I}&j^{0}_{k_{i_1}}&0&0&0&0&0&0\\
0&0&0&0&0&0&0&0\\
0&0&0&0&0&0&0&0\\
0&0&0&0&0&0&0&0\\
j^{00}_{k_I}&j^{00}_{k_{i_1}}&j^{00}_{m_1}&j^{00}_{m_2}&j^{00}_{n}&j^{00}_{n_2}&j^{00}_{x_1}&j^{00}_{x_2}\\
j^{10}_{k_I}&j^{10}_{k_{i_1}}&0&0&0&0&0&0\\
0&0&0&0&0&0&0&0\\
0&0&0&0&0&0&0&0\\
\vdots&\vdots&\vdots&\vdots&\vdots&\vdots&\vdots&\vdots
\end{pmatrix},
\end{align*}
with
\begin{align*}
&j^{0}_{k_I}=\frac{\partial\mathcal L_{f_0}\bold H}{\partial k_I}=-(I(0)+I_b)G(0),\\
&j^{0}_{k_{i_1}}=\frac{\partial\mathcal L_{f_0}\bold H}{\partial k_{i_1}}=-i_1(0)G(0),\\
&j^{00}_{k_I}=\frac{\partial\mathcal L_{f_0}\mathcal L_{f_0}\bold H}{\partial k_I}=\\
&(-[I(0) + I_b ])G(0)\cdot(-[k_1 + k_I\cdot ( I(0) + I_b ) + k_{i_1}\cdot i_1(0)])\\
&+(-[k_1 + k_I\cdot ( I(0) + I_b ) + k_{i_1}\cdot i_1(0)]\cdot G(0)\\
& + k_H\cdot( H(0) +H_b )\cdot\xi(0))\cdot(-(I(0) + I_b))\\
&+(- m_1\cdot I(0) + m_2\cdot i_1(0)^p)\cdot(-G(0)), 
\end{align*}
\begin{align*}
&j^{00}_{k_{i_1}}=\frac{\partial\mathcal L_{f_0}\mathcal L_{f_0}\bold H}{\partial k_{i_1}}=\\
&[-i_1(0)]\cdot G(0)\cdot(-[k_1 + k_I\cdot ( I(0) + I_b ) + k_{i_1}\cdot i_1(0)])\\
&+(-[k_1 + k_I\cdot ( I(0) + I_b ) + k_{i_1}\cdot i_1(0)]\cdot G(0)\\
& + k_H\cdot( H(0) +H_b )\cdot\xi(0))\cdot(-i_1(0))\\
&+ (- m_3\cdot i_1(0)^q+ m_4\cdot i_2(0))\cdot(-G(0)),\\
&j^{10}_{k_I}=\frac{\partial\mathcal L_{f_1}\mathcal L_{f_0}\bold H}{\partial k_I}=-(I(0)+I_b)r_G,\\
&j^{10}_{k_{i_1}}=\frac{\partial\mathcal L_{f_1}\mathcal L_{f_0}\bold H}{\partial k_{i_1}}=-i_1r_G.
\end{align*}

If none of the initial conditions ${I(0)},{i_1(0)},{H(0)},{h_1(0)},\xi(0),{G(0)},$ is zero, the matrix $\frac{d\hat J}{d\boldsymbol\theta}(X(0),\boldsymbol\theta)$ has a nontrivial null--space since 
\begin{align*}
&\Big(0,0,\frac{1}{I(0)},\frac{1}{i_1^p(0)},0,0,0,0\Big)^T\in\mathcal N\Big(\frac{d\hat J}{d\boldsymbol\theta}(X(0),\boldsymbol\theta)\Big),\\
&\Big(0,0,0,0,\frac{1}{H(0)},\frac{1}{h_1(0)},0,0\Big)^T\in\mathcal N\Big(\frac{d\hat J}{d\boldsymbol\theta}(X(0),\boldsymbol\theta)\Big),\\
&\Big(0,0,0,0,0,0,\frac{1}{H(0)\cdot\xi(0)},\frac{1}{G(0)\cdot I(0)}\Big)^T\in\mathcal N\Big(\frac{d\hat J}{d\boldsymbol\theta}(X(0),\boldsymbol\theta)\Big).
\end{align*}

However, to find a vector with nonzero first two entries (which correspond to $k_I$ and $k_{i_1}$) in the null--space of $\frac{d\hat J}{d\boldsymbol\theta}(X(0),\boldsymbol\theta)$ is possible if an only if
\begin{align*}
0=&(i_1(0) - (I(0)+I_b))(k_1+k_I(I(0)+I_b)+k_{i_1}i_1(0))\\
&-i_1(0)(m_1I(0)-m_2i_1(0)^p)\\
&+(I(0)+I_b)(m_3i_1(0)^q-m_4i_2(0)).
\end{align*}But and initial condition satisfying the equation as above is not a regular point, i.e. the rank of the Jacobi matrix is not constant in a neighborhood of that initial point. Therefore, the rank test cannot be applied as in \cite{stigter2015fast}.}

\gdef\thesection{S.8}
\section{Model Reduction Based on Structural\\ Identifiability Analysis}\label{Section:reduced_system}

\gdef\thesection{S.8.1}
\section{Re--parametrization for Correlated Parameters}
{\color{black}
In case of having non--identifiable correlated parameters, a re--parametrization of the model may be found to combine or fix those parameters 
(see Section \ref{sec:SVD}). For the bihormonal--glucose system (1)--(7)
, the \emph{affine control system} is defined in \ref{app_Jacobi}. Furthermore, the Jacobi matrix is computed and it is found to have a non--trivial null--space, which guarantees that a locally structurally identifiable re--parametrization can be found for parameters $(m_1,m_2), (n,n_2)$ and $(x_1,x_2)$. 

To address the structural non--identifiability and respective correlation of parameters $(m_1,m_2)$, $(n,n_2)$ and $(x_1,x_2)$, the bihormonal--glucose model (1)--(7) 
can be re--parameterized defining 
$$\overline{k_I}=k_I\cdot m_2, \quad\overline{I_b}=I_b/m_2, \quad \overline{x_1}=x_1\cdot n_2,$$ 
$$\overline{k_H}=k_H\cdot x_2\cdot n_2\cdot m_2, \quad \overline{H_b}=H_b\cdot x_2\cdot m_2/n_2,$$
and the state transformations 
$$\overline{I}=I/m_2,\quad\overline H=H/n_2 \quad \text{and} \quad  \overline\xi= \xi/x_2.$$
This leads to have three parameter less in the model while keeping the same output for blood glucose concentration. Indeed, since only blood glucose concentration is compared to experimental data (i.e. $G$ is the only output state), the transformation of any of the other states together with a proper re--parametrization does not affect the system output.

\emph{Note.} For the bihormonal--glucose system (1)--(7)
, the SVD of larger sensitivity matrices that include several experimental cases were computed. The results  are presented in  \ref{SMcomb}. It is observed that parameters
$$(k_I,k_{i_1}),\quad (m_1,m_2),\quad (n,n_2) \quad\text{and} \quad (x_1,x_2)$$are more precisely determined as those corresponding to the lowest singular values (see Figure \ref{fig:SVD_comb}). However, when computing the Jacobi matrix and its null--space, a change of variable cannot be assured for parameters $(k_I, k_{i_1})$ (see  \ref{app_Jacobi_comb}). 

}

\gdef\thesection{S.9}
\section{Model Reduction Based on Practical \\Identifiability and Physical Constraints}\label{Section:reduced_system2}

\gdef\thesection{S.9.1}
\section{Merging Compartments}

{\color{black}
Since the profile likelihood analysis of Section \ref{sec:PL} shows the constant rates of glucose consumption $k_i$ and $k_{i_1}$ lack of practical identifiability, only one of the compartments $I$ and $i_1$ active in glucose removal is kept. 

This reduction is also motivated by the fact of insulin measurements, as well as glucagon measurements, can only be obtained with laboratory methods \cite{staal2019glucose}. Actually, in the first series of animal trials with intraperitoneal insulin boluses \cite{dirnena2019intraperitoneal}, blood samples were analyzed to obtain insulin concentrations, but these did not explain the consumption of glucose. Then, the hypothesis of having an auxiliary state ($i_1$) being active in glucose removal gave better fitting results (compare BIC values of Table 2 in \cite{lopez2019simple} with those in Table \ref{table:model_insulin}).

However, for the Artificial Pancreas blood insulin concentration might not be considered as a physical output state and only one state can be used to represent insulin--dependent consumption of glucose. 

Here, the state $i_1$ is preserved and $I$ is omitted (and $I$ is replaced by $i_1$ In the glucagon sensitivity (7)
), because the mean value of the term $k_{i_1}\cdot i_1$ is slightly larger than that of $k_I(I+I_b)$ (20.00 and 23.89, respectively). Moreover, from the perspective of an intraperitoneal AP, it is expected that not only the insulin in blood will induce a decrease in blood glucose levels, but also that insulin in the intraperitoneal cavity will have a major impact in reducing blood glucose levels.

The single insulin compartment kept will only represent a dimensionless amount of insulin active in glucose consumption and no physiological explanation might be given to it.
}

\gdef\thesection{S.9.2}
\section{No Hormone Measurements}
{\color{black}
 In the first series of animal trials with intraperitoneal insulin boluses \cite{dirnena2019intraperitoneal} and with subcutaneous and intraperitoneal glucagon boluses \cite{am2020intraperitoneal}, blood samples were analyzed to obtain insulin and glucagon concentrations, respectively. Then, the insulin--glucose and glucagon--glucose submodels were used to approximate glucose measurements as well as insulin and glucagon measurements, respectively (see Sections \ref{insulin_model} and \ref{glucagon_model}). 
However, the bihormonal--glucose model (1)--(7) 
is not expected to be calibrated with insulin nor glucagon data, because these hormone measurements can only be obtained with laboratory methods \cite{staal2019glucose} and they will not be used in a AP.

For this reason, insulin and glucagon compartments can be re--scaled without changing the glucose output if a proper re--parametrization is also performed. The objective of re--scaling is to eliminate parameters to be estimated, without affecting the glucose--output (which is the only state to be measured) and achieve locally unique parameter sets. Without re--scaling, the model would be over parameterized.

Blood insulin $I$ and blood glucagon $H$ states were re--scaled in Section \ref{Section:reduced_system}. State $i_1$ can also be re--scaled defining the state transformation
$$\overline{i_1}=i_1/m_4,$$
as well as  new parameters
$$ \overline{k_{i_1}}=k_{i_1}\cdot m_4, \quad \overline{m_3}=m_3\cdot m_4^{q-1}.$$

\emph{Note.} Re--scaling consists of a change of variable, where a state is represented by a multiple of what it might correspond to its actual measurements. For instance, if I is blood insulin concentration, a re--scaled state is $\hat I=a\cdot I$, were $a$ is a positive real number. The derivative of the re--scaled state is obtained using the formula $\hat I'=a\cdot I'$. For the rest of the equations, the term I has to be replaced by $\hat I/a$. 
This change of variable does not numerically affect the equations and the solutions of the rest of states. Moreover, the re--scaled state corresponds to a homothety (i.e. a contraction or a dilatation) of the original state and the original state can be recovered by computing $I=\hat I/a$.  
On the other hand, notice that the re--scaled state has different units that the original one. For example, if $I$ has units mU/L and $a=1000$, then $\hat I$ has units $\mu$U/L or mU/kL.
However, when there are no available measurements of the original state $I$, it is not possible to determine the value of $a$ and the units of $a$ and $\hat I$. For this reason, we consider that the states with no measurements are dimensionless.

\gdef\thesection{S.10}
\section{\color{black}Parameter Estimation of the Reduced\\ 
Model Using the Quasi--Newton Method}\label{newton}

{\color{black}The parameters of the reduced bihormonal--glucose model (8) 
were also estimated using the quasi--Newton method. The \emph{optim} routine in Scilab was used to find a minimum of the cost function 
\begin{align*}
F(\boldsymbol\theta)=&\sum_{t\in T_{BGA}} \big[BGA(t)-G(t,\boldsymbol\theta)\big]^2,
\end{align*}where $\boldsymbol\theta$ is the vector of parameters to be estimated, $BGA(t)$ blood glucose measurements, $T_{BGA}$ the set of time--points at which glucose was measured and $G(t,\boldsymbol\theta)$ the glucose state of the model (8) 
 with the parameters in $\boldsymbol\theta$.

The estimates obtained with the \emph{optim} routine are close to those obtained with the \emph{fminsearch} algorithm in Scilab. The mean of the relative difference (the absolute value of the difference between two values divided by their mean) of the parameters is $0.11$ (compare Table II 
 and Table \ref{table:parameters_reduced_optim}).}

\gdef\thetable{S.6}
\begin{table*}[!]
\begin{center}
\caption{\color{black}Parameters estimated using the quasi--Newton method for the reduced bihormonal--glucose model (8)
, MSE and BIC values. MSE is the Mean Square Error of the model approximation. The better the model performance to approximate the data, the lower the MSE and BIC values.}\label{table:parameters_reduced_optim}
{\color{black}\begin{tabular}{cccccc}
\hline
Parameter			&	Pig 1	&	Pig 2	&	Pig 3	&	Pig 4				&	Pig 5	\\
\hline
$k_1$				&  13.79  &  21.56 &   8.39  &   0.98&   1.41	\\
			
$\overline{k_{i_1}}$	& 171.65  	&  168.68  & 156.51  & 115.36&
   189.84	\\

$\overline{k_H}$		&    38.50	&  44.02& 36.34  & 28.63
 &   657.99   	\\
	
$r_G$				& 4.76	& 2.79 & 2.87   &  1.77  &1.41	\\

$\overline{m_3}$	& 4.80   	&  2.33&   2.46  &   8.74&   37.82	\\

$m_4$				& 27.84	&   85.05 	&47.89 & 9.48 &  12.51	\\

$q$					&  0.48&  0.811&    0.084&   0.96 &1.33 
 	\\
 				
${n}$				&	   110.34   &  177.30 &  163.43  &  37.85&
   709.38	\\
$n_1$				&	138.95  &  200.09& 214.47  &  38.51&
   2152.04	\\

$\overline{x_1}$ 		& 237.39	&
   196.85& 94.27 &  0.0039&   3393.97	\\
\hline
MSE	&	0.297&0.52&0.25& 3.40&0.19\\				
\hline
BIC					&-57.66 &-12.03&-64.80&173.83&	-111.71	\\
\hline
\end{tabular}}
\end{center}
\end{table*}


\end{document}